 \documentclass[11pt,reqno]{amsart}
 \usepackage{amsmath,amssymb,amsfonts,mathrsfs}
\usepackage{url}
\usepackage[letterpaper,margin=1in]{geometry}
\usepackage[expansion=false]{microtype}
\numberwithin{equation}{section}

 \theoremstyle{plain}
\newtheorem{theorem}{Theorem}[section]
\newtheorem{lemma}[theorem]{Lemma}
\newtheorem{prop}[theorem]{Proposition}
\newtheorem{corollary}[theorem]{Corollary}
\theoremstyle{definition}
\newtheorem{definition}[theorem]{Definition}
\newtheorem{example}[theorem]{Example}

\theoremstyle{remark}
\newtheorem{remark}[theorem]{Remark}

\def\si{{\sigma}}

\def\EE{\mathbb E}\def\PP{\mathbb P}\def\QQ{\mathbb Q}\def\RR{\mathbb R}

\def\wt{\widetilde}
\def\wh{\widehat}
\def\sC{{\mathcal C}}
\def\sE{{\mathcal E}}
\def\sF{{\mathcal F}}
\def\sG{{\mathcal G}}
\def\sH{{\mathcal H}}
\def\sL{{\mathcal L}}

\def\sN{{\mathcal N}}
\def\sS{{\mathcal S}}
\def\sP{ {\mathcal P}_{\mathcal M}}
\def\wsP{ \wt {\mathcal P}_{\mathcal M}}
\IfFileExists{bbm.sty}{\usepackage{bbm}\def\1{{\mathbbm 1}}}{\def\1{{\mathbf 1}}}
\def\eps{\varepsilon}
\def\dKL{D_{\mathrm{KL}}}

\def\<{\langle}
\def\>{\rangle}

\usepackage[hidelinks,pdftitle={Sub-Diffusive Black--Scholes Model and Structure-Preserving Girsanov Transform},pdfauthor={Shuaiqi Zhang, Zhen-Qing Chen}]{hyperref}

\makeatletter
\renewcommand{\l@section}{\@tocline{1}{0pt}{1pc}{1.5pc}{}}
\renewcommand{\l@subsection}{\@tocline{2}{0pt}{3pc}{4pc}{}}
\makeatother
\begin{document}

\makeatletter
\def\@settitle{\begin{center}\baselineskip14\p@\relax\bfseries\@title\par\end{center}}
\def\author@andify{\nxandlist{\unskip,\space\ignorespaces}{\unskip\hspace{1.67em}and\hspace{1.67em}}{\unskip,\space and\space}}
\def\@setauthors{\begingroup\def\thanks{\protect\thanks@warning}\trivlist\centering\normalsize\@topsep30\p@\relax\advance\@topsep by -\baselineskip\item\relax\author@andify\authors\def\\{\protect\linebreak}\authors\endtrivlist\endgroup}
\let\old@setaddresses\@setaddresses
\def\@setaddresses{\begingroup\let\scshape\relax\old@setaddresses\endgroup}
\renewcommand{\@bibtitlestyle}{{\let\@tocwriteb\@gobblethree\@xp\section\@xp*\@xp{\refname}}}
\makeatother
\title[Sub-Diffusive Black--Scholes Model]{Sub-Diffusive Black--Scholes Model and Structure-Preserving Girsanov Transform}

\author{Shuaiqi Zhang}
\address{School of Mathematics, China University of Mining and Technology, Xuzhou, Jiangsu, 221116, China.}
\email{shuaiqiz@hotmail.com}
\author{Zhen-Qing Chen}
\address{Department of Mathematics, University of Washington, Seattle, WA 98195, USA.}
\email{zqchen@uw.edu}
\thanks{SZ's research is partially supported by the Humanities and Social Sciences Foundation of the Ministry of Education of China (Application of Anomalous Sub-diffusion in Insurance and Finance, Grant No.~24YJA910008) and by the National Natural Science Foundation of China (Grant No.~12571519). ZC's research is partially supported by a Simons Foundation Grant. Corresponding author: Zhen-Qing Chen.}
\date{November 8, 2026}

\begin{abstract}
We propose a sub-diffusive Black--Scholes model in which stock prices are governed by stochastic differential equations driven by sub-diffusions.
The model captures the periods of low trading activity observed in bear markets, and it contains the classical Black--Scholes model as a special case.
The sub-diffusive spot market is arbitrage-free but, in general, incomplete.
We study the pricing of European contingent claims in this model.
To this end, we establish a structure-preserving Girsanov transform for sub-diffusions.
We use it to identify the minimal martingale measure of the model.
This measure is the unique equivalent martingale measure under which the discounted stock prices are martingales driven by sub-diffusions with the same law as the original ones.
Using this measure, we price every contingent claim in $L^p$ with $p>2$.
For each such claim, we construct a replicating strategy that is mean self-financing and minimizes the local risk, and whose value process gives the price process of the claim.
Finally, we derive an explicit formula for the price of a European call option.
We show that this price solves a partial differential equation with a fractional derivative in time, which we call the time-fractional Black--Scholes equation.


\end{abstract}

\maketitle
\markboth{Shuaiqi Zhang and Zhen-Qing Chen}{Sub-Diffusive Black--Scholes Model}

\begingroup\normalsize\linespread{1.0}\selectfont\setlength{\parskip}{0pt}\vspace{-7ex}
\tableofcontents
\endgroup

\section{Introduction}\label{S1}

A call option is a contract between a buyer and a seller.
It gives the buyer the right, but not the obligation, to buy an underlying asset at a predetermined price, called the strike price, on or before a specified date, called the expiration date.
A European call option can be exercised only on the expiration date.
The pricing of European options is a central problem in mathematical finance.

The Black--Scholes model of Black and Scholes \cite{1973} and Merton \cite{Merton1973} provides a foundational framework for pricing European options under the assumption of no arbitrage.
Its formulation in terms of stochastic analysis, developed by Harrison and Kreps \cite{HK} and Harrison and Pliska \cite{HP}, is based on equivalent martingale measures and self-financing replicating strategies.
In the classical Black--Scholes model, the prices of the underlying stocks are governed by stochastic differential equations driven by Brownian motion.
The model is widely used by practitioners in options markets, often with some modifications.
However, its spot-market assumption does not account for periods of low trading activity, such as those observed in bear markets; see \cite[Section~4]{JW}.
Regime-dependent appreciation rates \cite{zhang2018} can describe changes in returns, but by themselves they do not capture intermittent market activity.

Time changes by inverse subordinators provide a mechanism for modeling market stagnation.
Magdziarz \cite{Mag} studied European options for a single stock whose price is a geometric Brownian motion time-changed by an inverse stable subordinator.
Karipova and Magdziarz \cite{KM} extended this approach to several stocks.
Option pricing under time-changed fractional Brownian motion was studied in \cite{Sh}.
Motivated by these works, we formulate a sub-diffusive market model with adapted, time-dependent interest rates, appreciation rates and volatilities, driven by sub-diffusions with a general inverse subordinator as the time change.
The model contains the classical Black--Scholes model as a special case and is considerably more flexible.
To the best of our knowledge, the following results have not been investigated before for a sub-diffusive Black--Scholes model.
\begin{enumerate}
\item[(i)] We characterize all equivalent martingale measures and identify the minimal martingale measure.
\item[(ii)] We construct a replicating strategy that is mean self-financing and minimizes the local risk.
\item[(iii)] We justify the pricing rule: the price of a contingent claim is the initial value of this strategy.
\item[(iv)] We derive the time-fractional Black--Scholes equation for the price of a European call option when the coefficients $r_t$, $\bar\mu_t$ and $\sigma_t$ of the stock price process are constant and the inverse subordinator is general.
\end{enumerate}

Sub-diffusions arise naturally from particle trapping and waiting-time effects.
They have applications to transport in porous media, systems with memory, and other anomalous transport phenomena \cite{MS04,KlS,MK}.
We first recall their structure and some quantitative features, and then we introduce the financial model.
For simplicity, we describe the model for a single stock in this introduction.
All results of the paper are established for a market consisting of one bond and several stocks; see Section~\ref{S:4}.

Let $S=\{S_t;\, t\geq 0\}$ be a subordinator with drift $\kappa\geq 0$ and L\'evy measure $\nu$.
That is, $S_t=\kappa t+S_t^0$, where $S^0$ is a driftless subordinator with L\'evy measure $\nu$.
We assume that $S$ is strictly increasing $\PP$-almost surely.
This holds if and only if $\kappa>0$ or $\nu((0,\infty))=\infty$.
Let $L_t:=\inf\{r>0: S_r>t\}$, $t\geq 0$, be the right-continuous inverse of $S$.
Since $S$ is strictly increasing, the paths of $L$ are continuous.
If $\nu$ is nonzero, then $L$ is constant on infinitely many disjoint time intervals, each corresponding to a jump of $S$.
If $\kappa=0$, then $\PP$-almost surely the measure $dL_t$ is singular with respect to the Lebesgue measure $dt$.
Let $B$ be a $d$-dimensional Brownian motion independent of $S$.
For $x\in\RR^d$, we write $\PP_x$ for a probability measure under which $B_0=x$, and $\EE_x$ for the expectation under $\PP_x$.
The time-changed process $B_{L_t}$ is a prototype of sub-diffusions.
It is a continuous martingale, and each coordinate $B^i_L$ has quadratic variation $\langle B^i_L\rangle_t=L_t$.
The class of subordinators considered here includes the deterministic subordinator $S_t\equiv t$, for which $L_t\equiv t$ and $B_{L_t}=B_t$.
Hence our results cover the Brownian case.

If $\nu$ is nonzero, we call $B_{L_t}$ an anomalous sub-diffusion.
Unlike Brownian motion, an anomalous sub-diffusion is not a Markov process.
Suppose that $S$ is a $\beta$-stable subordinator with $0<\beta<1$, that is,  it is a non-decreasing L\'evy process whose Laplace exponent is $\phi(\lambda)=\lambda^\beta$.
By the independence of $B$ and $S$ and by \eqref{e:2.8} below, the mean square displacement of $B_{L_t}$ is
$$
\EE_x\bigl[|B_{L_t}-B_0|^2\bigr]=d\,\EE[L_t]=\frac{d\,t^{\beta}}{\beta\,\Gamma(\beta)}
\quad\text{for every } x\in\RR^d \text{ and } t\geq 0,
$$
where $\Gamma(z):=\int_0^\infty s^{z-1}e^{-s}\,ds$, $z>0$, is the Gamma function.
In particular, $\lim_{t\to\infty} t^{-1}\,\EE_x\bigl[|B_{L_t}-B_0|^2\bigr]=0$ for every $x\in\RR^d$, so the mean square displacement grows sublinearly.

To compare $B_{L_t}$ with $B$, we recall density estimates from \cite{CKKW1}.
For positive functions $f$ and $g$, we write $f\asymp g$ if there is a constant $c\geq 1$ such that $c^{-1}g\leq f\leq c\,g$.
By \cite[Theorem~1.3]{CKKW1}, if $S$ is a $\beta$-stable subordinator, then for every $t>0$ the random variable $B_{L_t}$ under $\PP_0$ has a density $q(t,y)$ with respect to the Lebesgue measure $dy$ on $\RR^d$.
Moreover, there are constants $C_1,\dots,C_4>0$, depending only on $d$ and $\beta$, such that the following estimates hold.
\begin{enumerate}
\item[(i)] If $|y|^2\leq t^{\beta}$, with $y\neq 0$ when $d\geq 2$, then
$$
q(t,y)\asymp
\begin{cases}
t^{-\beta/2} & \text{for } d=1,\\
t^{-\beta}\log\bigl(2t^{\beta}/|y|^2\bigr) & \text{for } d=2,\\
t^{-\beta}|y|^{-(d-2)} & \text{for } d\geq 3.
\end{cases}
$$
\item[(ii)] If $|y|^2>t^{\beta}$, then
$$
C_1 t^{-\beta d/2}\exp\Bigl(-C_2\bigl(|y|^2/t^{\beta}\bigr)^{1/(2-\beta)}\Bigr)
\leq q(t,y)\leq
C_3 t^{-\beta d/2}\exp\Bigl(-C_4\bigl(|y|^2/t^{\beta}\bigr)^{1/(2-\beta)}\Bigr).  
$$
\end{enumerate}
Since $2/(2-\beta)<2$, the tails of $q(t,\cdot)$ decay more slowly than Gaussian tails.
This reflects the anomalous spreading of $B_{L_t}$ relative to $B$.

\medskip

In the classical Black--Scholes model, the price process $\sS=\{\sS_t;\, t\geq 0\}$ of a stock satisfies
\begin{equation}\label{e:1.2}
d\sS_t=\sS_t\left(\mu_t\,dt+\sigma_t\,dB_t\right),
\end{equation}
where $\mu_t$ is the appreciation rate of the stock, $\sigma_t>0$ is its volatility, and $B$ is a one-dimensional standard Brownian motion.
During a bear market, a stock typically has a smaller or even negative appreciation rate, and it is also traded less actively.
It is therefore natural to model the stock price by a stochastic differential equation driven by a sub-diffusion instead of a Brownian motion.
We propose the dynamics
\begin{equation}\label{BSsub}
d\sS_t=\sS_t\left(r_t\,dt+\bar\mu_t\,dL_{(t-a)^+}+\sigma_t\,dB_{L_{(t-a)^+}}\right),
\end{equation}
where $(t-a)^+:=\max\{t-a,0\}$, the constant $a\geq 0$ is the time at which the market becomes active, and $r_t\geq 0$, $\bar\mu_t$ and $\sigma_t>0$ are adapted processes.
In this model, the stock has a base appreciation rate equal to the interest rate $r_t$.
In addition, it earns an excess return $\bar\mu_t\,dL_{(t-a)^+}$, which accrues only when the market is active.
If $a=0$ and $S_t\equiv t$, that is, $\kappa=1$ and $\nu=0$, then $L_t\equiv t$ and \eqref{BSsub} reduces to \eqref{e:1.2} with $\mu_t=r_t+\bar\mu_t$.

The model \eqref{BSsub} is designed to capture the prolonged periods of subdued activity that are often observed in bear markets.
A striking example is the S\&P~500 index during the bear market of 2000--2002 that followed the dot-com bubble.
The index peaked at about 1{,}528 in March 2000 and fell by about 49\% to roughly 777 in October 2002.\footnote{S\&P Dow Jones Indices, S\&P~500 historical index levels, \url{https://www.spglobal.com/spdji/en/indices/equity/sp-500/}.}
It did not recover quickly.
Instead, it moved sideways or downward for several years and did not regain its 2000 peak until 2007.
Such stagnation is described more naturally by \eqref{BSsub} with time-dependent coefficients than by \eqref{e:1.2}.
Indeed, in \eqref{BSsub}, long stretches on which $B_{L_t}$ is constant can alternate with periods of negative drift.

We point out that the dynamics in \eqref{BSsub} are, in general, not obtained by a time change of the classical model \eqref{e:1.2}.
The reason is that the $dt$ and $dL_{(t-a)^+}$ terms evolve on different time scales.
Recall that if $\kappa=0$, then $\PP$-almost surely $dL_t$ is singular with respect to $dt$.
Our objective is to price and hedge European contingent claims in this market, which is incomplete in general.
For option pricing in other incomplete markets, we refer the reader to \cite{Bj, FS, GJ, JW, Sch, Sh} and the references therein.

Since the sub-diffusive market is incomplete, a central issue for pricing and hedging is the choice of a martingale measure.
Theorems~\ref{T:EMM} and~\ref{T:REMM} characterize the equivalent martingale measures and the generalized equivalent martingale measures, respectively.
If the subordinator is not deterministic, each of these classes contains infinitely many measures.
Nevertheless, Theorem~\ref{T:4.1} identifies a unique measure $\QQ^*$ that belongs to both classes and under which the discounted stock prices are martingales driven by a sub-diffusion $\widetilde B_{L_{(t-a)^+}}$ whose law under $\QQ^*$ equals the law of $B_{L_{(t-a)^+}}$ under $\PP$.
We call $\QQ^*$ the \emph{sub-diffusion equivalent martingale measure}.
Its construction relies on a structure-preserving Girsanov theorem (Theorem~\ref{G1}), which retains the representation by a Brownian motion and an independent inverse subordinator.
The general Girsanov theorem for continuous local martingales does not guarantee this property.
Moreover, $\QQ^*$ is the minimal martingale measure in the sense of F\"ollmer and Schweizer \cite{FS} (Theorem~\ref{T:MMM}), and it is the unique equivalent martingale measure that minimizes the reverse relative entropy $\dKL(\PP\Vert\QQ)$ (Theorem~\ref{T:4.6}).
Consequently, for every contingent claim in $L^p$ with $p>2$, $\QQ^*$ yields a replicating strategy that is mean self-financing and minimizes the local risk. The martingale representation theorem established in \cite{ZC2} for sub-diffusions plays a crucial role in our study. 

\medskip

The main contributions of this paper are as follows.
\begin{enumerate}
\item[(i)] We introduce the spot-market dynamics \eqref{BSsub}, which involve two time scales, a general strictly increasing subordinator, and adapted coefficients.
The multidimensional market is developed in Section~\ref{S:4}.

\item[(ii)] Theorems~\ref{G1} and~\ref{G2} establish changes of measure that add or remove drifts of the form $\int_0^t H_s\,dL_{(s-a)^+}$ and retain the sub-diffusive structure, including the independence of the Brownian motion and the subordinator.
This holds for every strictly increasing subordinator $S$.
Theorems~\ref{T:3.5} and~\ref{T:uniqueness} give the corresponding uniqueness results.
For $\kappa=1$ and $\nu=0$, we recover the Girsanov theorem for Brownian motion.

\item[(iii)] Using Theorems~\ref{G1}, \ref{T:3.5} and~\ref{T:uniqueness}, we show in Theorem~\ref{T:4.1} that $\QQ^*$ is the unique equivalent martingale measure under which the discounted stock prices are martingales driven by sub-diffusions with the same law as the original ones.
We further show that $\QQ^*$ is the minimal martingale measure (Theorem~\ref{T:MMM}) and the unique minimizer of the reverse relative entropy $\dKL(\PP\Vert\QQ)$ among all equivalent martingale measures (Theorem~\ref{T:4.6}).
Hence the sub-diffusive Black--Scholes model fits into the framework of F\"ollmer and Schweizer \cite{FS} for incomplete markets.
Using $\QQ^*$, we price contingent claims and construct replicating strategies that are mean self-financing and minimize the local risk (Proposition~\ref{P:5.5} and \S\ref{S:5.2}).

\item[(iv)] Suppose that additional information is available, described by an enlarged filtration $\{\widetilde\sF^a_t\}_{t\in[0,T]}$, where $T>0$ is the maturity.
Then every contingent claim in $L^p$ for some $p>2$ can be replicated by a self-financing strategy, although the initial endowment is random.
In Theorem~\ref{T:5.7}, we relate the resulting price to the locally risk-minimizing strategy under the original filtration.

\item[(v)] For constant $r_t$, $\bar\mu_t$ and $\sigma_t$, we derive an explicit formula for the price of a European call option (Theorem~\ref{T:5.1}, \eqref{e:5.18n} and \eqref{e:5.45}), which is the sub-diffusive counterpart of the classical Black--Scholes formula.
In Theorem~\ref{T:5.6}, we show that this price solves a partial differential equation with a fractional derivative in time, which we call the time-fractional Black--Scholes equation.
\end{enumerate}

\medskip

The rest of this paper is organized as follows.
Section~\ref{S:2} contains preliminaries on inverse subordinators, anomalous sub-diffusions and the filtrations involved.
Section~\ref{S:3} is devoted to the   Girsanov transform for sub-diffusions.
In Section~\ref{S:4}, we describe the spot market consisting of one bond and $d$ stocks, and we establish the existence and uniqueness of the sub-diffusion equivalent martingale measure $\QQ^*$, also called the risk-neutral measure, together with its minimality properties.
In Section~\ref{S:5}, we use $\QQ^*$ to price and hedge European options, and we derive the explicit call price formula and the time-fractional Black--Scholes equation.

Throughout this paper, the symbol $:=$ denotes a definition.
For $u,v\in\RR$, we write $u\vee v:=\max\{u,v\}$, $u\wedge v:=\min\{u,v\}$ and $u^+:=u\vee 0$.
The abbreviation rcll stands for ``right-continuous with left limits''.

 \section{Preliminaries}\label{S:2}
 
 A subordinator $S=\{S_t; t\geq 0\}$ is a non-negative real-valued L\'evy process with $S_0=0$.
 It is characterized by its  Laplace exponent $\phi$, a non-negative function defined on $[0, \infty)$ so that   
   \begin{equation}\label{e:2.2}
   \EE \left[ e^{-\lambda S_t} \right] = e^{-t \phi (\lambda)}
     \quad \hbox{for every } \lambda >0 \hbox{ and }   t\ge0.
 \end{equation}
 By the L\'evy-Khintchine formula, a non-negative function $\phi$ defined on $[0, \infty)$  is the Laplace exponent of
 a subordinator if and only if 
  there is a constant $\kappa\geq 0$ and a $\sigma$-finite measure  $\nu$  on $(0, \infty)$ satisfying
 $\int_0^\infty (1\wedge x) \nu (dx)<\infty$ so that   
  \begin{equation}\label{e:2.2a}
 \phi (\lambda) =\kappa \lambda + \int_0^\infty (1-e^{-\lambda x}) \nu ( dx).
 \end{equation}
 The constant $\kappa$ is called the drift and $\nu$ is called the L\'evy measure of the subordinator. 
  They are in one-to-one correspondence with the Laplace exponent of the subordinator. 
 So a subordinator $S$ with drift $\kappa \geq 0$ and  a 
L\'evy measure $\nu$ can be expressed as
 $$ 
 S_t=\kappa t+ S_t^0, \quad t\geq 0,
 $$
 where $S_t^0$ is a driftless subordinator with L\'evy measure $\nu$. 
 
 \subsection{Inverse Subordinator}
 
 The inverse $L=\{L_t; t\geq 0\}$ of the subordinator $S$ is defined to be 
 $$
 L_t = \inf\{r\geq 0:   S _r  >t      \}  \quad \hbox{for } t\geq 0.
 $$
 
 The following exponential integrability result  was first proved in \cite[Lemma 8]{JK}.
 Here we give an alternative proof.

\begin{theorem}\label{T:2.1} 
Let $S$ be a subordinator that is not identically zero.  For every $t>0$, $p>0$   and $0<q<\phi (\infty) = \kappa + \nu (0, \infty)$, 
$\EE [ L_t^p]<\infty$ and $\EE \left[ e^{q L_t} \right] <\infty$.
 Moreover,  their Laplace transforms are given by
 \begin{eqnarray}
 \int_0^\infty e^{-\lambda t} \EE \left[ L_t^p \right] dt &=& \frac{\Gamma (p+1)}{\lambda \phi (\lambda)^p} \qquad \hbox{for every } \lambda >0,
 \label{e:Ep} \\
  \int_0^\infty e^{-\lambda t} \EE \left[ e^{\gamma L_t}\right]  dt &=&\frac{\phi (\lambda)}{\lambda (\phi (\lambda )-\gamma)}
 \qquad  \hbox{for every } \gamma  < \phi (\infty)  \hbox{ and } \lambda >\phi^{-1}(\gamma\vee 0) \label{e:eEq} .
  \end{eqnarray} 
\end{theorem}

\begin{proof} Let $S$ be a subordinator with Laplace exponent $\phi$ and let $L$ be its inverse. First observe that, for any $r>0$ and $\lambda >0$, integration by parts gives 
\begin{eqnarray}
\int_0^\infty  e^{-\lambda t} \PP (S_r\leq t)  dt 
&=&  -\frac{1}\lambda \int_0^\infty   \PP (S_r\leq t)  d_t ( e^{-\lambda t})
=  \frac1{\lambda} \int_0^\infty     e^{-\lambda t}  d_t  \PP (S_r\leq t)   \nonumber \\
&=&\frac1{\lambda} \EE \left[ e^{-\lambda S_r}\right]= \frac{1} {\lambda} e^{-r \phi(\lambda)} . 
\label{e:2.7a}
\end{eqnarray} 
As for any $p>0$ and $t>0$, 
$$
\EE \left[ L_t^p\right] = \EE \int_0^\infty p r^{p-1} \1_{\{ S_r\leq t\}} dr
=p \int_0^\infty r^{p-1} \PP (S_r\leq t) dr ,
$$
we have by the Tonelli theorem and \eqref{e:2.7a} that for every $\lambda >0$,   
\begin{eqnarray*}
\int_0^\infty e^{-\lambda t} \EE \left[ L_t^p \right] dt 
&=& p \int_0^\infty  r^{p-1} \left( \int_0^\infty  e^{-\lambda t} \PP (S_r\leq t)  dt \right) dr \\
 &=& \frac{p}{\lambda} \int_0^\infty  r^{p-1}  e^{- r\phi (\lambda)}  dr \\
&=& \frac{p \Gamma (p) }{\lambda \phi (\lambda)^p} = \frac{  \Gamma (p+1)}{\lambda \phi (\lambda)^p} ,
\end{eqnarray*} 
establishing \eqref{e:Ep}. Similarly,  as for any $\gamma\in \RR $ and $t>0$, 
 \begin{equation}
\EE \left[ e^{\gamma L_t} \right] = 1+ \EE \int_0^\infty \gamma e^{\gamma r} \1_{\{ S_r\leq t\}} dr 
=1+  \gamma \, \int_0^\infty   e^{\gamma r} \, \PP (  S_r\leq t ) dr.
\end{equation}
We have by the Tonelli theorem and \eqref{e:2.7a}  that for any   $\gamma < \phi (\infty)$ and $\lambda >\phi^{-1}(\gamma \vee 0)$, 
 \begin{eqnarray*}
\int_0^\infty e^{-\lambda t} \EE \left[ e^{\gamma L_t}  \right] dt &=&
\frac1{\lambda}  +   \gamma \int_0^\infty  e^{\gamma r} \left( \int_0^\infty  e^{-\lambda t} \PP (S_r\leq t)  dt \right) dr \\
  &=& \frac1\lambda+ \frac{  \gamma} {\lambda} \int_0^\infty   e^{-(\phi (\lambda) -\gamma)  r } dr \\
&=& \frac{\phi (\lambda) }{\lambda (\phi (\lambda)-\gamma)},
\end{eqnarray*} 
which is \eqref{e:eEq}. 
From  \eqref{e:Ep}-\eqref{e:eEq},  we have for every $p>0$   and $0<\gamma < \phi (\infty)$,
$\EE [ L_t^p]<\infty$ and $\EE \left[ e^{\gamma  L_t} \right] <\infty$ for a.e. $t>0$
and hence for every $t>0$ by the increasing monotonicity of $t\mapsto L_t$.
\end{proof} 

 \medskip
 
 \begin{remark}\rm 
 \begin{enumerate}
 \item[(i)]
When $S$ is a $\beta$-subordinator with $\beta \in (0, 1)$, that is, when $\phi (r)=r^\beta$,
we have by \eqref{e:Ep} that for any $p>0$, 
\begin{equation} \label{e:2.7}
 \int_0^\infty e^{-\lambda t} \EE \left[ L_t^p \right] dt  =\frac{\Gamma (p+1)}{\lambda^{\beta p +1}  } \quad \hbox{for every } \lambda >0 . 
\end{equation} 
Taking the inverse Laplace transform yields 
\begin{equation} \label{e:2.8} 
   \EE \left[ L_t^p \right]   =\frac{\Gamma (p+1)}{\Gamma (\beta p +1)  } t^{\beta p}  \quad \hbox{for every } t \geq 0.
\end{equation} 
This recovers the moment formula of the inverse $\beta$-stable subordinator obtained in \cite[(2.4) and (2.7)]{PSW2005}.    

\item[(ii)] With the exponential integrability of the inverse subordinator from Theorem \ref{T:2.1},  
we can readily drop  the uniform boundedness condition  in 
  \cite{C1, C2, CKKW1, CKKW2} imposed  on the  strongly continuous semigroup $\{T_t; t\geq 0\}$ or the Markov process $X$ whose transition semigroup is strongly continuous in some Banach space $({\mathbb B}, \| \cdot \|_{\mathbb B})$.  \qed
\end{enumerate}
 \end{remark}

\subsection{Sub-diffusions and filtrations} \label{S:2.2}
 
 In the rest of this paper,  we  assume the subordinator $S$ is strictly increasing; that is, either $S$ has a positive drift $\kappa >0$
 or it has infinite L\'evy measure $\nu$ or both.  Under this assumption, the  inverse  subordinator $L$ is continuous.
 Let $B$ be a standard Brownian motion on $\RR^d $ starting from the origin $\mathbf 0$ that  is  independent of $S$.  Let $(\Omega,\mathcal F,\mathbb P)$ be the canonical probability space of $(B,S)$.
That is, $\Omega=C([0,\infty);\RR^d)\times \mathbb D([0,\infty);[0,\infty))$, where $\mathbb D$ is the Skorokhod space of rcll paths, and $\mathbb P$ is the product of the Wiener measure and the law of $S$.
Apart from $\nu$ and the potential measure $U$ of Section~\ref{S:5.6}, which are $\sigma$-finite, all measures below are probability measures.
For a probability measure $\mathbb Q$ on a $\sigma$-field $\mathcal G$ and a $\sigma$-field $\mathcal G'\subset\mathcal G$, we write $\mathbb Q|_{\mathcal G'}$ for the restriction of $\mathbb Q$ to $\mathcal G'$.  For $a\geq 0$,
 the sub-diffusion $t\mapsto B_{L_{(t-a)^+}}$ itself is not a Markov process. However if we add an overshoot process $R_t$, it will make the pair a 
 time-homogeneous Markov process. 
 The following result is established in \cite[Theorem 3.1]{ZC}, 
 which holds in multidimensional case by the same proof  although it  was stated there for one-dimension. 
    
 \begin{theorem}[Theorem 3.1 of  \cite{ZC}] \label{T:2.3} 
  Suppose that $B$  is a  standard Brownian motion on $\RR^d $ starting from $\mathbf 0$
  and  $S$ is any subordinator that is strictly increasing and independent of $B$ with $S_0=0$.
  For each $t\geq 0$, define  $L_t:=\inf\{r>0:S_r>t\}$.
Then
\begin{equation} \label{e:2.1}
\wt X_t:=(X_t, \, R_t):= \left(X_0+ B_{L_{(t-R_0)^+}}, 
 \, R_0+ S_{L_{(t-R_0)^+}} - t  \right), \quad t\geq 0,
\end{equation}
with $\wt X_0= (X_0, R_0) \in \RR^{d} \times [0, \infty)$   is a time-homogeneous Markov process taking values in the upper half-space $\mathcal X:=\RR^{d}\times [0, \infty)$ of $\RR^{d+1}$.
\end{theorem}

The overshoot process $R_t:=R_0+ S_{L_{(t-R_0)^+}} - t $  can be interpreted as the wake-up or holding time for the sub-diffusion process $X_t$ at time $t$. 
 
\medskip

  Denote by $\{\sG^B_t\}$ the natural augmented filtration generated by the Brownian motion $B$, 
that is, $\sG^B_t :=   \sG_{t+}^{B, 0} \vee {\mathcal N}$, where $\sG^{B, 0}_t:=\sigma (B_s; s\in [0, t]) $,
$ \sG_{t+}^{B, 0}  := \cap_{s>t} \sG^{B, 0}_s$, 
 and $\sN:=\{A\in \sG^{B, 0}_\infty: \PP (A)=0\}$. The notation $ \sG_{t+}^{B, 0} \vee {\mathcal N}$ 
  stands for the  $\sigma$-field generated by $ \sG_{t+}^{B,0} \cup  {\mathcal N}$. 
   The natural augmented filtration $\{\sG_t^B\}$ is right-continuous. 
   Similarly,   we denote by   $\{\sG^L_t\}$ the natural augmented filtration generated by 
   the inverse subordinator $L$.  
  Let  $\{\sF'_t\}_{t\geq 0}$ be the natural augmented filtration generated by $X$ and 
     $\{\sF_t\}_{t\geq 0}$ be     the minimum augmented filtration generated by 
  $\wt X=(X, R)$.   Clearly, $\sF'_t\subset \sF_t$ for every $t\geq 0$. 
Note that process $\wt X_t=(X_t, R_t)$  
  depends on the initial $a\geq 0 $ of $R_0$, so do the filtrations $\{\sF_t\}_{t\geq 0}$ and $\{\sF'_t\}_{t\geq 0}$. When $R_0=a$ for some deterministic $a\geq 0$,  sometimes for emphasis we denote them
by $\wt X^a_t=(X^a_t, R^a_t)$, $\{\sF^a_t\}_{t\geq 0}$ and $\{ \sF^{a,\prime}_t\}_{t\geq 0}$, respectively.  Clearly, $\sF^a_t$ and $ \sF^{a,\prime}_t $ are trivial for $t\in [0, a]$.

    \bigskip

 Let  $\wt \sG^B_t := \sG^B_t \vee \sG^L_{\infty}$  for $t\geq 0$ and for $a\geq 0$, 
   $$
  \wt \sF_t^{a}:=\wt  \sG^B_{L_{(t-a)^+} } \quad \hbox{for } t\geq 0.
 $$
 As $B_{L_{(t-a)^+}} \in \wt \sF_t^{a} $, we have $ \sF^{a, \prime}_t \subset \wt \sF_t^{a}$ for $t\in [0, T]$. 
 Note that  $R^a_t= a-t$ for $t\in [0, a]$ and for $t>a$,
 \begin{equation}\label{e:2.10}
  \{R^a_t > b\} = \{ S_{L_{t-a}} > t-a +b\} = \bigcup_{n\geq 1}\{ L_{t-a}=L_{t-a+b+1/n}\} \in \sG^L_\infty
 \quad \hbox{for }  b\geq 0.
 \end{equation} 
    Thus for each $t\geq 0$, $R^a_t$ is $\sG^L_\infty$-measurable. 
  Consequently, $\wt X^a= (X^a_t, R^a_t)$ is $\wt \sF^a_t$-measurable and so
  \begin{equation}\label{e:2.11}
  \sF^{a, \prime}_t \subset \sF^a_t \subset \wt \sF^a_t
  \quad \hbox{for every } t\geq 0.
  \end{equation} 
  Let $T>a\geq 0$ be constants. For $b<T-t$, the events $\{L_{t-a}=L_{t-a+b+1/n}\}$ increase in $n$ and belong to $\sG^L_{T-a}$ for large $n$. Hence, by \eqref{e:2.10}, $R^a_t\wedge (T-t) $ is $\sG^L_{T-a}$-measurable
  for every $t\in [0, T]$.

   Since $B$ and $L$ are independent, $B$ is an $\{  \wt \sG^B_t\}_{t\geq 0}$-martingale and 
 so  $B_{L_{(t-a)^+}}$ is a square-integrable $\{\wt \sF^a_t\}_{0\leq t\leq T}$-martingale in view of Theorem \ref{T:2.1}. 
On the other hand, $S_t= \inf\{r>0: L_r>t\}$ is $\sG^L_\infty$-measurable and $ \sG^L_\infty = \wt \sG^B_0= \wt \sF^a_0 $.  Since $\Omega$ is Polish, a family of probability measures $\mathbb Q_t$ on $\widetilde{\mathcal F}^a_t$, $t\geq 0$, with $\mathbb Q_t=\mathbb Q_s|_{\widetilde{\mathcal F}^a_t}$ for $t\leq s$, extends uniquely to a probability measure on $\widetilde{\mathcal F}^a_\infty:=\sigma( \cup_{t\geq 0} \widetilde{\mathcal F}^a_t )$; see \cite[Proposition~VIII.1.13]{RY}.
Two measures that are equivalent on each $\widetilde{\mathcal F}^a_t$ need not be equivalent on $\widetilde{\mathcal F}^a_\infty$.
Hence ``$\mathbb Q$-a.s.\ and hence $\mathbb P$-a.s.'' is used only for events in some $\widetilde{\mathcal F}^a_t$. 
  
  \section{Girsanov transform for sub-diffusions}\label{S:3}

 In this section, we study Girsanov transform for the sub-diffusion $X$.  Suppose $R_0=a$. Note that for each $t>0$, 
$$
a+S_t=\inf \{r>0: L_{(r-a)^+} > t\}=\inf \{r>0: \langle  X^{(1)}\rangle _r > t\} 
$$
is an $\{\sF^{a, \prime}_t\}_{t\geq 0}$-stopping time.  Here $X^{(1)}$ denote  the first coordinate of the multidimensional 
sub-diffusion $X^a$  and $\langle  X^{(1)}\rangle $ is the quadratic variation process of  the continuous martingale $X^{(1)}$. 
   Since $L_{S_t}=t$,   we have 
   \begin{equation} \label{e:3.2} 
   B_t=B_{L_{(a+S_t-a)^+}}  
   \quad  \hbox{ for } t\geq 0.
   \end{equation} 
Thus $B_t $ is  a standard Brownian motion on $\RR^{d}$ adapted to   the filtration $\{\sF^{a, \prime}_{R_0+S_t}\}_{t\geq 0}$. 
Moreover, 
$S_t= (a  + S_t)-a$ is $\sF^{a, \prime}_{a+S_t}$-measurable. 
Denote by $\{\sG_t; t\geq 0\}$ the minimum augmented  filtration generated by the L\'evy process $(B, S)$. 
We have from the above that 
 \begin{equation} \label{e:3.3a}
\sG_t \subset  \sF^{a, \prime}_{a+S_t}   \quad \hbox{for } t\geq 0. 
   \end{equation}
Observe that $L_{(t-a)^+}=\inf\{ r>0:S_{r}>(t-a)^+\}$ is a $\{\sG_t; t\geq 0\}$-stopping time.  
Thus  $\wt X^a_t$ is $\sG_{L_{(t-a)^+}}$-measurable and so
$ \sF^{a, \prime}_t \subset  \sG_{L_{(t-a)^+}}$ for every $ t\geq 0$. Hence
$$
\sF^{a, \prime}_{a+S_t}  \subset  \sG_{L_{S_t}} = \sG_t  \quad \hbox{for } t\geq 0. 
 $$
 This together with \eqref{e:3.3a} yields that for $a\geq 0$, 
 \begin{equation} \label{e:3.4a}
\sG_t = \sF^{a, \prime}_{a+S_t}   \quad \hbox{for } t\geq 0. 
   \end{equation}
 Consequently, 
 \begin{equation} \label{e:3.5a}
\sG_{L_{(t-a)^+}} = \sF^{a, \prime}_{a+S_{L_{(t-a)^+}}} = \sF^{a, \prime}_{t+R^a_t}   \quad \hbox{for } t\geq 0. 
   \end{equation}
Identity \eqref{e:3.4a} in particular implies that $\sG^B_\infty  \cup \sG^L_\infty \subset \sG_\infty \subset \sF^{a, \prime}_\infty$.  
 Thus we have 
 \begin{equation}\label{e:equal} 
 \sF^{a, \prime}_\infty = \sF^a_\infty = \wt \sF^a_\infty =\sG_\infty. 
 \end{equation} 
\medskip

For   semimartingales $N$ and $K$, we denote their quadratic  covariation  denoted by $[ N, K]$
and  the quadratic variation of $N$ by $[ N ]$, which is $[ N, N] $. When $[N, K]$ has locally integrable variation, 
we use $\langle N, K\rangle $ to denote its  \emph{compensator},
which  is the unique predictable finite-variation process $A$ with $[N, K]-A$ a local
martingale and $A_0=0$. It is known that
$$
[N, K]_t = \langle N^c, K^c\rangle _t + \sum_{0<s\leq t} (N_s -N_{s-}) (K_s -K_{s-}) \quad \hbox{for } t\geq 0,
$$ 
where $N^c$ and $K^c$ are the continuous local martingale parts of $N$ and $K$, respectively. 
For locally square-integrable martingales $N$ and $K$,  
 $NK -\langle N,K\rangle$ is a local martingale. We denote $\langle  N, N \rangle $ by $\langle  N\rangle $.  
See, e.g., \cite{HWY}, for above facts.

 \subsection{Structure-Preserving Girsanov transform}

\begin{theorem} \label{G1}
Let $a\geq 0$. 
 Suppose  $b(t, \omega)$ is an $\RR^d$-valued random process that is progressively measurable with respect to 
 the filtration $\{   \sF^{a, \prime}_t; t\geq 0\}$ so that  
 \begin{equation} \label{e:Novikov}
 \EE \exp \left( \frac12 \int_0^t   |b(s, \omega) |^2   dL_{(s-a)^+}  (\omega) \right)<\infty
 \quad \hbox{for every } t>0.
 \end{equation}  
 Define 
 \begin{equation}  \label{e:3.4}
M_t=\exp\left(  - \int_0^t b(s,\omega)     \cdot  d B_{L_{(s-a)^+}}  - \frac 1 2  \int _0^t |b  (s,\omega) |^2      d L_{(s-a)^+}     \right), \quad t\geq 0.
\end{equation}
Then $\{M_t; t\geq 0\}$  is a continuous  martingale with respect to the filtration $\{\wt \sF^a_t; t\geq 0\}$.   
 It uniquely determines a  probability measure $\QQ^*$ on $\wt \sF^a_\infty :=\sigma (\wt \sF^a_t; t\geq 0)$  so that 
  \begin{equation}  \label{e:3.5}
d\QQ^*=M_t d \PP   \quad  \hbox{on } \wt \sF^a_t   \ \hbox{ for each } t\geq 0.   
\end{equation}
   Under $\QQ^*$, $Z_t=B_t + \int_0^{a+S_t}  b(s,\omega)      d  {L_{(s-a)^+}}$ is a  standard Brownian motion on $\RR^d$
    with respect  to the filtration $\{\wt \sF^{a}_{a+S_t}\}_{t\geq 0}$,
     and $S_t$ is a subordinator with drift $\kappa \geq 0$
 and L\'evy measure $\nu$ that is independent of  $Z$.
  In particular, 
  $$
  Z_{L_{(t-a)^+}} = B_{L_{(t-a)^+}} + \int_0^t b(s,  \omega) dL_{(s-a)^+} 
  $$ 
  is a $d$-dimensional sub-diffusion under $\QQ^*$ that has the same distribution as $B_{L_{(t-a)^+}}$ under $\PP$, 
  and is a square-integrable $\{\wt \sF^a_t\}_{t\geq 0}$-martingale.  
     \end{theorem}

\begin{proof}    Note that $t\mapsto B_{L_{(t-a)^+}}$ is a square-integrable martingale with respect to the filtration $\{\wt \sF^a_t\}_{t\geq 0}$
having $\langle  B_{L_{(\cdot -a)^+}} \rangle _t =  L_{( t -a)^+}$.
Under condition \eqref{e:Novikov}, it follows from Novikov's criterion (see, e.g., \cite[Proposition VIII.1.15]{RY})
that $\{M_t; t\geq 0\}$ is a continuous  martingale with respect to the filtration $\{\wt \sF^a_t; t\geq 0\}$.   
So by \cite[Proposition VIII.1.13]{RY}, 
$d\QQ^*=M_t d \PP$   on $\wt \sF^a_t$ for each $t>0$ uniquely defines a probability measure $\QQ^*$ on $\wt \sF^a_\infty $.

Let $ Z_t := B_t+ \int_0^{a+S_t}  b(s,\omega)      d  {L_{(s-a)^+}}$.
   Write  $Z_t=(Z^{(1)}_t, \cdots, Z^{(d)}_t)$ and $B_t=(B^{(1)}_t, \cdots, B^{(d)}_t)$,  
$b(t, \omega)= (b^{(1)} (t, \omega), \cdots, b^{(d)} (t, \omega) )$. 
   As $B$ is a standard Brownian motion on $\RR^d$ with respect to the filtration   $\{  \wt \sG^B_t\}_{t\geq 0} = 
 \{\wt \sF^{a}_{a+S_t}\}_{t\geq 0}$, by  Girsanov's theorem for continuous martingales (see, e.g., \cite[Corollary 12.15]{HWY} or 
 \cite[Theorem VIII.1.4]{RY}),  
   under the probability measure $\QQ^*$, for each $1\leq i \leq d$,  
 \begin{eqnarray*}
&&  B^{(i)}_t- \Big\langle  B^{(i)}_\cdot  , \,  - \int_0^{a+S_\cdot}  b(s,\omega)      d B_{L_{(s-a)^+}}    \Big\rangle _t   \\
& =& B^{(i)}_t +  \Big\langle   \int_0^{a+S_\cdot} dB^{(i)}_{L_{(s-a)^+}} ,   \,  \int_0^{a+S_\cdot}  b(s,\omega)      d B_{L_{(s-a)^+}}    \Big\rangle _t  \\
& =&
B^{(i)}_t + \int_0^{a+S_t} b^{(i)}(s,\omega)      d  {L_{(s-a)^+}}=Z^{(i)}_t
  \end{eqnarray*}
is a continuous local martingale with respect to the filtration $\{\wt \sF^a_{a+S_t}\}_{t\geq 0}$
 with $\langle Z^{(i)}, Z^{(k)} \rangle_t =\langle B^{(i)}, B^{(k)} \rangle_t = \delta_{ik}  t$ for any $1\leq i, k\leq d$.
 Here $\delta_{ik}$ is the Kronecker delta, which equals $1$ if $i=k$ and $0$ if $i\not= k$. 
By L\'evy's characterization of Brownian motion (see, e.g., \cite[Theorem IV.3.6]{RY}), $Z$ is a $d$-dimensional Brownian motion
with respect to the filtration $\{\wt \sF^a_{a+S_t}\}_{t\geq 0}$.
Note that $ \int_0^{a+S_t}  b(s,\omega)      d  {L_{(s-a)^+}}$ is $\sF^{a, \prime}_{a+S_t} $-measurable so 
by \eqref{e:3.4a}, $Z_t$ is $\sG_t$-measurable. 
By \eqref{e:3.4a}, $\sG_t = \sF^{a, \prime}_{a+S_t} \subset \wt \sF^a_t$. Hence 
$Z$ is a $d$-dimensional Brownian motion with respect to the smaller filtration $\{\sG_t \}_{t\geq 0}$ as well.

  For any $f(x, y)\in C_c^2 (\RR^d  \times \RR)$,  by Ito's formula and the L\'evy system for the subordinator $S$, 
  under probability measure $\PP$, 
\begin{eqnarray}
&& f(Z_t,  S_t) \nonumber \\
&=& f(0, 0)+ \int_0^t \nabla_x f  (Z_r,   S_r) dZ_r  
  +   \int_0^t  f_y (Z_r,      S_{r-}) dS_r     +\frac12 \int_0^t \Delta_x f  (Z_r,  S_r) d \langle  Z\rangle _r  
   \nonumber \\
&&  +\sum_{r\in (0, t]} \left( f  (Z_r,  S_r) -  f  (Z_r,  S_{r-}) - f_y (Z_r,  S_{r-}) (S_r-S_{r-}) \right) \nonumber \\
&=& f(0, 0)+ \int_0^t  \nabla_x f (Z_r,  S_r) dZ_r   
  + \kappa \int_0^t f_y (Z_r,  S_r) d r     +\frac12 \int_0^t \Delta_x f  (Z_r,  S_r) d r       \nonumber \\
&&  +    M^d_t  + \int_0^t  \int_{(0, \infty)} ( f  (Z_r,  S_r +z ) -  f  (Z_r,  S_r)   ) \nu (dz) dr ,    \label{e:3.3}
 \end{eqnarray} 
 where  
 $\{M^d_t; t\geq 0\}$ is a purely discontinuous  $\{ \sG_t\}_{t\geq 0}$-martingale with 
 $$
 M^d_t -M^d_{t-} =  f  (Z_t,    S_t) -  f  (Z_t,   S_{t-}) \quad \hbox{for } t>0.
 $$
  Here $\nabla_x f(x, y):= (\frac{\partial}{\partial x_1} f(x, y), \cdots, \frac{\partial}{\partial x_d} f(x, y))$, 
     $\Delta_x f(x,  y):=\sum_{i=1}^d \frac{\partial^2}{\partial x_i^2}  f(x, y)$ and 
$f_y(x, y):= \frac{\partial}{\partial y}  f(x, y)$.  

By  Girsanov's theorem for purely discontinuous martingales (see, e.g., \cite[Corollary 12.16]{HWY}), 
under $\QQ^*$, $M^d$ remains to be a purely discontinuous  martingale with respect to the filtration $\{ \sG_t\}_{t\geq 0}= 
\{\sF^{a, \prime}_{a+S_t}\}_{t\geq 0} $.
Thus it follows from   \eqref{e:3.3} that  under $\QQ^*$,   
 \begin{eqnarray*}
&&  f(Z_t, S_t)  -   \frac12 \int_0^t  \Delta_x   f (Z_r,  S_r) d r  
   -  \kappa \int_0^t  f_y (Z_r,  S_{r-}) d r      \\ 
   && \quad     - \int_0^t  \int_{(0, \infty)} ( f  (Z_r,  S_r +z ) -  f  (Z_r,  S_r)   ) \nu (dz) dr 
 \end{eqnarray*} 
 is a bounded $\{\sG_t\}_{t\geq 0}$-martingale on $[0,T]$ for every $T>0$ and $f\in C_c^2(\RR^d   \times \RR)$.  By the martingale characterization of L\'evy processes, 
 we conclude that under $\QQ^*$,   $S$ is a subordinator with drift $\kappa \geq 0$
 and L\'evy measure $\nu$ that is independent of $Z$.  Consequently,  
  as $t\mapsto L_t$ does not change on  each random time interval  $[S_{t-}, S_t]$, we conclude 
   that 
  $$        
     Z_{L_{(t-a)^+}} =  B_{L_{(t-a)^+}} + \int_0^{a+S_{L_{(t-a)^+}}}  b(s,  \omega) dL_{(s-a)^+} 
     = B_{L_{(t-a)^+}} + \int_0^t b(s,  \omega) dL_{(s-a)^+}
  $$ 
     is a $d$-dimensional sub-diffusion under $\QQ^*$,
  which has the same distribution as $B_{L_{(t-a)^+}}$ under $\PP$. 
  Moreover, as $Z$ is   a $d$-dimensional  Brownian motion
with respect to the filtration $\{\wt \sF^a_{a+S_t}\}_{t\geq 0}$,
$\{Z_{L_{(t-a)^+}}; t\geq 0\}$ is a square-integrable martingale with respect to the filtration
$\big\{ \wt \sF^a_{a+S_{L_{(t-a)^+}}} \big\}_{t\geq 0} = \{\wt \sF^a_t\}_{t\geq 0}$.
\end{proof}

 \medskip

 \begin{remark} \label{R:3.2} \rm
 \begin{enumerate}
 \item[(i)]    In view of Theorem \ref{T:2.1} and $\phi(\infty)=\infty$ (which holds because $S$ is strictly increasing), condition \eqref{e:Novikov} is satisfied 
 if  $b(s, \omega)$ is bounded on $[0, t] \times \Omega$ for each $t>0$.
In fact, when $b(s, \omega)$ is bounded on $[0, T] \times \Omega$, it follows from Theorem \ref{T:2.1} and Cauchy-Schwarz inequality that 
 the exponential martingale $\{M_t; t\in [0, T]\}$ defined by \eqref{e:3.4}  has finite moment of any order. Indeed, 
 for each integer $k\geq 1$ and $T>0$,
 
 \begin{align*}
 \sup_{t\in [0, T]}  \EE [ M_t^k ]  
 &= \sup_{t\in [0, T]}  \EE \exp\Big(  -k  \int_0^t b(s,\omega)     \cdot  d B_{L_{(s-a)^+}} 
      - \frac k 2  \int _0^t |b  (s,\omega) |^2      d L_{(s-a)^+}     \Big) \\
  &\leq  \sup_{t\in [0, T]}  
  \Big( \EE  \exp\Big(  -2 k  \int_0^t b(s,\omega)     \cdot  d B_{L_{(s-a)^+}}  -\frac {{(2k)}^2} 2  \int _0^t |b  (s,\omega) |^2  d L_{(s-a)^+}  \Big)
  \Big)^{1/2}  \\ 
  &\qquad \times   \Big( \EE  \exp\Big(  (2k^2 -k) \int _0^t |b  (s,\omega) |^2      d L_{(s-a)^+}  \Big)  \Big)^{1/2}\\
  &\leq   \Big( \EE  \exp\Big(  (2k^2 -k) \int _0^T |b  (s,\omega) |^2      d L_{(s-a)^+} \Big)  \Big)^{1/2} <\infty. 
   \end{align*} 
   
   \item[(ii)] The exponential martingale $\{ M_t; t\geq 0 \}$ of \eqref{e:3.4} is also a continuous  
   $ \{ \sF^{a, \prime}_t; t\geq 0\}$-martingale because  each $M_t$ is $\sF^{a, \prime}_t$-measurable.  
   It     uniquely determines a  probability measure $\wh \QQ^*$ on $ \sF^{a, \prime}_\infty :=\sigma (\sF^{a, \prime}_t; t\geq 0)$  so that 
  \begin{equation}  \label{e:3.9}
d \wh \QQ^*=M_t d \PP   \quad  \hbox{on }   \sF^{a, \prime}_t   \ \hbox{ for each } t\geq 0.   
\end{equation}
Clearly, $\wh \QQ^* = \QQ^*$ on  $\sF^{a, \prime}_\infty$. 
   
       As noted in the proof, $\{Z_t; t\geq 0\}$ is  progressively measurable with respect to the filtration
   $\{ \sG_t \}_{t\geq 0} = \{ \sF^{a, \prime}_{a+S_t} \}_{t\geq 0}$,  and 
   $\{  Z_{L_{(t-a)^+}} ; t\geq 0\}$ is $\{ \sF^{a, \prime}_t\}_{t\geq 0}$-progressively measurable.
   So Theorem  \ref{G1} in particular implies that  under $\wh \QQ^*$, 
    $\{Z_t; t\geq 0\}$ is an $\{ \sG_t \}_{t\geq 0}$-Brownian motion and $\{  Z_{L_{(t-a)^+}} ; t\geq 0\}$ is a 
    square-integrable $\{ \sF^{a, \prime}_t\}_{t\geq 0}$-martingale. 
  
  \item[(iii)]
   Beyond being  a sub-diffusion, an important conclusion of Theorem  \ref{G1} is that  under $\QQ^*,$ 
      $\{  Z_{L_{(t-a)^+}} ; t\geq 0\}$  is a square-integrable martingale with respect to the enlarged filtration $\{\wt \sF^a_t\}_{t\geq 0}$. 
      \qed
   \end{enumerate}
    \end{remark}

 \medskip

 \begin{remark} \label{R:3.3} \rm
 The change of measure theorem for sub-diffusions in Theorem \ref{G1} is stronger than the usual
Girsanov theorem for continuous martingales in that under the transformed measure, the martingale part 
of the original Brownian motion
and the subordinator again have the required independence structure, so that the martingale part of 
 the original sub-diffusion is again a sub-diffusion. 
 In Theorem \ref{T:3.5} below, we further establish the  uniqueness of the measure change within this
structure-preserving class. 
 \qed
\end{remark}
   
 \medskip
    
 We recall the following martingale representation theorem from \cite[Theorem 2.1]{ZC2} by taking $T_0=\infty$ there. 
 
 \begin{theorem}\label{T:MRT}
Let $0\leq a<T$. 
  For every  $\eta \in L^2( \wt \sF^a_T)$,  
there exists an $\RR^d$-valued $\{\wt \sF^a_s \}_{0\leq s\leq T}$-progressively measurable process   $\{  \xi_s;  s\in [0, T]\}$   
 having  $\EE \int_0^T |   \xi_s  |^2 dL_{(s-a)^+}  <\infty$
so that  
\begin{equation}\label{e:M1.1}
\eta  =  \EE \big[ \eta \big| \wt \sF^a_0\big]  +  \int_0^{T}    \xi_s  \cdot  dB_{L_{(s-a)^+} }  .
\end{equation}
Such $ \xi$ is unique in the sense that if $ \xi'$ is another $\RR^d$-valued $\{ \wt \sF^a_t\}_{t\in [0, T]}$-progressively measurable process
having  $\EE \int_0^T |  \xi'_s  |^2 dL_{(s-a)^+}  <\infty$ so that \eqref{e:M1.1} holds with $\xi'$ in place of $\xi$, then
$\EE  \int_0^T |  \xi_s-   \xi'_s|^2 dL_{(s-a)^+} =0$.
\end{theorem}

 \begin{remark} \label{R:3.5} \rm
 For any continuous local martingale $M$,  stochastic integrals are well defined for progressively measurable integrands;
 see for instance \cite{KS, RY}. In fact, it is known (see, e.g., \cite[Theorem 3.10]{ChW}), for any progressively measurable
 process $\xi$, there is a predictable process $xi'$ so that $\EE \int_0^\infty |\xi_s -\xi'_s|^2  d \langle  M \rangle _s =0$. 
 Thus in Theorem \ref{T:MRT}, the integrand $\xi$ can be replaced by its predictable modification whenever needed. 
 \qed
\end{remark}

A consequence of Theorem \ref{T:MRT} is that every $\{ \wt \sF^a_t\}_{t\geq 0}$-martingale has a continuous modification.

\begin{corollary}\label{C:3.5}
 Suppose that $\{K_t; t\geq 0\}$ is an  $\{ \wt \sF^a_t\}_{t\geq 0}$-martingale.
 Then there is a continuous $\{ \wt \sF^a_t\}_{t\geq 0}$-martingale $\{\wt K_t; t\geq 0\}$ so that $\wt K_t=K_t$ $\PP$-a.s. for every $t\geq 0$.
\end{corollary}

\begin{proof} By \cite[Theorem 2.44 and Corollary 2.48]{HWY}, there is an  
  $\{\wt \sF^a_t\}_{t\geq 0}$-martingale  $\{\wt K_t; t\geq 0\}$ that is right continuous having left limits such that 
$\PP(\wt K_t=K_t)=1$   for each $t\geq 0$.  
  Fix any $T>a \geq 0$. For each integer $n\geq 1$, by Theorem \ref{T:MRT}, the bounded martingale $\EE [  K_T  \1_{\{|K_T| \leq n\}} | \wt \sF^a_t]$ has a continuous modification $K^n_t$ for $t\in [0, T]$
that has the expression
$$
K^n_t = \EE \big[ K_T  \1_{\{|K_T| \leq n\}} \big| \wt \sF^a_0 \big] + \int_0^t  \xi^n_s \cdot d B_{L_{(s-a)^+}},  \quad t\in [0, T]
$$
for some  $\{\wt \sF^a_s \}_{0\leq s\leq T}$-progressively measurable process   $\{  \xi^n_s;  s\in [0, T]\}$  
 having  $\EE \int_0^T |   \xi^n_s  |^2 dL_{(s-a)^+}  <\infty$.
  By Doob's $L^1$-maximal weak-type inequality (see, e.g., \cite[Theorem 2.42]{HWY}) and the Lebesgue dominated convergence theorem,
  for any $\eps >0$, 
 $$
 \lim_{n\to \infty}   \PP \Big( \max_{t\in [0, T]} |K^n_t - \wt K_t| \geq \eps \Big) 
 \leq \lim_{n\to \infty}  \eps^{-1} \EE | K^n_T- \wt K_T| 
 = \eps^{-1}  \lim_{n\to \infty}  \EE \left[ |K_T| \1_{\{|K_T|>n\}} \right]  =0.
$$
Hence there is a subsequence $\{n_k; k\geq 1\}$ so that $\lim_{k\to \infty} \max_{t\in [0, T]} |K^{n_k}_t -\wt K_t| =0$.
This proves that $\PP$-a.s., $\wt K_t$ is continuous in $t\in [0, T]$ and hence in  $t\in [0, \infty)$.
\end{proof}

\medskip

Using a stopping argument, one deduces immediately from the above corollary that 
every $\{ \wt \sF^a_t\}_{t\geq 0}$-local-martingale has a continuous modification.

\medskip

 Fix $b$ as in Theorem~\ref{G1}. Denote by $ \wsP$ the family of probability measures $\QQ$ on $\wt \sF^a_\infty$ that 
  is equivalent to $ \PP$   on $\wt \sF^a_t$ for each $t\geq 0$
 under which
  \begin{equation}\label{e:4.1}
  Y^a_t:= B_{L_{(t-a)^+}} + \int_0^t b(s,  \omega) dL_{(s-a)^+} , \quad t\geq 0, 
  \end{equation} 
   is a continuous $\{\wt \sF^a_t\}_{t\geq 0}$-local-martingale. 
The next theorem   characterizes the family $\wsP$.

  \begin{theorem}\label{T:EMM}
 Let $a\geq 0$, and   $b(t, \omega)$ be an $\RR^d$-valued random process that is progressively measurable with respect to $\{ \sF^{a, \prime}_t\}_{ t\geq 0}$ satisfying 
 the Novikov condition \eqref{e:Novikov}. A  probability measure  $\QQ$ on $(\Omega, \wt \sF^a_\infty)$ 
 is in $\wsP$  if and only if   $d\QQ= K_0 d\QQ^*$ for some strictly positive random variable $K_0 \in \wt \sF^a_0 = \sG^L_\infty$ having $\EE [ K_0]=1$.
 Here $\QQ^*$ is the probability measure on $(\Omega, \wt \sF^a_\infty)$ in Theorem \ref{G1}. 
 \end{theorem}

\begin{proof}  (i) Suppose that  $\QQ$  is a  probability measure  on $(\Omega, \wt \sF^a_\infty)$ that is 
in $\wsP$.  By Corollary \ref{C:3.5}, there is a continuous 
  $\{\wt \sF^a_t\}_{t\geq 0}$-martingale  $\{K_t; t\geq 0\}$   so that 
$K_t=\frac{d\QQ}{d \PP}|_{\wt \sF^a_t}$ $\PP$-a.s. for every  $t\geq 0$.  

 Fix an arbitrary $T>a\geq 0$. For an integer $n\geq 1$, define $T_n:=\inf\{t\geq 0: K_t \geq n\}\wedge T$ and $\Omega_n:=\{K_0<n\}\in\wt\sF^a_0$.
Then $\lim_{n\to \infty} T_n=T$ $\PP$-a.s., and $K_{ T_n}\1_{\Omega_n}$ is bounded and $\wt \sF^a_T$-measurable.
Since $K_0<\infty$, we have $ \Omega_n $ increases to $\Omega$  $\PP$-a.s.\ as $n\to\infty$.
 By Theorem \ref{T:MRT},  
 there exists an $\{\wt \sF^a_s \}_{0\leq s\leq T}$-progressively measurable process   $\{  \xi^n_s;  s\in [0, T]\}$   
 having  $\EE \int_0^T |   \xi^n_s  |^2 dL_{(s-a)^+}  <\infty$
so that  
\begin{equation}\label{e:1.1}
\1_{\Omega_n}K_{  T_n}  =  \1_{\Omega_n}\EE \big[ K_{T_n} \big| \wt \sF^a_0\big]  +  \int_0^{T}    \xi^n_s    \cdot d B_{L_{(s-a)^+} }  .
\end{equation}
Note that by a   Doob's optional stopping theorem (see \cite[Theorem 2.58]{HWY}),
 $\{K_{t\wedge T_n}; t\geq 0\}$ is a $\PP$-martingale with respect to  the filtration $\{\wt \sF^a_{t\wedge T_n}; t\geq 0\}$.
 We thus have from \eqref{e:1.1} that   $\PP$-a.s., 
\begin{equation}\label{e:n3.13} 
\1_{\Omega_n}K_{t\wedge T_n} = \EE \left[ \1_{\Omega_n}K_{T_n} | \wt \sF^a_{t\wedge T_n} \right] 
=  \1_{\Omega_n}\EE \big[ K_{T_n} \big| \wt \sF^a_0\big]  +  \int_0^{t \wedge T_n}    \xi^n_s    \cdot d B_{L_{(s-a)^+} } \quad \hbox{for } t\in [0, T].
\end{equation} 
  By Girsanov theorem (see, e.g., \cite[Corollary 12.15]{HWY}),  under $\QQ$, 
$$
V_t:= B_{L_{(t-a)^+}}-\int_0^t \frac{1}{K_s} d \langle  B_{L_{\cdot-a)^+}}, K \rangle _s  
$$
is a $\QQ$-local martingale with respect to the filtration $\{\wt \sF^a_t\}_{t\in [0, T]}$. 
On the other hand,  by assumption under $\QQ$, 
$$
Y_t:= B_{L_{(t-a)^+}} + \int_0^t b(s,  \omega) dL_{(s-a)^+}, \quad t\geq 0,
$$
 is a continuous local martingale with respect to $\{\wt \sF^{a }_t\}$.  
It follows from the last two displays that under $\QQ$
$$
  Y_{t\wedge T_n} -V_{t\wedge T_n}  = \int_0^t b(s,  \omega) \1_{\{s\leq T_n\} } dL_{(s-a)^+} 
  +  \int_0^t \frac{\xi^n_s}{K_s}  \1_{\{s\leq T_n\} } dL_{(s-a)^+}
 $$
is both a continuous local martingale with respect to $\{\wt \sF^a_{t\wedge T_n}\}_{t\in [0, T]}$
 and a continuous process of finite variation. 
Hence we have $\QQ$-a.s.,
$$ 
\int_0^t \left( b(s,  \omega)   +   \frac{\xi^n_s}{K_s}   \right)   \1_{\{s\leq T_n\} }  dL_{(s-a)^+}=0 \quad \hbox{for every } t\in [0, T].
$$
This implies that $\QQ$-a.s. and hence $\PP$-a.s. that  
$$
\int_0^T |  b(s,  \omega)   +  (\xi^n_s/K_s) |  \1_{\{s\leq T_n\} }    dL_{(s-a)^+}=0. 
$$
Thus we see from \eqref{e:n3.13} that $K_{t\wedge T_n}$ satisfies the equation
\begin{equation}\label{e:n3.14} 
\1_{\Omega_n}K_{t\wedge T_n}= \1_{\Omega_n}K_0 -  \int_0^t    \1_{\Omega_n}K_{s\wedge T_n} b(s, \omega)   \1_{\{s\leq T_n\} }  \cdot d B_{L_{(s-a)^+} }
\quad \hbox{for } t\in [0, T],
\end{equation} 
where $K_0 = \EE \big[ K_T \big| \wt \sF^a_0\big]$. 
It is well known (see, e.g., \cite[Theorem 9.39]{HWY}) that the above stochastic integral equation has a unique solution  given by 
\begin{align*}
  \1_{\Omega_n}K_{t\wedge T_n}  =  \1_{\Omega_n}K_0  \exp\Big(  - \int_0^t b(s,\omega)   \1_{\{s\leq T_n\} }   \cdot  d B_{L_{(s-a)^+}}   - \frac 1 2  \int _0^t |b  (s,\omega) |^2   \1_{\{s\leq T_n\} }    d L_{(s-a)^+}     \Big)   
  \end{align*}
 for $ t\in [0, T]$.  Taking $n\to \infty$ (so that $\1_{\Omega_n}\uparrow 1$ and $T_n\to T$) yields that $K_t = K_0 M_t$ for $t\in [0, T]$, where $M_t$ is 
 the Dol\'eans-Dade exponential martingale defined by \eqref{e:3.4}.
    Since this holds for any $T>a$, we have $K_t = K_0M_t$ for every $t\geq 0$. 
Consequently, $d \QQ = K_0 d \QQ^*$.
\medskip

(ii) Conversely, suppose that $d\QQ =K_0 d \QQ^*$ with $K_0 \in \wt \sF^a_0 = \sG^L_\infty$ is strictly positive $\PP$-a.s. having $\EE [K_0]=1$. 
  By Theorem \ref{G1}, $M_t$ is a strictly positive $\{\wt \sF^a_t\}_{t\geq 0}$-martingale with $M_0=1$.
   We know from Theorem \ref{G1} that under $\QQ^*$,
   $Y^a_t$ of \eqref{e:4.1} 
   is a square-integrable $\{\wt \sF^a_t\}_{t\geq 0}$-martingale. As $K_0 \in \wt \sF^a_0$,   
    it follows readily that  $Y^a_t$ is an  $\{\wt \sF^a_t\}_{t\geq 0}$-martingale under $\QQ$. 
    So $\QQ \in \wsP$. 
    This completes the proof of the theorem.
\end{proof}

\begin{theorem} \label{T:3.5} 
Let $a\geq 0$.  Under the assumption of Theorem \ref{G1}, 
the exponential martingale   $\{M_t; t\in [0, \infty)\}$ defined by \eqref{e:3.4}
   is the unique strictly positive   $\{\wt \sF^a_t\}_{t\geq 0}$-martingale  with $\EE M_0=1$ such that 
   under the probability measure $\QQ^*$ defined by \eqref{e:3.5}, 
   the process $\{Y^a_t; t\geq 0\}$ of  \eqref{e:4.1} 
  is a $d$-dimensional sub-diffusion that has the same distribution as $B_{L_{(t-a)^+}}$ under $\PP$
  and is an  $\{\wt \sF^a_t\}_{t\geq 0}$-local-martingale. 
  \end{theorem}
  
  \begin{proof} Suppose that $\{K_t; t\geq 0\}$ is a strictly positive  $\{\wt \sF^a_t\}_{t\geq 0}$-martingale  with $\EE K_0=1$
such that under the probability measure $\QQ$ defined by $d\QQ = K_t d\PP$ on $\wt \sF^a_t$ for each $t\geq 0$,
$Y^a_t$ is a $d$-dimensional sub-diffusion that has the same distribution as $B_{L_{(t-a)^+}}$ under $\PP$
and is a  $\{\wt \sF^a_t\}_{t\geq 0}$-local martingale.  By Theorem \ref{T:EMM}, $d\QQ = K_0 d \QQ^*$
and $K_t=K_0 M_t$ for $t\geq 0$.
Then in view of Theorem \ref{G1}, $Y^a$ has the same law under $\QQ^*$ and $\QQ$. 
Thus for any non-negative $\xi \in \sigma (Y^a):=\sigma(Y^a_s; s\geq 0)$,
$$ 
\EE^{\QQ^*} [ \xi] = \EE^{\QQ} [ \xi] = \EE^{\QQ^*} [ K_0 \xi ].
$$
As $\langle  Y^a\rangle _t = L_{(t-a)^+}$, $\wt \sF^a_0= \sG^L_\infty \subset \sigma (Y^a)$. 
Taking $\xi =K_0$ in the above display yields
 $$
  \EE^{\QQ^*} [ K_0^2] =  \EE^{\QQ^*} [K_0]=1.
 $$ 
It follows that $ \EE^{\QQ^*}  \big[ (K_0 -1)^2] =0$, that is,  $K_0=1$ $\QQ^*$-a.s. and hence $\PP$-a.s.
This shows $K_t=M_t$ for every $t\geq 0$, which establishes the theorem.
\end{proof}

\medskip

As an application of Theorem \ref{G1}, we have the following. 

\begin{theorem} \label{G2}
Let $k\geq 1$ be an integer and $a \geq 0$.
Suppose that   $\beta (t,\omega) $ is an $\RR^k$-valued stochastic process that is $\{\sF^{a, \prime}_t\}_{t\geq 0}$-progressively measurable
and  $\theta (t,\omega) $ is an $ \RR^{k\times d}$-valued  $\{\sF^{a, \prime}_t\}_{t\geq 0}$-progressively measurable stochastic process
so that for each $t>0$, $\int_0^t ( |\beta (s, \omega)| + | \theta (s, \omega)|^2  ) \, dL_{(s-a)^+} <\infty$ 
$\PP$-a.s.

Let 
$$ 
d \Xi _t= \beta (t,\omega) dL_{(t-a)^+}  + \theta (t, \omega) \cdot d B_{L_{(t-a)^+}},  \quad t \geq 0 .
 $$
Suppose there exist   $\{\sF^{a, \prime}_t\}_{t\geq 0}$-progressively measurable $\RR^d$-valued process $u(t,\omega)$ and  $\RR^k$-valued process $\alpha (t,\omega)$ such that
\begin{eqnarray}
\beta(t,\omega)= \theta(t,\omega)u(t,\omega)      + \alpha(t,\omega)
\end{eqnarray}
and $\EE \exp (\frac12 \int_0^t   |u (s, \omega ) |^2    d L_{(s-a)^+} )<\infty$  for every $t>0$. 
Define
\begin{eqnarray} \label{MG2}
M_t=\exp\left(  - \int_0^t   u(s,\omega)    \cdot   d B_{L_{(s-a)^+}}  - \frac 1 2  \int _0^t   |u  (s,\omega)|^2       dL_{(s-a)^+}   \right), \quad t\geq 0, 
\end{eqnarray}
Then $\{M_t; t\geq 0\}$  is a continuous  martingale with respect to the filtration $\{ \wt \sF^a_t; t\geq 0\}$.   
 Define  the probability measure $\QQ$ on $\wt \sF^a_\infty  $  by
  \begin{eqnarray} \label{QG2}  d\QQ=M_t d \PP   \qquad   \hbox{on } \wt \sF^a_t  \  \hbox{ for  each } t\geq 0.  
  \end{eqnarray}
  Under $\QQ$, $\widehat B_t :=B_t + \int_0^{a+S_t}  u(s,\omega)      d  {L_{(s-a)^+}}$ is a  standard Brownian motion on $\RR^d$
  with respect to the filtration $\{\wt \sF^a_{a+S_t}\}_{t\geq 0}$,
     and $S_t$ is a subordinator with drift $\kappa \geq 0$
 and L\'evy measure $\nu$ that is independent of  $\widehat B$.
 Consequently,  
 \begin{eqnarray}\label{hatBG2}
\widehat  B_{L_{(t-a)^+}}=\int _0^t u(s,\omega )  dL_{(s-a)^+}  +B_{L_{(t-a)^+}},
\end{eqnarray}
 is  a $d$-dimensional sub-diffusion under $\QQ$ that has the same distribution   as $\{B_{L_{(t-a)^+}} \}_{t\geq 0}$ under $\PP$, 
  and $\Xi_t$ has the representation
 \begin{eqnarray}
 d\Xi_t=\alpha (t,\omega) dL_{(t-a)^+}+\theta(t,\omega )d\widehat  B_{L_{(t-a)^+}}.
 \end{eqnarray}
\end{theorem}

\begin{proof} Let $\beta,\theta,u,\alpha$ be as in the statement.  By Theorem \ref{G1} (applied with $b=u$),  we know that  $M=\{M_t; t\geq 0\}$ is a continuous $\{\wt \sF^a_t\}_{t\geq 0}$-martingale
and under probability measure $\QQ$, $S$ is a subordinator and 
$$
\widehat B_t :=B_t + \int_0^{a+S_t}  u(s,\omega)      d  {L_{(s-a)^+}}, \quad t\geq 0,
$$
 is an $\RR^d$-valued $\{ \sF^a_{a+S_t} \}_{t\geq 0}$-Brownian motion   independent of $ S$. 
 So   $\widehat  B_{L_{(t-a)^+}}=\int _0^t u(s,\omega )  dL_{(s-a)^+}  +B_{L_{(t-a)^+}}$
 is a sub-diffusion under $\QQ$ that has the same distribution   as $\{B_{L_{(t-a)^+}} \}_{t\geq 0}$ under $\PP$.
 Thus 
\begin{eqnarray*}
      d\Xi_t &=& \beta (t,\omega) dL_{(t-a)^+}  + \theta (t, \omega)( d \widehat  B_{L_{(t-a)^+}} -u(t,\omega)  d L_{(t-a)^+} )\nonumber\\
   &=& \left(\beta (t,\omega) -  \theta (t, \omega)u  (t, \omega)   \right)  d L_{(t-a)^+}  + \theta (t, \omega) d \widehat  B_{L_{(t-a)^+}}\nonumber\\
      &=& \alpha(t,\omega) d  L_{(t-a)^+}    + \theta (t, \omega) d\widehat  B_{L_{(t-a)^+}}.
         \end{eqnarray*}
  This establishes the theorem.
\end{proof} 
   
    \subsection{Characterization of equivalent martingale measures}
    
    We first recall the following Kunita-Watanabe decomposition theorem from Case 3 of \cite{AS}. 
    It is established in \cite{AS} for $d=1$ but the proof works for the $d\geq 2$ case as well. 
    For square-integrable martingales, this result is initially established in \cite{KW}. 
    
\begin{lemma}[Kunita-Watanabe decomposition] \label{L:KW}
Let $M=(M^{(1)}, \cdots, M^{(d)})$ 
be an $\RR^d$-valued  continuous local martingale with $\langle M^{(i)}, M^{(j)}\rangle  \equiv 0$ for any $i\not= j$, 
and $K$ an arbitrary real-valued local martingale with respect to some filtration
$\{\sH_t\}$. Then there is an $\RR^d$-valued $\{\sH_t\}$-predictable process $\{\beta = (\beta^{(1)}, \cdots, \beta^{(d)}) \}$ having 
$\sum_{i=1}^d \int_{0}^t |\beta^{(i)}_s |^2 \,d\langle M^{(i)}\rangle_s \le\langle K^{c}\rangle_t<\infty$, and   a rcll 
local martingale $N$ with $N_{0}=0$ and $[N,M]\equiv 0$  so that 
$$
  K_t =K_{0}+\int_0^t \beta_s \cdot dM_s +N_t  \quad \hbox{for } t\geq 0, 
$$
where $K^c$ is the continuous local martingale part of $K$.  The decomposition is
unique in the sense that $\beta$ is unique $d\langle M\rangle\otimes d\PP$-a.e., $N$ is unique up to indistinguishability, 
and $[K,M]_t =\langle K^{c},M\rangle_t =\int_0^t \beta_s \,d\langle M\rangle_s$.
\end{lemma}

Suppose that the condition of Theorem \ref{G1} holds. 
Theorem \ref{T:EMM} characterized the space $\wsP$  of 
all the probability measures $\QQ$ on $\wt \sF^a_\infty$ that is  equivalent to $\PP$ on $\wt \sF^a_t$ for every $t\geq 0$
under which the process   $\{  Y^a_t; t\geq 0\}$ of \eqref{e:4.1}    is a continuous  $\{\wt \sF^a_t\}_{t\geq 0}$-local martingale.
In many applications such as finance, one is interested in the martingale properties with respect to the smaller filtration
$\{\sF^{a, \prime}_t\}_{t\geq 0}$ that is generated by the sub-diffusion $X^a_t=X^a_0+ B_{L_{(t-a)^+}}$.
Denote by $\sP$ the set of probability measures $\QQ$ on $\sF^{a, \prime}_\infty$ 
that is  equivalent to $\PP$ on each $\sF^{a, \prime}_t$
under which    $\{  Y^a_t; t\geq 0\}$ of \eqref{e:4.1}    is a continuous  $\{\sF^{a, \prime}_t\}_{t\geq 0}$-local martingale.
Recall from \eqref{e:equal} that  $\sF^{a, \prime}_\infty = \wt \sF^a_\infty$. 
Thus  $\wsP\subset \sP$.    The next theorem characterizes $\sP$.  

\begin{theorem} \label{T:REMM}
Let $a\geq 0$. 
Suppose the condition of Theorem \ref{G1} holds and let $M$ be the exponential martingale defined by \eqref{e:3.4}. 
Then a probability measure $\QQ$ on $\sF^{a, \prime}_\infty$ is in  $\sP$ if and only if there is an   
 $\{\sF^{a, \prime}_t\}_{ t\geq 0}$-local-martingale $N$  with $N_0=0$, $\big\langle  N, B_{L_{(\cdot -a)^+}} \big\rangle  \equiv 0$
 and    $N_t-N_{t-} > -1$ for every $t>0$
so that 
\begin{equation} 
\frac{d\QQ}{d \PP} \Big|_{\sF^{a, \prime}_t} = M_t  \, \sE (N)_t \quad \hbox{for } t\geq 0.
\end{equation} 
Here $\sE (N)$ denotes the Dol\'eans-Dade exponential local martingale of $N$, that is, $\sE (N)$ is the unique solution of 
$dZ_t = Z_{t-}\, d N_t$ with $Z_0=1$. 
\end{theorem} 
 
\begin{proof} (i) Suppose that $\QQ \in \sP$.  Let $K_t:= \frac{d\QQ}{d \PP} \big|_{\sF^{a, \prime}_t}. $
Since the filtration $\{ \sF^{a, \prime}_t\}_{t\geq 0}$ is right continuous, it is well known that the strictly positive
$\{ \sF^{a, \prime}_t\}_{t\geq 0}$ -martingale $\{K_t; t\geq 0\}$ has a right continuous modification having left limits.
We represent $K$ by its rcll version. Define 
\begin{equation}\label{e:3.21}
A_t:= \int_0^t \frac{1}{K_{s-}} dK_s, \quad t\geq 0, 
\end{equation}
  which is a rcll  $\{ \sF^{a, \prime}_t\}_{t\geq 0}$-local-martingale. 
By Lemma \ref{L:KW} applied to the continuous local martingale $B_{L_{(t-a)^+}}$ and the local martingale $A$, 
there is an $\{ \sF^{a, \prime}_t\}_{t\geq 0}$-predictable process $\{\beta_s\}$ having 
$\int_{0}^t\beta^{2}_s \,d L_{(s-a)^+} \le\langle A^{c}\rangle_t<\infty$, and   a rcll 
local martingale $N$ with $N_{0}=0$, $\big[N, B_{L_{(\cdot -a)^+}} \big]\equiv 0$  so that 
$$
  A_t = \int_0^t \beta_s \cdot dB_{L_{(s-a)^+}}  +N_t  \quad \hbox{for } t\geq 0.
$$
Let 
\begin{equation}\label{e:3.22}
Z_t:= \exp (A_t -\tfrac12 \langle A^c\rangle _t ) \prod_{ 0<s\leq t} (1+   A_s-A_{s-}) e^{-(A_s-A_{s-} )}, \quad  t\geq 0.
\end{equation}
By \cite[Theorem 9.39]{HWY},  for almost every $\omega \in \Omega$, 
the infinite product  on the right hand side is absolutely convergent 
and $Z_t$ is the unique rcll local martingale satisfying  $Z_t=1+\int_0^t Z_{s-} dA_s$ (that is, 
$Z= \sE (A)$). 
 Observe that 
 $$
 N_t-N_{t-}= A_t-A_{t-} = \frac{1}{K_{t-}} (K_t-K_{t-}) > -1 \quad \hbox{for all } t>0.
 $$ 
 Consequently, $\PP (Z_t>0 \hbox{ for all } t\geq 0)=1$.   
On the other hand, $dK_t = K_{t-} dA_t$ by   \eqref{e:3.21} with $K_0=1$.  Thus $K_t = 1+ \int_0^t K_{s-} dA_s$
for all $t\geq 0$.  We conclude by uniqueness that $K_t=Z_t$ for all $t\geq 0$. 
Note that the continuous local martingale part $K^c$ of $K$ is 
$$
K^c_t = \int_0^t K_{s-} dA^c_s = \int_0^t K_{s-} \beta_s \cdot dB_{L_{(s-a)^+}} + \int_0^t K_{s-} dN^c_s.
$$
Since $t\mapsto B_{L_{(t-a)^+}}$ is continuous, 
\begin{equation}\label{e:3.23}
0\equiv \big[N, B_{L_{(\cdot -a)^+}} \big] =  \big<N^c , B_{L_{(\cdot -a)^+}} \big\rangle .
 \end{equation}
Thus 
$$
\langle B_{L_{(\cdot -a)^+}}, K^c\rangle _t = \int_0^t K_{s-} \beta_s d L_{(s-a)^+}.
$$
 By Girsanov theorem (see, e.g., \cite[Theorem 12.13]{HWY}), under $\QQ$, 
$$
B_{L_{(t-a)^+}} - \int_0^t \frac{1}{K_{s-}} d\langle B_{L_{(\cdot -a)^+}}, K^c\rangle _s  
= B_{L_{(t-a)^+}} - \int_0^t \beta_s d L_{(s-a)^+}
$$
is an $\{ \sF^{a, \prime}_t\}_{t\geq 0}$-local martingale under $\QQ$. 
 We conclude from $\QQ\in \sP$ that $\beta_s = - b(s, \omega)$ for all $s\geq 0$.  
Thus it follows from  \eqref{e:3.22} and \eqref{e:3.23} that 
$$
K_t = M_t  \exp (N_t -\tfrac12 \langle N^c\rangle _t ) \prod_{ 0<s\leq t} (1+   N_s-N_{s-}) e^{-(N_s-N_{s-} )}
= M_t  \, \sE (N)_t \quad \hbox{for } t\geq 0.
$$

(ii) Next suppose that $N$ is an   
 $\{\sF^{a, \prime}_t\}_{ t\geq 0}$-local martingale    with $N_0=0$, $\big\langle  N, B_{L_{(\cdot -a)^+}} \big\rangle  \equiv 0$
 and    $N_t-N_{t-} > -1$ for every $t>0$ so that $K_t:= M_t  \, \sE (N)_t$ is an $\{\sF^{a, \prime}_t\}_{ t\geq 0}$-martingale.
 Note that under the above assumption, the exponential local martingale $\sE (N)$ is strictly positive.  
 By the explicit formula for the Dol\'eans-Dade exponential local martingale, we have
 $$
 K  = \sE \Big( - \int_0^\cdot b (s, \omega) d B_{L_{(s-a)^+}} + N \Big)  .
 $$
  Let $\QQ$ be the probability measure $\QQ$ on  $\sF^{a, \prime}_\infty$ determined by 
   $\frac{d\QQ}{d \PP} \big|_{\sF^{a, \prime}_t}=K_t$ for each $t\geq 0$. 
 By Girsanov theorem,  
 $$
B_{L_{(t-a)^+}} -  \big\langle B_{L_{(\cdot -a)^+}},  - \int_0^\cdot b (s, \omega) d B_{L_{(s-a)^+}} + N^c \big\rangle _t  
= B_{L_{(t-a)^+}} + \int_0^t  b (s, \omega)  d L_{(s-a)^+}
$$
is an $\{\sF^{a, \prime}_t\}_{ t\geq 0}$-local martingale under $\QQ$. In other words, $\QQ\in \sP$.
This completes the proof of the theorem
\end{proof} 
 
 \medskip

The following uniqueness result is the counterpart of Theorem \ref{T:3.5} 
but with respect to the smaller filtration $\{ \sF^{a, \prime}_t\}_{t\geq 0}$.
   
\begin{theorem} \label{T:uniqueness} 
Let $a\geq 0$.  Under the assumption of Theorem \ref{G1}, 
the probability measure $\QQ^* $   
   is the unique probability measure in $\sP$
   under which   the process $\{Y^a_t; t\geq 0\}$ of  \eqref{e:4.1} 
   is a $d$-dimensional sub-diffusion that has the same distribution as $B_{L_{(t-a)^+}}$ under $\PP$. 
     \end{theorem}
 
 \begin{proof} By Theorem \ref{G1} and Remark \ref{R:3.2}(ii), $ \QQ^* \in \sP$ 
 and $\{Y^a_t; t\geq 0\}$ of  \eqref{e:4.1} 
 under $ \QQ^*$  is a $d$-dimensional sub-diffusion  that has the same distribution as $B_{L_{(t-a)^+}}$ under $\PP$.
 
 (i) Suppose $\QQ\in \sP$ is another probability measure  on $\sF^{a, \prime}_\infty$ that has the above stated property. 
 Let $K_t$ be the strictly positive rcll $\{\sF^{a, \prime}_t\}_{t\geq 0}$-martingale so that 
  $\frac{d \QQ}{d \PP}=K_t$ on $\sF^{a, \prime}_t$  for $t\geq 0$.
  By Girsanov theorem, $Y^a$ is a local martingale under $\QQ$ with respect to the filtration
 $\{\sF^{a, \prime}_t\}_{t\geq 0}$. 
  Define  
   \begin{equation}
   Z_t := Y^a_{a+S_t} = B_t + \int_0^{a+S_t}  b(s,\omega)      d  {L_{(s-a)^+}}, \quad t\geq 0, 
   \end{equation} 
   where the last equality is due to the fact that $t\mapsto L_t$ does not change on  each random time interval  $[S_{t-}, S_t]$.
   By the proof of Theorem \ref{G1},  under $\QQ$, $Z$ is  a $d$-dimensional Brownian motion     with respect to the filtration 
    $\{ \sF^{a, \prime}_{a+S_t} \}_{t\geq 0}$, which is $\{\sG_t\}_{t\geq 0}$ by \eqref{e:3.4a},      
    and $S$ is still a subordinator with drift $\kappa \geq 0$ and L\'evy measure $\nu$ and is independent of $Z$. Indeed, $Z$ is a $\{\sG_t\}$-Brownian motion by L\'evy's characterization, since $\langle Y^a\rangle=L_{(\cdot-a)^+}$.
Moreover, $R^a_{a+S_t}=0$ for every $t\geq 0$, and $(Y^a,R^a)$ has the same law under $\QQ$ and $\PP$.
So by the strong Markov property (Theorem~\ref{T:2.3}) at the $\{\sF^{a,\prime}_s\}$-stopping time $a+S_t$, the increments of $(Z,S)$ after time $t$ are independent of $\sG_t=\sF^{a,\prime}_{a+S_t}$, and have the same law as under $\PP$.
     We conclude from the independence between $Z$ and $S$  under $\QQ$ that  
    $Z$ is a $d$-dimensional  Brownian motion   with respect to the filtration 
    $\{ \sG_t \vee \sG^L_\infty \}_{t\geq 0} = \{ \wt \sG^B_t\}_{t\geq 0}$ under $\QQ$.   
    Hence under $\QQ$, 
    $Y^a_t = Z_{L_{(t-a)^+}}$ is a local martingale with respect to the filtration  $\{ \wt \sG^B_{L_{(t-a)^+}} \}_{t\geq 0}
    =\{ \wt \sF^a_t \}_{t\geq 0} $.
    
   (ii)  Take an arbitrary $T_0>a$. Using the notations in \S \ref{S:2.2}, let  
       $$
   \wt \sG^{B, T_0}_t:=\sG^B_t \vee \sG^L_{T_0-a}  \quad \hbox{for } t\geq 0 
   \quad \hbox{ and } \quad 
     \wt \sF^{a, T_0}_t:= \wt \sG^{B, T_0}_{L_{(t-a)^+}}  \quad \hbox{for } t\in [0, T_0].
   $$
   We have $\sF^{a, \prime}_t \subset   \wt \sF^{a, T_0}_t \subset \wt \sF^a_t$ for $t\in [0, T_0]$. 
   Note that   for $ t \in [0, \infty)$, 
   $$
   t+R^a_t= a +S_{L_{(t-a)^+} }=\inf\{r>0: L_{(r-a)^+} > L_{(t-a)^+} \} 
   =\inf \big\{r>0: \< B_{L_{(\cdot-a)^+}} \>_r > \< B_{L_{(\cdot-a)^+}} \>_t \big\} .
      $$
   Thus $t+R^a_t$ is a finite $\{\sF^{a, \prime}_t\}_{t\geq 0}$-stopping time.
   Put $\tau:=T_0+R^a_{T_0}$.
Since $Y^a$ has the same law under $\QQ$ as $B_{L_{(\cdot-a)^+}}$ under $\PP$, and $\tau$ is a functional of $Y^a$, we have $\QQ(\tau<\infty)=\PP(\tau<\infty)=1$.
For $A\in\sF^{a,\prime}_\tau$ and $u>0$, we have $A\cap\{\tau\leq u\}\in\sF^{a,\prime}_u$, so $\QQ(A\cap\{\tau\leq u\})=0$ if and only if $\PP(A\cap\{\tau\leq u\})=0$.
Letting $u\to\infty$ shows that $\QQ$ and $\PP$ are equivalent on $\sF^{a, \prime}_\tau$.
 Note that  by \eqref{e:3.5a},
$$ \sG^B_{L_{(T_0-a)^+}} \subset  \sG_{L_{(T_0 -a)^+}}  = \sF^{a, \prime}_{T_0+R^a_{T_0}}  ,
$$
and $\sG^L_{T_0-a} \subset \sF^{a, \prime}_{T_0} \subset  \sF^{a, \prime}_{T_0+R^a_{T_0}}$. 
It follows that   $\wt \sF^{a, T_0}_t \subset \sF^{a, \prime}_{T_0+R^a_{T_0}} = \sF^{a, \prime}_{\tau}$.
 Hence $\QQ$ and $\PP$ are equivalent on $\wt \sF^{a, T_0}_{T_0}$. 
    Let $K_t$ be the strictly positive continuous $\{\wt \sF^{a, T_0}_t\}_{t\in [0, T_0]}$-martingale (continuity follows as in Corollary~\ref{C:3.5}) so that 
  $\frac{d \QQ}{d \PP}=K_t$ on $\wt \sF^{a, T_0}_t$  for $t\in [0, T_0]$.
    
    For an integer $n\geq 1$, define $T_n:=\inf\{t\geq 0: K_t \geq n\}\wedge T_0$  and $\Omega_n:=\{K_0<n\}\in \wt \sF^{a, T_0}_0$.
Then $\lim_{n\to \infty} T_n=T_0$ $\PP$-a.s., and $K_{  T_n}\1_{\Omega_n}$ is bounded and $\wt \sF^{a, T_0}_{T_0}$-measurable.
Since $K_0<\infty$, $\PP$-a.s. $ \Omega_n $ increases to $\Omega$   as $n\to\infty$. 
 By   \cite[Theorem 1.1]{ZC2}, 
 there exists an $\{\wt \sF^{a, T_0}_t\}_{t\in [0, T_0]}$--progressively measurable process   $\{  \xi^n_s;  s\in [0, T_0]\}$   
 having  $\EE \int_0^{T_0} |   \xi^n_s  |^2  dL_{(s-a)^+}   <\infty$
so that  
\begin{equation}\label{e:aa}
 \1_{\Omega_n} K_{T_n}  =  \1_{\Omega_n} \EE \big[ K_{T_n}  \mid \wt \sF^{a, T_0}_0 \big]  +  \int_0^{T_0}    \xi^n_s    \cdot d B_{L_{(s-a)^+} }  .
\end{equation}
 As $Y^a_t = Z_{L_{(t-a)^+}}$ is $\sF^{a, \prime}_t$-measurable, we have by (i) above that $Y^a$  
 is a local martingale with respect to the filtration    $\{ \wt \sF^{a, T_0}_t\}_{t\in [0, T_0]}$.
Now proceeding by the same argument as that for Theorem \ref{T:EMM} yields that 
$K_t=K_0 M_t$ for $t\in [0, T_0]$, where $M$ is the exponential martingale defined by \eqref{e:3.4}. 
  Then we deduce by the same proof as that for Theorem \ref{T:3.5} that $K_0=1$.
 This shows that $\QQ=\QQ^*$ on $\wt \sF^{a, T_0}_t =\sG^B_{L_{(t-a)^+}}  \vee \sG^L_{T_0-a}  $ for every $t\in [0, T_0]$ and every $T_0>a$. 
  Fixing  $t \geq 0$ first and taking $T_0\to \infty$ yields that $\QQ=\QQ^*$ on $\sG^B_{L_{(t-a)^+}}  \vee \sG^L_{\infty} = \wt \sF^a_t  $.
   Consequently, $\QQ=\QQ^*$ on $\wt \sF^a_\infty = \sF^{a, \prime}_\infty$ in view of \eqref{e:equal}. 
   This completes the proof of the theorem. 
 \end{proof}

    \section{Spot market}\label{S:4} 
  
  In this section, we introduce a spot market model that is suitable for bear markets but also contains 
  the classical spot market model as its special case, 
  and study the pricing  of contingent claims for this model.

We use the setup from  Section  \ref{S:2}; that is,    $B=\{B_t; t\geq 0\} = \big\{(B^{(1)}_t, \cdots ,  B^{(d)}_t)^{tr};  t\geq 0 \big\} $ 
is a $d$-dimensional standard Brownian motion,
    and $S$ a  subordinator that is strictly increasing and independent of $B$ with $S_0=0$, 
 and the processes $\wt X_t =(X_t, R_t)$ is defined as in \eqref{e:2.1}. 
 Here the  superscript {\it tr} stands for vector transpose. 
 We emphasize that   the subordinator $S$ can possibly have zero drift. 
  Recall also the filtrations  
 $\{\sF_t^{a, \prime}; t\geq 0\}$, $\{\sF^a_t; t\geq 0\}$ and    $\{\wt \sF^a_t; t\geq 0\}$ defined at the end of Section  \ref{S:2}.

\subsection{Sub-diffusion model for stock prices}\label{S:4.1}

   Consider a spot market that consists of a risk-free bond (or money market fund)
   and $d$ common stocks. 
           Denote by $A_t$ the price process of one share of the risk-free asset at time $t$ with $A_0=1$ and $r(t)$ the interest rate at time $t$,
       which is progressively measurable with respect to $\{\sF^{a,\prime}_t\}_{t\geq 0}$. Then 
\begin{eqnarray*}
d A_t  &=&   r (t)   A_t d  t,\quad t> 0
\end{eqnarray*}
 and $A_t=\exp( \int_0^t r(s) ds )$. 
Define the discount process by
 \begin{eqnarray}\label{D}
 A_t^{-1}  = \exp{\left(-\int_0^t r(s) ds \right)} ,
 \end{eqnarray} 
 which is the current value of one dollar at the future time $t$.

Denote by $\sS^{(i)}_t$ the  unit price of   the $i$th stock at time $t$, and 
  ${\mathcal S}_t:=(\sS_t^{(1)}, \cdots , \sS^{(d)}_t)^{tr}$.
     We model the stock price process $\sS$ by the following stochastic differential equation driven by sub-diffusion:
     for each $1\leq i \leq d$, 
  \begin{equation}\label{e:4.2} 
  d\sS^{(i)}_t = \sS^{(i)}_t \Big( r(t) dt +  \bar \mu_i (t) dL_{(t-a)^+}  +  \sum_{j=1}^d \sigma_{ij}(t) dB^{(j)}_{L_{(t-a)^+}}  \Big), \quad  t\geq 0 ,
  \end{equation}  
  where  $a\geq 0$ is the initial wake up time of the stock,  $\bar \mu (t)=(\bar \mu_1(t), \cdots, \bar\mu_d (t))^{tr} $ 
   and $\sigma (t)=(\sigma_{ij}(t))_{1\leq i, j\leq d}$ are all $\{\sF^{a, \prime}_t\}_{t\geq 0}$-progressively measurable.
   For simplicity, we assume throughout the rest of the paper that 
   \begin{equation} \label{e:con1}
 \hbox{the   processes $\{r(t); t\geq 0\}$, $\{\bar \mu (t); t\geq 0\}$ and $\{\sigma (t); t\geq 0\}$ are all bounded}
 \end{equation}
   and there is a constant $\lambda_\sigma>0$ so that $\PP$-a.s.,
   \begin{equation} \label{e:con2}
   \xi \cdot \sigma (t) \xi \geq \lambda_\sigma |\xi |^2 \quad \hbox{for every } t>0 \hbox{ and } \xi \in \RR^d .
   \end{equation} 
   Under these conditions, the system of  the SDEs \eqref{e:4.2} has a unique strong solution for any given initial value $\sS_0\in \RR^d$;
   see \eqref{e:4.5} below. 
     The discounted stock price process  $A_t^{-1} \sS^{(i)}_t$ for the $i$th stock  satisfies 
   \begin{eqnarray}\label{e:4.3} 
      d(A_t^{-1} {\sS^{(i)}_t} ) &=&   -r(t)A_t^{-1}{\sS^{(i)}_t} dt+ A_t^{-1} {\sS^{(i)}_t} 
  \Big( r(t)  dt +    \bar\mu_i (t)  dL_{(t-a)^+}+  \sum_{j=1}^d \sigma_{ij}(t)     dB^{(j)}_{L_{(t-a)^+}}  \Big)  \nonumber\\
  &=& A_t^{-1} {\sS^{(i)}_t}    \Big(     \bar\mu_i (t)  dL_{(t-a)^+}+  \sum_{j=1}^d \sigma_{ij}(t)     dB^{(j)}_{L_{(t-a)^+}}  \Big)  \nonumber\\ 
  &=& A_t^{-1} {\sS^{(i)}_t}  \sum_{j=1}^d \sigma_{ij}(t)    \Big(   \sum_{k=1}^d     \sigma^{jk}(t)  \bar\mu_k (t)  dL_{(t-a)^+} 
  +      dB^{(j)}_{L_{(t-a)^+}}  \Big)  \nonumber  \\
  &=& A_t^{-1} {\sS^{(i)}_t}  \sum_{j=1}^d \sigma_{ij}(t)   d \widetilde B^{(j)}_{L_{(t-a)^+}},
   \end{eqnarray}
   where $(\sigma^{ij}(t))$ is the inverse of the matrix $(\sigma_{ij}(t))$ (by \eqref{e:con2}, $|\sigma(t)\xi|\geq \lambda_\sigma|\xi|$, so $\sigma^{-1}$ is bounded by $1/\lambda_\sigma$) and 
  \begin{equation}\label{e:4.6n} 
   \widetilde B^{(j)}_{L_{(t-a)^+}} :=  B^{(j)}_{L_{(t-a)^+}}  + \int_0^t  \sum_{k=1}^d     \sigma^{jk}(s)  \bar \mu_k (s)  dL_{(s-a)^+} .
\end{equation}
   Thus for each $1\leq i\leq d$, 
  \begin{equation}\label{e:4.5}
   \sS^{(i)}_t =   \sS_0^{(i)}   A_t \exp \Big( \int_0^t  \sum_{j=1}^d \sigma_{ij}(s)   d \widetilde B^{(j)}_{L_{(s-a)^+}} 
   -\frac12 \int_0^t \sum_{j=1}^d \sigma_{ij}(s)^2 dL_{(s-a)^+} 
   \Big) .
   \end{equation}

\subsection{Sub-diffusion martingale equivalent   measure} \label{S:4.2}

\begin{definition} \label{D:4.4} \rm 
 Let $a\geq 0$.
 \begin{enumerate}
\item[(i)] We call a probability measure $\QQ$ on $(\Omega,   \sF^{a, \prime}_\infty)$  
{\it  an  equivalent martingale measure}  or  {\it  a risk-neutral probability measure} if 
$\QQ$  is  equivalent to $\PP$ on   $\sF^{a, \prime}_t$ for every $t\geq 0$ and 
$ A^{-1} \sS$ is a local martingale with respect to the filtration 
$\{ \sF^{a, \prime}_t\}_{t\geq 0}$
under $\QQ$. The   latter is equivalent to  $\{\widetilde B_{L_{(t-a)^+}} ; t\geq 0\}$ defined by \eqref{e:4.6n}   being an
 $\{ \sF^{a, \prime}_t\}_{t\geq 0}$-local martingale   under $\QQ$. 
We denote the totality of all equivalent martingale measures by  $\sC_1$.
 
 \item[(ii)] We call a probability measure $\QQ$ on $(\Omega, \wt \sF^a_\infty)$  
 a {\it  generalized  equivalent martingale measure}  or  a {\it  generalized risk-neutral probability measure}
 if  $\QQ$ t is   equivalent to $\PP$  on    $\wt \sF^a_t$ for every $t\geq 0$ and 
 $A^{-1} \sS$ is a local martingale with respect to the filtration $\{\wt \sF^a_t\}_{t\geq 0}$
under $\QQ$. The   latter is equivalent to  $\{\widetilde B_{L_{(t-a)^+}} ; t\geq 0\}$ defined by \eqref{e:4.6n}   being an
 $\{\wt \sF^a_t\}_{t\geq 0}$-local martingale   under $\QQ$. 
 We denote the totality of all generalized equivalent martingale measures by  $\sC_2$.

  \end{enumerate}
 \end{definition}

 Recall  that   $ \sF^{a, \prime}_\infty= \wt \sF^a_\infty$ by \eqref{e:equal}.

\begin{theorem} \label{T:4.1} 
For each $i\in\{1, 2\}$, there exists a   unique  probability measure $\QQ^* \in \sC_i$ so  that the discounted stock price process 
$$
 A^{-1} \sS :=  \big\{A_t^{-1} \sS_t=(A_t^{-1} \sS^{(1)}_t, \cdots, A_t^{-1} \sS^{(d)}_t)^{tr}; \  t\geq 0 \big\}
$$ 
is an  $\{{\wt \sF^a_t} \}_{ t\geq 0}$-martingale  
 driven by the sub-diffusion   $\wt B_{L_{(t-a)^+}}$ under $\QQ^*$
that has the same distribution as $B_{L_{(t-a)^+}}$ under $\PP$.   
These two probability measures $\QQ^*$ are the same on  $ \sF^{a, \prime}_\infty= \wt \sF^a_\infty$, independent of $i=1, 2$. 
\end{theorem}

\begin{proof} Let $\sigma$, $\bar\mu$ and $r$ satisfy \eqref{e:con1}--\eqref{e:con2}. Define 
\begin{equation}
\widehat \mu_j (t):=  \sum_{k=1}^d     \sigma^{jk}(t)  \bar \mu_k (t)  \quad \hbox{and} \quad 
\widehat \mu (t) = (\widehat \mu_1 (t), \cdots, \widehat \mu_d (t))^{tr}.
\end{equation}
Under assumptions \eqref{e:con1}-\eqref{e:con2}, $\{\wh \mu (t)\}$ is a bounded process so by 
  Theorem~\ref{G1} and Remark \ref{R:3.2},  
\begin{equation}\label{e:4.9}
M_t=\exp\left(  - \int_0^t  \widehat \mu (s)   \cdot  d B_{L_{(s-a)^+}}  - \frac 1 2  \int _0^t |\widehat \mu (s)  |^2      d L_{(s-a)^+}     \right), \quad t\geq 0.
\end{equation} 
  is a strictly positive continuous square-integrable martingale with respect to the filtration $\{{\wt \sF^a_t}\}_{t\geq 0}$ as well as 
   $\{{ \sF^{a, \prime}_t}\}_{t\geq 0}$.     
 It uniquely defines a   probability measure $\QQ^*$ on $\wt  \sF^{a}_\infty = \sF^{a, \prime}_\infty$  by
 \begin{equation}\label{e:4.10} 
    d \QQ^*=M_t d \PP   \quad  \hbox{on }\wt  \sF^{a}_t   \ \hbox{ for each } t\geq 0.   
 \end{equation} 
  Moreover by Theorem \ref{T:3.5} (respectively, Theorem \ref{T:uniqueness}),  
  $\QQ^*$ is the unique probability measure in $\sC_2$ (respectively, in $\sC_1$)  under which 
 \begin{equation}\label{e:4.11a} 
 \widetilde B_{L_{(t-a)^+}} =  B_{L_{(t-a)^+}}  + \int_0^t   \widehat \mu (s)  dL_{(s-a)^+}, \quad t\geq 0,
\end{equation} 
is a $d$-dimensional sub-diffusion that has the same distribution as $\{B_{L_{(t-a)^+}}; t\geq 0\}$ under $\PP$. Finally, by \eqref{e:4.14n} and the boundedness of $\sigma$, Novikov's criterion applies under $\QQ^*$.
Indeed, $\EE^{\QQ^*}\exp\big(\frac12\int_0^T|\sigma_{i\cdot}(s)|^2dL_{(s-a)^+}\big)\leq \EE\exp(c L_T)<\infty$ for some constant $c$ by Theorem~\ref{T:2.1}, because $L$ has the same law under $\QQ^*$ and $\PP$.
Hence each $\wh\sS^{(i)}$ is a $\{\wt\sF^a_t\}$-martingale under $\QQ^*$.
\end{proof}

Theorem \ref{T:4.1} in particular gives the existence of an equivalent martingale measure as well as a generalized  equivalent martingale measure.   Since the discounted stock price process 
$$
A^{-1} \sS :=  \big\{A_t^{-1} \sS_t=(A_t^{-1} \sS^{(1)}_t, \cdots, A_t^{-1} \sS^{(d)}_t)^{tr}; t\geq 0 \big\}
$$ 
 is continuous on $[0, \infty)$, it is  locally bounded on $[0, \infty)$.
Fix $T>0$. The measure $\QQ^*$ is equivalent to $\PP$ on $\sF^{a,\prime}_T$. Thus by \cite[Corollary 1.2]{DS} applied on $[0,T]$, the spot market satisfies ``no free lunch with vanishing risk" property in the sense of  \cite[Definition 2.8]{DS}. 
Hence the spot market admits no arbitrage opportunity in the above sense.

\begin{remark} \label{R:4.2} \rm 
In the sub-diffusion model, unless the subordinator $S_t  =\kappa t$ for some constant $\kappa >0$, 
 there are infinitely many generalized equivalent martingale measures and infinitely many  equivalent martingale measures for   $A^{-1} \sS$ 
 as shown Theorem   \ref{T:EMM} and  Theorem \ref{T:REMM}, respectively.  
However, by Theorem \ref{T:3.5} and Theorem \ref{T:uniqueness}, 
the probability measure $\QQ^*$ defined by \eqref{e:4.9}-\eqref{e:4.10} is the unique probability measure in $\sC_1$ (respectively, 
 in $\sC_2$)  that makes    $A^{-1} \sS$   an    $\{ \sF^{a, \prime}_t  \}_{ t\geq 0}$-local-martingale (respectively, 
  an    $\{ \wt \sF^a_t  \}_{ t\geq 0}$-local-martingale)
driven by a sub-diffusion.
We call this particular measure  $\QQ^*$  the {\it sub-diffusion equivalent martingale measure}. The notation  
 $\QQ^*$ from this point on specifically refers to this probability measure. 
  \qed
\end{remark}
 
 In the next two sub-sections, we give further properties of the sub-diffusion equivalent martingale measure $\QQ^*$
 that are relevant  for option pricing and hedging in the incomplete sub-diffusive markets.

 \subsection{Minimal martingale measure}

For simplicity, denote by $\wh \sS$ the discounted stock price process $ A^{-1} \sS$,
 and $J= (J^{(1)}, \cdots, J^{(d)})$ the $\PP$-martingale part of $\wh \sS$. By  \eqref{e:4.3} and \eqref{e:4.5},     
 for $1\leq i \leq d$, 
\begin{equation} \label{e:n4.12} 
J^{(i)}_t =\int_0^t  \wh \sS^{(i)}_r  \sum_{j=1}^d \sigma_{ij}(r)   d   B^{(j)}_{L_{(r-a)^+}}, \quad t\geq 0
\end{equation} 
and
\begin{equation}\label{e:4.14n}
 \wh   \sS^{(i)}_t =   \sS_0^{(i)}    \exp \Big( \int_0^t  \sum_{j=1}^d \sigma_{ij}(s)   d \widetilde B^{(j)}_{L_{(s-a)^+}} 
   -\frac12 \int_0^t \sum_{j=1}^d \sigma_{ij}(s)^2 dL_{(s-a)^+} 
   \Big) ,
   \end{equation} 
   with $\widetilde B^{(j)}$ defined by \eqref{e:4.6n}.
    The following definition is adapted from  that of \cite[Definition 3.2]{FS} to the multi-dimensional stock market case. 

\begin{definition} \rm 
A probability measure $\QQ$ on $\sF^{a,\prime}_\infty$ that is  equivalent to $\PP$ on each $\sF^{a,\prime}_t$ for $t\geq 0$
is said to be a  {\it minimal martingale measure (MMM)} if
\begin{enumerate}
\item[(i)] $\QQ=\PP$ on $\sF^{a, \prime}_{0}$ and $\wh  \sS$ is a $\QQ$-local martingale with respect to the filtration
$\{\sF^{a, \prime}_t\}_{t\geq 0}$; and

\item[(ii)] every square-integrable $\PP$-martingale $K$ with respect to the filtration $\{\sF^{a, \prime}_t\}_{t\geq 0}$
having $[ K, \wh \sS^{(i)} ] \equiv 0$   for every $1\leq i \leq d$ is a $\QQ$-martingale. 
\end{enumerate} 
\end{definition}

The following theorem can be established by utilizing \cite[Theorem 3.5]{FS}. Here we give a direct and self-contained proof.

\begin{theorem} \label{T:MMM} 
The sub-diffusion equivalent martingale measure $\QQ^*$  
is the unique minimal martingale measure.
\end{theorem}

\begin{proof} (i)  (Existence) 
 By Theorem \ref{T:4.1}, 
 $\wh \sS$ is a $ \QQ^*$-local martingale with respect to the filtration
$\{\sF^{a, \prime}_t\}_{t\geq 0}$. 
Suppose that $K$ is a  square-integrable $\PP$-martingale  with respect to the filtration $\{\sF^{a, \prime}_t\}_{t\geq 0}$
having $[ K, \wh \sS^{(i)} ]\equiv 0$ for every $1\leq i\leq d$. By \eqref{e:n4.12}, the latter is equivalent to 
$$
\int_0^t \wh \sS^{(i)}_r  \sum_{j=1}^d \sigma_{ij}(r)   d  \big\langle K,  B^{(j)}_{L_{(\cdot -a)^+}} \big\rangle _r =0, \quad t\geq 0,
$$
for every $1\leq i\leq d$. By the ellipticity of the volatility matrix \eqref{e:con2} and the strict positivity of  $\wh \sS^{(i)}_s$,
the above is equivalent to 
    \begin{equation} \label{e:4.13n}
\big\langle  K,  B^{(i)}_{L_{(\cdot -a)^+}} \big\rangle  \equiv 0 \quad \hbox{for every } 1\leq i\leq d. 
 \end{equation}
Consequently, $\langle  K, M\rangle  \equiv 0$, where $\{M_t; t\geq 0\}$ is the exponential martingale defined by \eqref{e:4.9}.
 By Girsanov theorem (see, e.g., \cite[Theorem 12.13]{HWY}), 
$K$ is an $\{\sF^{a, \prime}_t\}_{t\geq 0}$-local-martingale under $\wh \QQ^*$. 
As $K_t \in L^2(\PP)$ for every $t>0$, by Cauchy-Schwarz inequality and Doob's $L^2$-maximal inequality, we have for each $T>0$,  
$$
\EE^{\wh \QQ^* } \big[\sup_{t\in [0, T] }|K_t| \big] \le  \Big( \EE [ M_{T}^2] \, \EE \Big[ \sup_{t\in [0, T]} K_t^2 \Big] \Big)^{1/2}  
\leq 2  \Big( \EE [ M_{T}^2] \, \EE  [  K_T^2 ] \Big)^{1/2} <\infty. 
$$
Thus  $K$ is a $  \QQ^*$-martingale. This proves that $  \QQ^*$ is a minimal martingale measure. 
 
  (ii) (Uniqueness) Suppose $\QQ$ is another  minimal martingale measure. 
By the ellipticity assumption \eqref{e:con2}  and \eqref{e:4.3}
  $\wh \sS$ is a $\QQ$-local martingale with respect to the filtration
$\{\sF^{a, \prime}_t\}_{t\geq 0}$ if and only if   $\{  Y^a_t; t\geq 0\}$ of \eqref{e:4.1} 
is a $\QQ$-local martingale with respect to the filtration
$\{\sF^{a, \prime}_t\}_{t\geq 0}$. 
 By Theorem \ref{T:REMM},   there is an   
 $\{\sF^{a, \prime}_t\}_{ t\geq 0}$-local martingale $N$  with $N_0=0$, $\big\langle  N, B^{(i)}_{L_{(\cdot -a)^+}} \big\rangle  \equiv 0$
 for every $1\leq i\leq d$ and    $N_t-N_{t-} > -1$ for every $t>0$ so that 
 \begin{equation} \label{e:4.14}
 \frac{d \QQ}{d \PP} =  M_t  \sE (N)_t  =  \sE \Big( - \int_0^\cdot \widehat\mu (s) d B_{L_{(s-a)^+}} + N \Big)_t
 \quad \hbox{on each } \sF^{a, \prime}_t . 
\end{equation}
Let $K$ be an arbitrary  square-integrable $\PP$-martingale  with respect to the filtration $\{\sF^{a, \prime}_t\}_{t\geq 0}$
having $\big\langle  K, B^{(i)}_{L_{(\cdot -a)^+}} \big\rangle  \equiv 0$  for every $1\leq i \leq d$.
By the Girsanov theorem and \eqref{e:4.13n}, $K_t-\langle K,N\rangle _t$ is an $\{\sF^{a, \prime}_t\}_{t\geq 0}$-local martingale under $\QQ$.
Since $\QQ$ is a minimal martingale measure, $K$ is a $\QQ$-martingale.
Hence $\langle K,N\rangle $ is a predictable local martingale of finite variation under $\QQ$, so that $\langle K,N\rangle \equiv 0$.
We next show that $\langle K,N\rangle \equiv 0$ in fact holds for any square-integrable $\PP$-martingale $K$ with respect to the filtration $\{\sF^{a, \prime}_t\}_{t\geq 0}$.

(a) Note that $B_{L_{(t-a)^+}} $ is a square-integrable $\{\sF^{a, \prime}_t\}_{t\geq 0}$-martingale 
under $\PP$ with 
$$
\big\langle  B^{(i)}_{L_{(\cdot -a)^+}}, B^{(j)}_{L_{(\cdot -a)^+}} \big\rangle _t = 0
\quad \hbox{for all } i\not= j
$$
and $t\geq 0$.  Let  $K$ be an arbitrary square-integrable $\{\sF^{a, \prime}_t\}_{t\geq 0}$-martingale under $\PP$.
By Kunita-Watanabe decomposition, Lemma \ref{L:KW}, there 
an $\RR^d$-valued $\{\sF^{a, \prime}_t\}_{t\geq 0}$-predictable  process $\{\xi_s; s\geq 0\}$ which is unique
$dL_{(t-a)^+}\times \PP$-a.s. having $\EE \int_0^T |\xi_s|^2  d L_{(s-a)^+} <\infty$ for each $T>0$ and 
a rcll local  martingale $I$  having $\big\langle  I, B^{(i)}_{L_{(\cdot -a)^+}} \big\rangle _t = 0$ for every $t\geq 0$
and $1\leq i\leq d$ so that 
\begin{equation}\label{e:4.16} 
K_t = \int_0^t \xi_s \cdot d B_{L_{(s-a)^+}} + I_t \quad \hbox{for } t\geq 0.
\end{equation} 
 It follows that both $I$ and $\int_0^\cdot \xi_s \cdot d B_{L_{(s-a)^+}}$ are square-integrable martingales. 
 Since $\langle B^{(i)}_{L_{(\cdot -a)^+}},N\rangle \equiv 0$ and, by the preceding paragraph applied to $I$, $\langle I,N\rangle \equiv 0$, we have
$$
\langle K, N\rangle _t = \Big\langle  \int_0^\cdot  \xi_s \cdot d B_{L_{(s-a)^+}}, \,  N \Big\rangle _t + \langle I, N\rangle _t =0 \quad \hbox{for every  } t\geq 0.
$$

(b) Let $\{T_n\}_{n\geq 1}$ be an increasing sequence of $\{\sF^{a, \prime}_t\}_{t\geq 0}$-stopping times so that 
$t\mapsto N_{t\wedge T_n}$ is a uniformly  integrable $\PP$-martingale for every $n\geq 1$.  Fix $T>0$ and $n\geq 1$. 
Let 
$$
K_t = \begin{cases} 
\EE \big[  {\rm sgn} (N_{T\wedge T_n} ) | \sF^{a, \prime}_t \big]       \quad  &\hbox{for }  t\in [0, T]\\
   {\rm sgn} (N_{T\wedge T_n} )  &\hbox{for }  t >T.
 \end{cases}
 $$
  Then $K$ is a bounded $\{\sF^{a, \prime}_t\}_{t\geq 0}$-martingale under $\PP$ so we have $\langle K, N\rangle \equiv  0$
  by (a).
Since $KN-\langle K,N\rangle $ is a $\PP$-local martingale, $KN$ is a $\PP$-local martingale. As $ | K_{t\wedge T_n} N_{t\wedge T_n}| \leq |N_{t\wedge T_n}|$,
$t\to K_{t\wedge T_n} N_{t\wedge T_n}$ is uniformly integrable and hence is a  $\PP$-martingale. It follows that 
$$ 
\EE \left[  |N_{T\wedge T_n} | \right]  =
\EE \left[    K_{T\wedge T_n} N_{T\wedge T_n} \right] = \EE [ K_0 N_0]  =0.
 $$
 Thus $N_{T\wedge T_n}=0$ $\PP$-a.s.. Taking $n\to \infty$ yields $N_T=0$ $\PP$-a.s.. 
 This holds for every deterministic $T>0$. 
 So we have by \eqref{e:4.14} that  $\frac{d \QQ}{d \PP} =  M_t $  on each $ \sF^{a, \prime}_t$;
 that is, $\QQ = \QQ^*$. This completes the proof of the theorem.
\end{proof}  

\medskip

 For probability measures $\QQ \ll \PP$ on $\sF^{a, \prime}_t$,
 the relative entropy  (or Kullback-Leibler divergence)  of $\QQ$ with respect to $\PP$ on $\sF^{a, \prime}_t$ is 
   $$  \dKL (\QQ\Vert  \PP; \, \sF^{a, \prime}_t) := \EE^\PP \left[ \varrho_t\log \varrho_t \right], \qquad \varrho_t:=\frac{ d \QQ}{d \PP} \Big|_{\sF^{a, \prime}_t }.
 $$
Note that   the relative entropy is always nonnegative, and that $\dKL (\QQ\Vert  \PP; \, \sF^{a, \prime}_t) =0$
if and only if $\QQ =\PP$ on $\sF^{a, \prime}_t$. By \cite[Theorem 3.11]{FS},
among all  equivalent martingale measures $\QQ \in \sC_1$, the minimal martingale measure $\QQ^*$
is the unique equivalent martingale measure  that minimizes the functional
\begin{equation} \label{e:4.17}
 \dKL (\QQ\Vert  \PP; \, \sF^{a, \prime}_t) - \frac12 \EE^{\QQ}  \int_0^t | \widehat\mu_s|^2 dL_{(s-a)^+}  
 \quad \hbox{for every } t\geq 0,
 \end{equation} 
 where $\{\widehat\mu_t\}_{t\geq 0}$ is the bounded $\{\sF^{a, \prime}_t\}_{t\geq 0}$-progressively measurable process defined in the proof of Theorem~\ref{T:4.1}. 
In next subsection, we characterize $\QQ^*$ using reverse  entropy,  without an additional term.

   \subsection{Minimal reverse  entropy martingale measure}
 
 For probability measures $\PP \ll \QQ$ on $\wt \sF^a_t$, the relative entropy  (or Kullback-Leibler divergence)
 of $\PP$ with respect to $\QQ$ on $\wt \sF^a_t$ is 
$$  \dKL (\PP\Vert  \QQ; \, \wt \sF^a_t) := \EE^\PP \left[ \log \frac{ d \PP}{d \QQ} \Big|_{\wt \sF^a_t} \right]=\EE^\PP \left[ -\log \frac{ d \QQ}{d \PP} \Big|_{\wt \sF^a_t} \right].
 $$
  
 Recall that $\sC_1$ and $\sC_2$ denote the collection of all 
   equivalent martingale measures and of all generalized  equivalent martingale measures as defined in Definition \ref{D:4.4}.

 \begin{theorem} \label{T:4.6}
 Let $a\geq 0$. 
The sub-diffusion equivalent martingale measure $\QQ^*$ on $\sF^{a, \prime}_\infty = \wt \sF^a_\infty$
 is the unique minimizer of  the reverse relative entropy
$ \dKL (\PP\Vert  \QQ; \, \sF^{a, \prime}_t)$ for every  $t\geq 0 $ over $\QQ\in \sC_1$
as well as the unique minimizer of  the reverse relative entropy
$ \dKL (\PP\Vert  \QQ; \, \wt \sF^a_t)$  for every  $t\geq 0$ over $\QQ\in \sC_2$, that is, for every $t\geq 0$,

\begin{align*}
\dKL (\PP\Vert  \QQ^*; \, \sF^{a, \prime}_t) &= \inf_{\QQ \in \sC_1} \dKL (\PP\Vert  \QQ; \, \sF^{a, \prime}_t ),\\
\dKL (\PP\Vert  \QQ^*; \,  \wt \sF^a_t) &= \inf_{\QQ \in \sC_2} \dKL (\PP\Vert  \QQ; \, \wt \sF^a_t ).
\end{align*} 
\end{theorem}

\begin{proof}  Let   $M$ be the exponential martingale defined by \eqref{e:4.9}.

 (i) By Theorem \ref{T:REMM},  a probability measure $\QQ$ on $\sF^{a, \prime}_\infty$ is in  $\sC_1$ if and only if there is an   
 $\{\sF^{a, \prime}_t\}_{ t\geq 0}$-local-martingale $N$  with $N_0=0$, $\big\langle  N, B_{L_{(\cdot -a)^+}} \big\rangle  \equiv 0$
 and    $N_t-N_{t-} > -1$ for every $t>0$
so that $\{M_t  \, \sE (N)_t; t\geq 0\}$ is a $\{\sF^{a, \prime}_t\}_{ t\geq 0}$-martingale under $\PP$ and 
$$
\frac{d\QQ}{d \PP} \Big|_{\sF^{a, \prime}_t} = M_t  \, \sE (N)_t \quad \hbox{for } t\geq 0.
$$ 
   Thus for $t>0$,
\begin{eqnarray*} 
\dKL (\PP\Vert  \QQ; \, \sF^{a, \prime}_t ) &= & \EE^\PP  [-  \log  (  M_t  \sE (N)_t  )  ] =  
 \EE^\PP  [ - \log    M_t  ]  +   \EE^\PP  [ - \log     \sE (N)_t  ] \\
&=&    \dKL (\PP\Vert  \QQ^*; \, \sF^{a, \prime}_t)  +   \EE^\PP  [ - \log     \sE (N)_t  ] .
  \end{eqnarray*}
  Note that $-\log x$ is   convex and decreasing  on $(0, \infty)$ and $\sE (N)_t $ is a positive supermartingale. Thus by Jensen's inequality,
$$ \EE^\PP  [ - \log     \sE (N)_t  ] \geq  -\log (  \EE^\PP  [   \sE (N)_t  ])  \geq -\log (  \EE^\PP  [   \sE (N)_0  ] ) =0,
$$
with $ \EE^\PP  [ - \log     \sE (N)_t  ]=0$ if and only if $\sE (N)_t = \EE^\PP  [   \sE (N)_t  ]=1$.
Thus for every $t>0$ and $\QQ\in \sC_1$,  
$$
\dKL (\PP\Vert  \QQ; \, \sF^{a, \prime}_t)   \geq    \dKL (\PP\Vert  \QQ^*; \, \sF^{a, \prime}_t)
$$
with equality if and only if $\QQ=\QQ^*$ on $\sF^{a, \prime}_t$. 

\smallskip

  (ii) By    Theorem \ref{T:EMM},
 \begin{equation} \label{e:n4.10} 
 \sC_2 = \{ \QQ: d \QQ = K_0 d \QQ^* \hbox{ with }   K_0\in \sG^L_\infty, K_0>0 \hbox{ and } \EE [K_0] =1\},
 \end{equation}
 Let  $\QQ =K_0 \QQ^*\in \sC_2$.  Then for any $t>0$,  
\begin{eqnarray*} 
\dKL (\PP\Vert  \QQ; \, \wt \sF^a_t  ) &= & \EE^\PP  [-  \log  (K_0 M_t)  ] =  \EE^\PP [- \log K_0] 
+ \EE^\PP  [ - \log    M_t  ]  \\
&=&    \EE^\PP  [ - \log K_0] + \dKL (\PP\Vert  \QQ^*; \, \wt \sF^a_t )  .
  \end{eqnarray*}
Note that $-\log x$ is a convex function on $(0, \infty)$. By Jensen's inequality,
 $$
 \EE^\PP  [ - \log K_0] \geq \EE^\PP  [ - \log (\EE K_0)]  = 0
 $$ 
 with equality if and only if   $K_0=\EE K_0 =1$.
 Thus for every $t>0$ and $\QQ\in \sC_2$,  
$$
\dKL (\PP\Vert  \QQ; \, \wt \sF^a_t)   \geq    \dKL (\PP\Vert  \QQ^*; \, \wt \sF^a_t)
$$
with equality if and only if $\QQ=\QQ^*$ on $\wt \sF^a_\infty$.
\end{proof}

\section{Hedging and option pricing in sub-diffusive markets} \label{S:5} 

 In this section, we   study hedging  contingent claims and option pricing in sub-diffusive market.
In particular, we study  Black-Scholes model \eqref{BSsub} during bear markets and derive the price for its European call option explicitly. 
  We use the setup and notations at the beginning of Section \ref{S:4} including those in \S \ref{S:4.1}.
  
To avoid triviality, assume that $T>a\geq 0$.
   The European call option  gives the buyer the right to purchase 1 share of a stock, say stock $i$,  at time $T$ at a fixed price $K$.
So the value of this option at maturity is $({\sS^{(i)}_T}-K)^+$. One can consider a more general contingent claim $H$ which is 
an $\sF^{a, \prime}_T$-random variable and represents the payoff at maturity $T$.

The goal of this section is to study hedging and pricing of such a contingent claim under the sub-diffusive Black-Scholes model.
We then specialize it to the case when $d=1$ and $H=  \psi  (\sS_T)$, where $\psi $ is a non-negative function on $\RR$. 
This in particular covers the European call option with $\psi(x)=(x-K)^+$ and the European put option with $\psi  (x)= (K-x)^+$.

We will solve this problem by utilizing the sub-diffusion equivalent martingale measure $\QQ^*$ from Theorem \ref{T:4.1}.
Its   minimal martingale measure property from Theorem \ref{T:MMM} allows us to find an optimal replicating strategy that minimizes
the local risk and is mean self-financing.   
 To find the value or price process explicitly in $d=1$ case, 
we assume $r_t$, $\bar\mu_t$ and $\sigma_t$ are deterministic functions in $t\in [0, T]$.
 We emphasize that  the subordinator $S$ is quite general: the only condition imposed is that it is  strictly increasing. 
  In particular,     the subordinator $S$ can possibly have zero drift.

\subsection{Hedging and option pricing in incomplete markets} \label{S:5.1}

Recall that $A_t$   is the price of one share of the risk-free asset at time $t$ with $A_0=1$, that is, 
 $A_t= \exp (\int_0^t r_s ds)$. Suppose $T>a\geq 0$. 
As pointed out in Remark \ref{R:4.2}, when the subordinator $S$ is not deterministic, the spot market modeled 
   by sub-diffusions is arbitrage-free but incomplete as it may have infinitely many equivalent martingale measures under which
   the  discounted stock price process $\wh \sS_t= A_t^{-1}\sS_t$  is a martingale.  Consequently, not every European  option
    can be replicated by a self-financing strategy with  deterministic initial endowment.
     A typical claim will carry an intrinsic risk, and one tries to find a dynamic portfolio strategy
which reduces the actual risk.          We   show in this sub-section via F\"ollmer-Schweizer decomposition 
    that pricing using the sub-diffusion equivalent martingale measure $\QQ^*$
    minimizes  the local risk.
    
    Recall that by \eqref{e:4.3}, 
$$
d  \wh   \sS^{(i)}_t  =  \wh   \sS^{(i)}_t   \sum_{j=1}^d \sigma_{ij}(t)   dB^{(j)}_{L_{(t-a)^+}} 
 +    \wh   \sS^{(i)}_s \bar \mu_i (s) dL_{(s-a)^+}.
$$
So $  d  \langle    \wh   \sS^{(i)}\rangle _t=    (\wh   \sS^{(i)}_t )^2 \sum_{j=1}^d \sigma_{ij} (t)^2 d L_{(t-a)^+}$
and the drift term of $\wh   \sS^{(i)}_t$ is $\int_0^t   \wh   \sS^{(i)}_s \bar \mu_i (s) dL_{(s-a)^+}$.

\medskip

        Fix  constants $T>a \geq 0$.  Let 
        \begin{eqnarray}
\Theta &:=&  \Big\{  \xi =( \xi^{(1)}, \cdots, \xi^{(d)}): \hbox{   an $\RR^d$-valued $\{\sF^{a, \prime}_t\}_{t\in [0, T]}$-predictable 
process  having} \nonumber \\
&& \hskip 0.4 truein   \sum_{i=1}^d \EE \Big[ \int_{0}^{T} |\xi^{(i)}_r|^{2} d\langle \wh \sS^{(i)} \rangle_r 
+ \Big( \int_0^T    |\xi^{(i)}_s  \bar \mu_i (s)  |   \wh   \sS^{(i)}_s    dL_{(s-a)^+} \Big)^2 \Big]  <\infty \Big\}  .
\label{e:5.1a}
\end{eqnarray}
 
The integrability condition in \eqref{e:5.1a} is to ensure that $\int_0^t \xi_r \cdot d \wh  \sS_r$ is a square-integrable $\PP$-semimartingale
for $t\in [0, T]$. 

\begin{remark}\label{R:5.1}  \rm 
 Suppose that $\xi =( \xi^{(1)}, \cdots, \xi^{(d)})$ is  an $\RR^d$-valued $\{\sF^{a, \prime}_t\}_{t\in [0, T]}$-predictable 
process  having   
$\sum_{i=1}^d \int_{0}^{T} |\xi^{(i)}_r|^{2}\,d\langle \wh \sS^{(i)}\rangle_r  \in L^r (\sF^{a, \prime}_T, \PP)$
 for some $r>1$. Then by  \eqref{e:con1}-\eqref{e:con2}, Cauchy-Schwarz inequality,  H\"older's inequality and Theorem \ref{T:2.1},
 for each $1\leq i \leq d$, 
 \begin{eqnarray}
 &&   \EE \Big[  \Big( \int_0^T |\xi^{(i)}_s    \bar \mu_i (s) |    \wh   \sS^{(i)}_s   dL_{(s-a)^+}  \Big)^2 \Big]
    \leq      c \EE \Big[  L_{(T-a)^+}  \int_0^T |\xi^{(i)}_s |^2   (\wh   \sS^{(i)}_s)^2  dL_{(s-a)^+}   \Big]  \nonumber  \\
&\leq &  c \Big(\EE \Big[  \Big( \int_0^T |\xi_s^{(i)}|^2  d  \langle  \wh   \sS^{(i)} \rangle _s \Big)^r  \Big] \Big)^{1/r}
\left( \EE \left[  \left(  L_{(T-a)^+}  \right)^{r'}  \right] \right)^{1/r'}
  <\infty,
 \end{eqnarray}
 where $r'>1$ is the conjugate of $r$, that is, $1/r +1/r' =1$.  Thus $\xi \in \Theta$. 
 \qed
\end{remark}

     \begin{definition}  \rm Let $T>a\geq 0$. 
  An  $H\in L^{2}(\sF^{a, \prime}_T, \PP)$ has a {\it F\"ollmer-Schweizer   decomposition} if
\begin{equation}\label{e:FS} 
  H=H_{0}+\int_{0}^{T}\xi^{H}_t \cdot d \wh \sS_t+K^{H}_{T},
\end{equation} 
with $H_{0}\in\RR$, $\xi^{H}\in\Theta$, and $K^{H}$ a square-integrable $\PP$-martingale with respect to
the filtration $\{\sF^{a, \prime}_t\}_{t\in [0, T]}$ having 
$K^{H}_{0}=0$ and  $\langle K^{H}, \wh \sS^{(i)} \rangle\equiv 0$  for every $1\leq i\leq d$. 
\end{definition}
    
\begin{remark} \label{R:5.3} \rm
\begin{enumerate}
 \item[(i)] Recall that  $\QQ^*$ is the sub-diffusion equivalent martingale 
measure    defined via the exponential martingale $M$ through  \eqref{e:4.9}-\eqref{e:4.10}.
 For the F\"ollmer-Schweizer   decomposition \eqref{e:FS}, both  $\int_{0}^{\cdot}\xi^{H}_t\,d \wh \sS_t$
 and $K^H$ are  uniformly integrable $\{\sF^{a, \prime}_t\}_{t\in [0, T]}$-martingales under $\QQ^*$. 
This is because   by the ellipticity of the volatility matrix \eqref{e:con2} and the strict positivity of $  \wh \sS^{(i)}_s$,
$\langle K^{H}, \wh \sS^{(i)} \rangle\equiv 0$  for every $1\leq i\leq d$ is equivalent to 
    \begin{equation} \label{e:5.2}
\big\langle  K^H,  B^{(i)}_{L_{(\cdot -a)^+}} \big\rangle  \equiv 0 \quad \hbox{for every } 1\leq i\leq d. 
 \end{equation}
Thus by Girsanov theorem, $K^H$ is a local martingale under $\QQ^*$.  
On the other hand, by Cauchy-Schwarz inequality, Doob's maximal inequality,  Theorem \ref{T:2.1} and Remark \ref{R:3.2}, 
 \begin{eqnarray*}
 \EE^{\QQ^*}\Bigl[ \sup_{t\in [0, T]}|K^H_t| \Bigr]
 &=&  \EE \Bigl[ M_T \sup_{t\in [0, T]}|K^H_t|  \Bigr] 
 \leq  \Big( \EE \left[ M_T^2  \right]  \EE \Bigl[   \sup_{t\in [0, T]}|K^H_t|^2  \Bigr] \Big)^{1/2} \\
& \leq & 2 \Big( \EE \left[ M_T^2  \right]  \EE \Bigl[   |K^H_T|^2  \Bigr] \Big)^{1/2}
<\infty
\end{eqnarray*}
This establishes that $K^H$ is a uniformly integrable $\{\sF^{a, \prime}_t\}_{t\in [0, T]}$-martingale under $\QQ^*$. 
Furthermore, by Theorem \ref{T:4.1}, \eqref{e:4.14n} and Theorem \ref{T:2.1},
  $\wh \sS$ is a continuous square-integrable  $\{\sF^{a, \prime}_t\}_{t\in [0, T]}$-martingale  under $\QQ^*$.
Thus by the Burkholder-Davis-Gundy inequality and \eqref{e:con1}-\eqref{e:con2}, 
\begin{eqnarray*}
 \EE^{\QQ^*}\Bigl[ \sup_{t\in [0, T]} \Big| \int_0^t \xi^{H}_u   \cdot  d \wh \sS_u \Big|  \Bigr]
 &\leq & c  \EE^{\QQ^*}  \Bigl[  \Big( \int_0^T  |\xi^{H}_u|^2 d \langle \wh \sS \rangle _u\Big)^{1/2}  \Bigr] \\
 &=& c \EE  \Bigl[ M_T \Big( \int_0^T |\xi^{H}_u|^2 d \langle \wh \sS \rangle _u\Big)^{1/2}  \Bigr]  \\
& \leq &  c  \Big( \EE \left[ M_T^2  \right]  \EE \Bigl[     \int_0^T  |\xi^{H}_u|^2 d \langle \wh \sS \rangle _u  \Bigr] \Big)^{1/2} < \infty.
  \end{eqnarray*}
Thus $\int_{0}^{\cdot}\xi^{H}_t\,d \wh \sS_t$
 is a  continuous uniformly integrable $\{\sF^{a, \prime}_t\}_{t\in [0, T]}$-martingale under $\QQ^*$. 
 
\item[(ii)] The F\"ollmer-Schweizer   decomposition of $H\in L^{2}(\sF^{a, \prime}_T, \PP)$, if exists,  is unique in the sense that if
there are some $\wt H_{0}\in\RR$, $ \wt \xi^H \in\Theta$, and   a square-integrable $\PP$-martingale $\wt  K^H$ with respect to
the filtration $\{\sF^{a, \prime}_t\}_{t\in [0, T]}$ having 
 $ \wt K^H_{0}=0$ and  $\langle \wt K^H, \wh \sS \rangle\equiv 0$ so that  
$$
 H=\wt H_{0}+\int_{0}^{T}\wt \xi^{H}_t \cdot d \wh \sS_t+ \wt K^{H}_{T},
 $$
 then $\PP$-a.s., 
 $$
 \wt H_0=H_0, \quad \wt K^H_t=K^H_t \hbox{  for all } t\in [0, T] \quad \hbox{and} \quad  
 \int_0^T   | \wt \xi^{H}_u -  \xi^H_u |^2 d L_{(u-a)^+} =0.
 $$ 
    Indeed, taking a difference with \eqref{e:FS} yields 
 $$
H_{0}-\wt H_0 +\int_{0}^{T} ( \xi^{H}_t - \wt \xi^H_t) \,d \wh \sS_t = \wt K^H_T - K^{H}_T,
$$
 Taking conditional expectation on $\sF^{a, \prime}_t$ under $\QQ^*$, we have by (i), 
\begin{equation}\label{e:5.4}
 H_{0}-\wt H_0 +\int_{0}^{t} ( \xi^{H}_u - \wt \xi^H_u) \cdot d \wh \sS_u = \wt K^H_t - K^{H}_t
 \quad \hbox{for every } t\in [0, T].
\end{equation} 
Since $\wt K^H_t - K^{H}_t$ is a uniformly integrable $\{\sF^{a, \prime}_t\}_{t\in [0, T]}$-martingale under $\QQ^*$
with $\wt K^H_0-K^H_0=0$ and $[ \wh \sS, \wt K^H - K^{H} ] \equiv 0$, both sides of \eqref{e:5.4} are identically zero. 
Hence we have $\QQ^*$-a.s. (which is equivalent to $\PP$-a.s.), 
$\wt H_{0}=H_0$, $\wt K^H_t = K^{H}_t$ for all $t\in [0, T]$ and 
$\big\langle  \int_{0}^{\cdot}   ( \xi^{H}_u - \wt \xi^H_u) \cdot d \wh \sS_u \big\rangle \equiv 0$.
The latter is equivalent to $\int_0^T   | \xi^{H}_u - \wt \xi^H_u  |^2 d L_{(u-a)^+} =0$ by
\eqref{e:4.3}   and \eqref{e:con1}-\eqref{e:con2}. 
This establishes the claim.
\end{enumerate} 
\qed
\end{remark}

\smallskip

The following theorem gives the F\"ollmer-Schweizer   decomposition via the sub-diffusion equivalent martingale measure $\QQ^*$. 
Typically such a property is   derived  from the minimal martingale property of $\QQ^*$ established in Theorem \ref{T:MMM};
see \cite{FS}. Here we establish it directly by using Kunita-Watanabe decomposition and Girsanov theorem, 
without invoking Theorem \ref{T:MMM}.

\begin{theorem} \label{T:5.2}
Let $T>a\geq 0$. 
Every  $H\in L^p ( \sF^{a, \prime}_T, \PP)$ for some $p>2$
 admits a F\"ollmer-Schweizer decomposition  \eqref{e:FS} with $H_0=\EE^{\QQ^*} [H]$. 
Moreover, $\EE^{\QQ^*} [H^2] <\infty$ and \eqref{e:FS} is exactly the Kunita-Watanabe decomposition
in Lemma \ref{L:KW} 
for the square-integrable martingale $\{\EE^{\QQ^*} [ H \mid \sF^{a, \prime}_t]; t\in [0, T]\}$ with respect to $\wh \sS$ 
under $\QQ^*$. 
 \end{theorem}

\begin{proof}  Recall that $M$ is the exponential martingale defined by \eqref{e:4.9}. 
Fix $q\in(2,p)$. By Remark \ref{R:3.2}(ii) and H\"{o}lder's inequality,   
\begin{equation}\label{e:Lq}
  \EE^{\QQ^*}[|H|^{q}]=\EE [M_{T}|H|^{q}]\le
  ( \EE [ M^{p/(p-q)}_T])^{(p-q)/p} \, ( \EE [ |H|^p] )^{q/p} <\infty.
 \end{equation}
   Let $H_t$ be the rcll version of the square-integrable martingale $\EE^{\QQ^*}[H \mid \sF^{a, \prime}_t]$. 
    The discounted stock price process $\wh \sS$ is a continuous square-integrable martingale under $\QQ^*$.
  So by Kunita-Watanabe decomposition, Lemma \ref{L:KW}, 
   there is an  $\RR^d$-valued $\{\sF^{a, \prime}_t\}_{t\in [0, T]}$-predictable
   process $\{\xi^H_r\}_{r\in [0, T]}$ having $\EE^{\QQ^*} \int_0^T |\xi_r^H|^2 d   \langle  \wh \sS\rangle _r <\infty$
  and a square-integrable $\QQ^*$-martingale $K^H$ so that  
  $$
  H_t = H_0 + \int_0^t  \xi^H_r \cdot d \wh \sS_r + K^H_t \quad \hbox{for } t\in [0, T]
  $$
  and $ [ K^H, \wh \sS^{(i)} ] \equiv  0$ for all $1\leq i\leq d$.
  Consequently, $[H]_t = \big[   \int_0^\cdot  \xi^H_r  \cdot d \wh \sS_r \big]_t + [ K^H ]_t$. 

    By Burkholder-Davis-Gundy inequality and Doob's maximal inequality,
\[
  \EE^{\QQ^*}\Bigl[\sup_{t\in [0, T]}|K^H_t|^{q}\Bigr]+\EE^{\QQ^*}\Bigl[(\textstyle\int_{0}^{T} |\xi^H_r|^{2}\,d \langle  \wh \sS \rangle _r )^{q/2}\Bigr]
  \le c_1 \,\EE^{\QQ^*}\bigl[[H]_{T}^{q/2}\bigr]\le c_2 \,\EE^{\QQ^*}[|H|^{q}]<\infty.
\]

 Note that 
 \begin{eqnarray} \label{e:5.3}
 1/ M_t  &= &  \exp\left(    \int_0^t  \widehat \mu (s)   \cdot  d B_{L_{(s-a)^+}}  +\frac 1 2  \int _0^t |\widehat \mu (s)  |^2      d L_{(s-a)^+} \right) 
 \nonumber \\
 &=& \exp\left(    \int_0^t  \widehat \mu (s)   \cdot  d \wt B_{L_{(s-a)^+}}  - \frac 1 2  \int _0^t |\widehat \mu (s)  |^2      d L_{(s-a)^+} \right) ,
 \end{eqnarray} 
 where $ \wt B_{L_{(s-a)^+}}$ defined by  \eqref{e:4.11a} is a sub-diffusion under $\QQ^*$. 
 Hence by Theorem \ref{G1} and Remark \ref{R:3.2}, $1/M_t$ is a  $\QQ^*$-martingale having finite moments of any order. 
Take $r\in (2, q)$. By H\"older's inequality with  $s=q/r>1$,
\begin{eqnarray}
  \EE^{\PP}\Bigl[ \sup_{t\in [0, T]}|K^H_t|^{r} \Bigr]
 &=& \EE^{\QQ^*}\Bigl[ (1/M_T)  \sup_{t\in [0, T]}|K^H_t|^{r} \Bigr]   \nonumber \\
 & \le& \Big( \EE^{\QQ^*}\Bigl[ (1/M_T)^{q/(q-r)}  \Bigr] \Big)^{(q-r)/q} 
 \Big( \EE^{\QQ^*}\Bigl[  \sup_{t\in [0, T]}|K^H_t|^{q}  \Bigr] \Big)^{r/q}  < \infty.  \label{e:5.6n}
\end{eqnarray} 
Similarly, we have 
$\EE^{\PP}[ (\int_{0}^{T} |\xi^H_s|^2 \,d \langle  \wh \sS \rangle _s )^r ]<\infty$. 
 So $\xi^{H}\in\Theta$ by Remark \ref{R:5.1}.

  As $\QQ^*$ is the minimal  martingale measure by Theorem \ref{T:MMM} and 
 $ [ K^H, \wh \sS^{(i)} ] \equiv  0$ for all $1\leq i\leq d$, we can readily conclude that 
 $K^H$ is a $\PP$-martingale. However, we present a direct proof of this fact without
 using Theorem  \ref{T:MMM}  as follows. 
 Observe that by   \eqref{e:4.11a} and \eqref{e:5.2}
$$
\big\langle  K^H,  \wt B^{(i)}_{L_{(\cdot -a)^+}} \big\rangle  = \big\langle  K^H,  B^{(i)}_{L_{(\cdot -a)^+}} \big\rangle  \equiv 0 \quad \hbox{for every } 1\leq i\leq d. 
 $$
 As $\frac{d \PP}{d \QQ^*} = 1/M_t$ on each $\sF^{a, \prime}_t$,
 we have  by Girsanov theorem and \eqref{e:5.3}, 
 $K^H$ is a $\PP$-local martingale with $[K^H, \wh \sS] \equiv 0$. Since  $\EE \left[ \sup_{t\in [0, T]} |K^H_t|^2 \right] <\infty$,  
  $K^H$ is a square-integrable $\PP$-martingale.
\end{proof} 

\medskip

Given $T>a\geq 0$, an {\it admissible  strategy} is a pair   of $\RR^d$-valued stochastic process  $\xi= (\xi^{(1)}, \cdots, \xi^{(d)}) \in \Theta$ and 
a real-valued $\{\sF^{a, \prime}_t\}_{t\in [0, T]}$-progressive measurable process $\eta$
 so that 
\begin{equation}\label{e:5.5}
 \wh V_t=\xi_t \cdot  \wh \sS_t+\eta_t,  \qquad \wh C_t= \wh V_t-\int_{0}^t\xi_u \cdot d \wh \sS_u.
 \end{equation} 
are right continuous and square-integrable. 
Here  $\xi^{(i)}_t$ stands for the number of shares invested in  the $i$th stock at time $t$ and  
 and $\eta_t$ is  the number of shares of the risk-free asset purchased at time $t$ with outside money.
 The values of $\xi^{(i)}_t$ and $\eta_t$ can be negative. 
 Thus $\wh V_t$ and $\wh C_t$ are the discounted value (that is, the value at the initial time)
 of the portfolio and discounted value  of the {\it cumulative cost} at time $t$.
 In other words, $\wh C_t$ is the discounted value of all the cash injected up to time $t$.
 A   strategy is called {\it self-financing}  if and only if the 
 cumulative cost process is constant, that is, $\wh C_t=\wh C_0$ for all $t\in [0, T]$.
 Note that   $\wh V_0 =\wh C_0$.
 Given a contingent claim $H \in L^2(\sF^{a, \prime}_T, \PP)$, we like to find
 an admissible strategy that  replicates the claim $H$ at the maturity $T$, that is, 
  $A_T \wh V_{T}=H$.  
   In an incomplete market, not every claim can be replicated by a self-financing strategy with
   {\it deterministic initial endowment}.
  We measure the remaining risk of the strategy at $t$ by
\begin{equation}\label{e:5.6}
  R_t=\EE\bigl[(\wh C_{T}-\wh C_t)^{2}\mid\sF^{a, \prime}_t\bigr]
\end{equation}
and   look for an admissible strategy that minimizes $R_t$ for  each $t\geq 0$,
 over all admissible continuations of this strategies from time $t$ on.  

An admissible strategy is called {\it locally risk-minimizing} if, for any $t \in [0, T)$, 
the remaining risk $R_t$ is minimal under all
infinitesimal perturbations of the strategy at time $t$; see \cite[Definition on p.349]{Sch} for a precise definition. 
 It is shown in \cite[Proposition 2.3]{Sch} to be  
equivalent to the following property of the associated discounted accumulative cost process
\begin{equation}\label{e:lrm}
\wh C \hbox{ is a square-integrable $\{\sF^{a, \prime}_t\}_{t\in [0, T]}$-martingale  under $\PP$ with } 
[\wh C, \wh \sS] \equiv 0.
\end{equation}

\smallskip

Following \cite[Definition 2.23]{FS}, we call an admissible strategy $(\xi, \eta)$ {\it optimal} 
if its associated discounted accumulative cost process $\wh C$ has the property \eqref{e:lrm}.
An optimal  strategy is   {\it mean-self-financing} in the sense that 
$\EE [ \wh C_T- \wh C_t \mid \sF^{a, \prime}_t] =0$ for every $t\in [0, T]$.

\smallskip

Suppose $(\xi, \eta)$ is an optimal replicating strategy for a contingent claim $H\in L^2(\sF^{a, \prime}_T, \PP)$. Then by 
\eqref{e:lrm},
\begin{equation}\label{e:5.9n}
A^{-1}_T  H= \wh V_T = \wh C_0 + \int_0^T \xi_u \cdot d \wh \sS_u + (\wh C_T -\wh C_0),
\end{equation}
with $\wh C -\wh C_0$ being a square-integrable $\PP$-martingale with zero initial values
and $[ \wh C -\wh C_0, \wh \sS ] \equiv 0$. This is exactly the F\"ollmer-Schweizer decomposition of $A_T^{-1}H$.
In fact, it  is shown in \cite[Proposition 2.24]{FS} that  the existence of an optimal strategy 
for a contingent claim $H\in L^2(\sF^{a, \prime}_T, \PP)$ 
is equivalent to the  F\"ollmer-Schweizer decomposition of $A_T^{-1}H$: 
\begin{equation}\label{e:5.7}
A^{-1}_T  H=H_0 + \int_{0}^{T}\xi^{H}_t\cdot d \wh \sS_t+K^{H}_{T},
\end{equation} 
with  $H_0\in \RR$,   $\xi^{H}\in\Theta$, and $K^{H}$ a square-integrable $\{\sF^{a, \prime}_t\}_{t\in [0, T]}$-martingale under  
both $\PP$  and $\QQ^*$  having  $K^{H}_{0}=0$ and  $\langle K^{H}, \wh \sS \rangle\equiv 0$.
For such a decomposition,  if it exists, the associated optimal strategy $(\xi , \eta )$ is given by
\begin{equation} \label{e:optim}
  \xi_t =\xi^{H}_t,
   \qquad
  \eta_t= \eta^H_t:= \EE^{\QQ^*}[ A^{-1}_T H \mid \sF^{a, \prime}_t] -\xi^{H}_t  \cdot \wh \sS_t,
\end{equation} 
Indeed, by the uniqueness of F\"ollmer-Schweizer decomposition (Remark \ref{R:5.3}(ii)), we have by \eqref{e:5.9n} and \eqref{e:5.7} that
$$
\wh C_0= H_0=\EE^{\QQ^*}[ A^{-1}_T H] , \quad \xi = \xi^H , \quad \hbox{and} \quad 
\wh C -\wh C_0 = K^H.
$$
Thus
\begin{eqnarray*}
\eta_t &=& \wh V_t - \xi^H_t \cdot \wh \sS_t= \wh C_t + \int_0^t \xi^H_u \cdot \wh \sS_u - \xi^H_u \cdot \wh \sS_t
=\wh C_0 +K^H_t +  \int_0^t \xi^H_u \cdot \wh \sS_u - \xi^H_u \cdot \wh \sS_t \\
&=& \EE^{\QQ^*}[ A^{-1}_T H \mid \sF^{a, \prime}_t] -\xi^{H}_t  \cdot \wh \sS_t. 
\end{eqnarray*}

\begin{prop} \label{P:5.5}
Every contingent claim $H\in L^p ( \sF^{a, \prime}_T, \PP)$ for some $p>2$
has an optimal replicating strategy $(\xi^H, \eta^H)$. The optimal replicating strategy is unique in the sense 
that if $(\wt \xi^H, \wt \eta^H)$ is another optimal replicating strategy, then $\PP$-a.s.
$$
 \int_0^T   | \wt \xi^{H}_u -  \xi^H_u |^2 d L_{(u-a)^+} =0
 \quad \hbox{and} \quad     \wt \eta^H _t=\eta ^H_t \hbox{  for all } t\in [0, T] . 
 $$ 
 The discounted value process for this optimal replicating strategy is 
 \begin{equation}\label{e:5.12} 
 \wh  V_t = \EE^{\QQ^*} \big[ A^{-1}_T H  \mid \sF^{a, \prime}_t \big] , \quad t\in [0, T].
 \end{equation}
In particular,  the initial value  for this optimal replicating strategy is $\EE^{\QQ^*} [ A^{-1}_T H ] $. 
\end{prop}

\begin{proof} Suppose $H\in L^p ( \sF^{a, \prime}_T, \PP)$ for some $p>2$. 

\smallskip

(i) (Existence) 
By Theorem \ref{T:5.2}, $A_T^{-1}H \in L^p ( \sF^{a, \prime}_T, \PP)$ has 
the F\"ollmer-Schweizer decomposition \eqref{e:5.7}  
with   $H_0=\EE^{\QQ^*} [ A^{-1}_T H]$,  $\xi^{H}\in\Theta$, and $K^{H}$ a square-integrable $\{\sF^{a, \prime}_t\}_{t\in [0, T]}$-martingale under both
$\PP$ and $\QQ^*$  having  $K^{H}_{0}=0$ and  $\langle K^{H}, \wh \sS \rangle\equiv 0$.
Let 
$$
\eta^H_t:= \EE^{\QQ^*}[ A^{-1}_T H\mid\sF^{a, \prime}_t] -\xi^{H}_t  \cdot \wh \sS_t \quad \hbox{for } t\in [0, T].
$$
As we noted in the above via  \cite[Proposition 2.24]{FS}, 
$(\xi^H, \eta^H)$ is an optimal replicating strategy for the contingent claim $H$.
Its discounted value process  for this optimal replicating strategy is 
$$
\wh V_t= \xi^H_t \cdot \wh \sS_t + \eta^H_t  
= \EE^{\QQ^*} \big[ A^{-1}_T H \mid\sF^{a, \prime}_t \big] , \quad t\in [0, T].
$$
   
\smallskip

(ii) (Uniqueness) Suppose that $(\wt \xi^H, \wt \eta^H)$ is an optimal replicating strategy for $H$.
Then its associated discounted accumulative cost process
$$
\wh C_t= \wt \xi^H_t \cdot  \wh \sS_t+ \wt \eta^H_t -\int_{0}^t \wt \xi^H_u \cdot d \wh \sS_u, \quad t\in [0, T],
$$
is a   a square-integrable $\{\sF^{a, \prime}_t\}_{t\in [0, T]}$-martingale  under $\PP$ with 
$[\wh C, \wh \sS] \equiv 0.$ Taking $t=T$ yields
$$
A^{-1}_T H = \wt \xi^H_T \cdot \wh \sS_T + \wt \eta^H_T =\wh C_0+  \int_{0}^T \wt \xi^H_u \cdot d \wh \sS_u + ( \wh C_T-\wh C_0). 
$$
The above is the   F\"ollmer-Schweizer decomposition for $A^{-1}_T H$.
Thus by the uniqueness of the  F\"ollmer-Schweizer decomposition from Remark \ref{R:5.3}(ii), 
we have $\PP$-a.s., 
\begin{equation}\label{e:5.11}
\wh C_0= H_0=\EE^{\QQ^*} [ A^{-1}_T H], \qquad \wh C_t -\wh C_0 = K^H_t \hbox{  for every } t\in [0, T],
\end{equation}
and $\int_0^T |\wt \xi^H_t - \xi^{H}_t|^2 dL_{(t-a)^+} =0$. 
Consequently, 
 \begin{eqnarray*}
 \wt \eta^H_t &=&   \wh C_0 +  (\wh C_t -\wh C_0) -  \wt \xi^H_t \cdot  \wh \sS_t+ \int_{0}^t \wt \xi^H_u \cdot d \wh \sS_u \\
& =&  \EE^{\QQ^*} [ A^{-1}_T H] +  K^H_t   + \int_{0}^t  \xi^H_u \cdot d \wh \sS_u  -  \xi^H_t \cdot  \wh \sS_t \\
&=& \EE^{\QQ^*}[ A^{-1}_T H\mid\sF^{a, \prime}_t] -\xi^{H}_t  \cdot \wh \sS_t \, = \, \eta^H_t
  \qquad \hbox{for every } t\in [0, T],
\end{eqnarray*}
This establishes the uniqueness of the optimal replicating strategy for $H$.
\end{proof}

\medskip

\begin{theorem}[Risk-minimization under $\QQ^*$] \label{thm:riskmin}
Suppose $T>a\geq 0$
and  $H\in L^p ( \sF^{a, \prime}_T, \PP)$ for some $p>2$.
Let  $\xi^H$ be the process in the  F\"ollmer-Schweizer decomposition \eqref{e:5.7} of $A_T^{-1}H$.  
 Then  among replicating strategies $(\xi, \eta)$  with
$\EE^{\QQ^*}[\int_{0}^{T} |\xi_u |^{2}\,d \langle  \wh \sS \rangle _u ]<\infty$ and $\wh V_t\in L^{2}(\QQ^*)$, the strategy
\begin{equation} 
  \xi_t =\xi^{H}_t,
   \qquad
  \eta_t= \eta^H_t:= \EE^{\QQ^*}[ A^{-1}_T H \mid \sF^{a, \prime}_t] -\xi^{H}_t  \cdot \wh \sS_t,
\end{equation} 
minimizes the risk $R^*_t=\EE^{\QQ^*} \big[ (\wh C_{T}- \wh C_t)^{2} \mid\sF^{a, \prime}_t  \big]$ simultaneously for all $t\in [0, T]$. 
 The optimal strategy is unique in the sense that if $(\wt \xi, \wt \eta)$ is an admissible replicating strategy that 
minimizes the risk $R^*_t $ simultaneously for all $t\in [0, T]$, then
$$
\PP \Big(  \int_0^{T} | \wt \xi_u-   \xi^{H}_u  |^2 d L_{(u -a)^+} =0 \hbox{ and } \wt \eta_t = \eta^H_t \hbox{ for every } t\in [0, T] \Big)=1.
$$
\end{theorem}

\begin{proof}    For any admissible replicating strategy $(\xi, \eta)$, by \eqref{e:5.5}
$$
\wh C_{T}- \wh C_t= \wh V_T-\wh V_t -\int_t^{T}\xi_u\cdot d \wh \sS_u
= A^{-1}_T H- \wh V_t-\int_t^{T}\xi_u\cdot d \wh \sS_u.
$$ 
By Theorem \ref{T:5.2}, $\EE^{\QQ^*} [ H^2 ] <\infty$.  
 Let $\wh H_t$ be the rcll version of the square-integrable $\QQ^*$-martingale $\EE^{\QQ^*}[ A_T^{-1} H \mid \sF^{a, \prime}_t]$. 
Then by \eqref{e:5.7}
$$ 
 A^{-1}_T H= \wh H_T=\wh H_t+\int_t^{T}\xi_u^{H}\cdot d\wh \sS_u +(K^{H}_{T}-K^{H}_t)
$$
and so 
  \[
 \wh  C_{T}- \wh C_t=(\wh H_t- \wh V_t)+\int_t^{T}(\xi^{H}_u-\xi_u)\cdot d \wh \sS_u+(K^{H}_{T}-K^{H}_t).
\]
The first term on the right hand-side is $\sF^{a, \prime}_t$-measurable; the last two terms are $\QQ^*$-martingale increments
with $[\wh\sS, K^{H}]\equiv 0$. Hence
\begin{align*}
R^*_t&= \EE^{\QQ^*} [(\wh C_{T}-\wh C_t)^{2} \mid \sF^{a, \prime}_t ] \\
&= (\wh H_t-\wh V_t)^{2}+\EE^{\QQ^*} \Bigl[ \Big( \int_t^{T}(\xi^{H}_u -\xi_u) \cdot d \wh \sS_u \Big)^2 \big| \sF^{a, \prime}_t\Bigr] \\
&\quad +\EE^{\QQ^*}\bigl[(K^{H}_{T}-K^{H}_t)^{2} \mid \sF^{a, \prime}_t\bigr],
\end{align*}
minimized exactly when $\wh V_t=\wh H_t$ and $\int_t^{T}(\xi^{H}_u -\xi_u) \cdot d \wh \sS_u =0$.
As by \eqref{e:n4.12}, 
$$
d \langle  \wh \sS^{(i)}, \wh \sS^{(j)}\rangle _t = \wh \sS^{(i)}_t \wh \sS^{(j)}_t \sum_{k=1}^d \sigma_ik (t) \sigma_{jk} (t) dL_{(t-a)^+},  
$$
it follows from  \eqref{e:con2} that $\int_t^{T}(\xi^{H}_u -\xi_u) \cdot d \wh \sS_u =0$ if and only if 
$ \int_t^{T} |\xi^{H}_u -\xi_u|^2 d L_{(u -a)^+} =0$, that is, if and only if  
$\xi^{H}_u =\xi_u$ $d L_{(\cdot -a)^+} \times \PP$-a.s. on $[t, T] \times \Omega$. 
The corresponding $\eta$ is given 
$$
\eta^H_u = \wh V_u -  \xi_u \cdot  \wh \sS_u=   \wh H_u -  \xi^H_u \cdot  \wh \sS_u  
= \EE^{\QQ^*} \big[ A^{-1}_T H \mid \sF^{a, \prime}_u \big] - \xi^H_u \cdot  \wh \sS_u 
\quad \hbox{for every } u\in [t, T].
$$
This establishes the theorem.
\end{proof}

\subsection{Pricing with the sub-diffusion equivalent martingale measure}\label{S:5.2}

Let $T>a\geq 0$. 
For a contingent claim $H\in L^p ( \sF^{a, \prime}_T, \PP)$ for some $p>2$, we know by Proposition~\ref{P:5.5} that 
it has the optimal replicating strategy $(\xi^H, \eta^H)$ given by   \eqref{e:optim} with 
the discounted value  process 
$$
  \wh  V_t = \EE^{\QQ^*} \big[ A^{-1}_T H  \mid \sF^{a, \prime}_t \big]  \quad \hbox{for } t\in [0, T].
$$
We price  this claim    at time $t\in [0, T]$ by the value process 
\begin{equation}\label{e:Vta}
V_t =A_t\wh V_t = A_t \,  \EE^{\QQ^*} \big[ A^{-1}_T H  \mid \sF^{a, \prime}_t \big]  .
\end{equation} 
In particular, we   price this claim at initial time  by
\begin{equation}\label{e:V0}
V_0=\EE^{\QQ^*} \left[ A^{-1}_T H \right].
\end{equation} 
  As we noted in \S \ref{S:5.1}, this pricing is locally risk-minimizing under the real-world probability $\PP$  
in the sense of \cite[Definition on p.349]{Sch}. 

 \smallskip
 
Since the cumulative cost process $\wh C^H$ is a square-integrable $\PP$-martingale,  we have
$\EE [ \wh C_t ] = \wh C_0$ for all $t\in [0, T]$. Hence a seller who charges $\wh C_0$ and follows
the optimal replicating strategy $(\xi^H, \eta^H)$ breaks even on average at the maturity time $T$ under the real-world probability measure $\PP$. 
By Proposition \ref{P:5.5}, the  F\"ollmer-Schweizer decomposition \eqref{e:5.7} of $A_T^{-1}H$ and \eqref{e:5.11},
 $\EE^{\QQ^*} [  A^{-1}_T  H] + \int_{0}^{T}\xi^{H}_t\cdot d \wh \sS_t$ is the attainable part of the claim $A^{-1}_T  H$,
 while $K^H_T= \wh C_T-\wh C_0$ is the unattainable part of the claim as $[ K^H,  \wh \sS] \equiv 0$. 
 As shown in Proposition \ref{P:5.5}, this pricing is the best possible in the sub-diffusive market setting
 as it resulted from the optimal replicating strategy that minimizes the local risk.

 In subsection \ref{S:5.4}, we will derive an explicit formula for  European   option $\psi (\sS_T)$ at maturity  $T$
 when $d=1$. 
 
 \subsection{Incomplete information} \label{S:5.3}

Let $T>a\geq 0$. 
In the enlarged filtration $\{\wt \sF^a_t\}_{t\in [0, T]}$, 
by the same proof as  the first part of that  for Theorem \ref{T:5.2}, for every $G \in L^p(\wt \sF^a_T , \PP)$ for some $p>2$,
 $\EE^{\QQ^*} [ |G|^q ] <\infty$ for every $q\in (2, p)$.
 Thus by Theorem \ref{T:MRT} and Remark \ref{R:3.5} (applied to the $\PP$-martingale $NM$, where $N_t:=\EE^{\QQ^*}[G\mid\wt\sF^a_t]$, together with the It\^o formula and \eqref{e:5.3}, which turn $dB_{L_{(\cdot-a)^+}}$ into $d\widetilde B_{L_{(\cdot-a)^+}}$), 
$G$  is attainable in the sense that 
there exists an $\RR^d$-valued $\{\wt \sF^a_s \}_{0\leq s\leq T}$-predictable process   $\{ \wt \xi_s;  s\in [0, T]\}$   
 having  $\EE^{\QQ^*}  \int_0^T |  \wt \xi_s  |^2 dL_{(s-a)^+}  <\infty$
so that  
$$
G  =  \EE^{\QQ^*} \big[ G \big| \wt \sF^a_0\big]  +  \int_0^{T}  \wt  \xi_s  \cdot  d\wt B_{L_{(s-a)^+} }  .
$$
In view of \eqref{e:con1}-\eqref{e:con2} and \eqref{e:4.3}, 
for every contingent claim $H\in L^p(\wt \sF^a_T )$ for some $p>2$, 
there is an $\RR^d$-valued $\{\wt \sF^a_s \}_{0\leq s\leq T}$-predictable process  
 $\{ \xi_s;  s\in [0, T]\}$   
 having  $\EE^{\QQ^*}  \int_0^T |  \xi_s  |^2 d \langle  \wh \sS\rangle _s  <\infty$
so that  
\begin{equation}\label{e:5.16a}
A_T^{-1}H  =  \EE^{\QQ^*} \big[ A_T^{-1}H  \big| \wt \sF^a_0\big]  +  \int_0^{T}    \xi_s  \cdot  d \wh \sS_s .
\end{equation} 
Thus with respect to the enlarged filtration $\{\wt \sF^a_t\}_{t\in [0, T]}$, 
every $H \in   L^p(\wt \sF^a_T , \PP)$ for some $p>2$ can be replicated by a self-finance trading strategy $\xi$ 
and initial endowment $\EE^{\QQ^*} \big[ G \big| \wt \sF^a_0\big] $. 
Note that typically $\EE^{\QQ^*} \big[ G \big| \wt \sF^a_0\big] $ is not deterministic. 
 In the real world, the available information is only the smaller filtration $\{\sF^{a, \prime}_t\}_{t\in [0, T]}$.

Fix some $q\in (2, p)$.  By the same proof as that for \eqref{e:Lq}, $H\in L^q(\wt \sF^a_T, \QQ^*)$. 
Consequently,  $\EE^{\QQ^*} \big[ A_T^{-1}H  \big| \wt \sF^a_0\big] \in L^q(\wt \sF^a_0, \QQ^*)$
and so $ \int_0^{T}    \xi_s  \cdot  d \wh \sS_s =A_T^{-1}H - \EE^{\QQ^*} \big[ A_T^{-1}H  \big| \wt \sF^a_0\big] 
  \in L^q(\wt \sF^a_T, \QQ^*)$.
By  Burkholder-Davis-Gundy inequality and the Doob's maximal inequality, 
$ \big( \int_0^T | \xi_s|^2 d\langle  \wh \sS\rangle _s  \big)^{1/2} \in L^q(\wt \sF^a_T, \QQ^*)$. 
Let $r\in (2, q)$. Then by \eqref{e:5.3} and the same argument as that for \eqref{e:5.6n}, we have
$$
\EE^{\QQ^*} \big[ A_T^{-1}H  \big| \wt \sF^a_0\big] \in L^r(\wt \sF^a_0, \PP)
\quad \hbox{and} \quad 
\Big( \int_0^T | \xi_s|^2 d\langle  \wh \sS\rangle _s  \Big)^{1/2} \in L^r(\wt \sF^a_T, \PP).
$$
By the same argument as that for Remark \ref{R:5.1}, $\xi$ satisfies the property \eqref{e:5.1a} but with the filtration
$\{\wt \sF^a_t\}_{t\in [0, T]}$ in place of the filtration $\{  \sF^{a, \prime}_t\}_{t\in [0, T]}$. 
Denote by $\wt {\mathcal P}$  be the space of all $\{\wt \sF^a_t\}_{t\in [0, T]}$-predictable processes $\xi$ having
$\EE \int_0^T |\xi_u|^2 d\langle \wh \sS\rangle _u <\infty$, equipped with the Hilbert norm $\| \xi\|:=\big( \EE \int_0^T |\xi_u|^2 d\langle \wh \sS\rangle _u\big)^{1/2}$;
and $ {\mathcal P}$  be the space of all $\{ \sF^{a, \prime}_t\}_{t\in [0, T]}$-predictable processes $\xi$ having
$\EE \int_0^T |\xi_u|^2 d\langle \wh \sS\rangle _u <\infty $, equipped with the Hilbert norm $\| \xi\|:=\big( \EE \int_0^T |\xi_u|^2 d\langle \wh \sS\rangle _u\big)^{1/2}$.
Clearly, $ {\mathcal P}$ is a closed subspace of  $\wt {\mathcal P}$.  For $\xi \in \wt {\mathcal P}$, denote its orthogonal projection to $ {\mathcal P}$ by
${\rm Proj}_{\mathcal P} (\xi)$. 
With these preparation, we have the following by \cite[Theorem 4.6]{FS}.

\begin{theorem}\label{T:5.7}
Let  $H\in L^p(\sF^{a, \prime}_T )$ for some $p>2$.
The discounted claim $A^{-1}_T H$ has F\"ollmer-Schweizer   decomposition 
$$
A^{-1}_TH= \EE^{\QQ^*} [A^{-1}_T H] + \int_0^T \xi^H_u \cdot d\wh \sS_u  +K^H_T,
$$
where $\xi^H ={\rm Proj}_{\mathcal P} (\xi)$ with $\xi$ being the $\{\wt \sF^a_t\}_{t\in [0, T]}$-predictable process in  \eqref{e:5.16a}
and

\begin{align*}
K^H_t&:=\EE^{\QQ^*}\big[K^H_T\mid\sF^{a,\prime}_t\big],\\
K^H_T&:=\EE^{\QQ^*} \big[ A_T^{-1}H  \big| \wt \sF^a_0\big]  -\EE^{\QQ^*} [A^{-1}_T H]  + \int_0^T (\xi_u-\xi^H_u)\cdot d\wh \sS_u.
\end{align*}
\end{theorem} 

Thus we can find the optimal replicating strategy $( \xi^H, \eta^H)$ of \eqref{e:optim}  from the $\{\wt \sF^a_t\}_{t\in [0, T]}$-predictable process 
$\xi$ in the martingale representation \eqref{e:5.16a} for $A^{-1}_T H$ under $\QQ^*$. 

\subsection{European  option pricing} \label{S:5.4}
 
 In this section, we consider the case $d=1$, that is, 
 there  a risk-free bond (or money market fund)
   and one common stock modeled by \eqref{BSsub}. 
   Let $\psi$ be a Borel measurable function $\RR$  that satisfies the polynomial growth condition, that is, 
   there is some $n\geq 1$ so that $|\psi (x)| \leq c(1+|x|^n)$ on $\RR$.
   Note that by Theorem \ref{T:2.1}  and \eqref{e:4.5}, $\psi (\sS_{T}) \in L^p (\sF^{a, \prime}_T, \PP)$ for every $p>2$.

   In view of Proposition \ref{P:5.5} and \eqref{e:Vta}-\eqref{e:V0} in \S \ref{S:5.2}, we price
   the European   option    that has payoff $\psi (\sS_T)$ at maturity  $T$ by 
   the value process for the optimal replicating strategy, that is, by 
\begin{eqnarray}\label{e:Vt}
 V_t = A_t\wh V_t   = \EE^{\QQ^*}  \left[    \exp{  \Big(    -\int_t^T r_s ds  \Big) }  \psi  (\sS_T ) \Big| \sF^{a,\prime}_t \right]   ,
  \quad t\in [0, T], 
  \end{eqnarray}   
  where $\QQ^*$ is the sub-diffusion   equivalent martingale measure defined by \eqref{e:4.9}-\eqref{e:4.10}, which is also the minimal
  martingale measure.

  We point out that formula \eqref{e:Vt}  is in exact analogy with the classical spot market model driven by Brownian motion. 
  In that setting, the equivalent martingale measure   is unique and thus the market is complete. 
  Under the equivalent martingale measure, 
    the discounted stock price process $A_t^{-1}\sS_t$ is a a local martingale  that is the stochastic integral with respect to 
  the Brownian motion $B$, and the value process $\{V_t; t\in [0, T]\}$
  of the European option
   that has value $\psi (\sS_T)$ at maturity  $T$ is given by \eqref{e:Vt}; see, e.g.,  \cite{HP, KS}. 
For the sub-diffusion  spot market model, 
by Theorem \ref{T:4.1}, $\QQ^*$ defined by \eqref{e:4.9}-\eqref{e:4.10} is the unique equivalent martingale measure 
  under which  the discounted stock price process $A_t^{-1}\sS_t$ is a stochastic integral of
   the sub-diffusion $\widetilde B_{L_{(t-a)^+}}$. 
  
 We next proceed to compute the value process $V_t$ in \eqref{e:Vt} more explicitly.  
 Recall that    $R_t:=R_0+ S_{L_{(t-R_0)^+}} - t $  is the  overshoot process defined in Theorem \ref{T:2.3}. 
 For emphasis, we denote $R_t$ by $R^a_t$ when $R_0=a\geq 0 $,  which is the initial wake-up time 
 or holding time for the market.  
   For notational simplicity, in the rest of this paper we assume that  $\si_t $,  $r_t $ and $\bar \mu_t$ are constants independent of $t$. 

We first compute the condition expectation with the
enlarged Markovian filtration $\{ \sF^a_t \}_{t\geq 0}$ containing the overshoot process  $R$, that is, we first compute
\begin{equation}\label{e:5.17n} 
\bar V_t:= A_t \EE^{\QQ^*} \big[ A^{-1}_T \psi  (\sS_T ) \Big| \sF^{a }_t \big] 
= \EE^{\QQ^*}  \left[    \exp{  \Big(    -\int_t^T r_s ds  \Big) }  \psi  (\sS_T ) \Big| \sF^{a}_t \right]   .
\end{equation} 
Clearly, 
\begin{equation}\label{e:5.18n}
V_t = \EE^{\QQ^*} \big[ \bar V_t \mid \sF^{a, \prime}_t  \big] \quad \hbox{for } t\in [0, T].
\end{equation} 
Note that $\sS_t$ and $L_{(t-a)^+} = \langle  X\rangle _t$ are $\sF^{a, \prime}_t$-measurable.
The overshoot  $R^a_t := a +S_{L_{(t-a)^+}}-t$ is $\sF^a_t$-measurable but in  general is not $\sF^{a, \prime}_t$-measurable. 
 Once the formula for $\bar V_t$ is derived, we can compute  $V_t = \EE^{\QQ^*} \big[ \bar V_t \mid \sF^{a, \prime}_t  \big] $
 by using the conditional distribution of the overshoot  $R^a_t$ given $L_{(t-a)^+}$ 
 (see Theorem \ref{thm:triple}).

\begin{theorem}\label{T:5.1} 
Let  $\bar V_t$ be an  $\sF^{a}_t$-measurable random variable defined by \eqref{e:5.17n}.
 Then  conditioning on   $(\sS_t, R_t)  = (x, \tilde a)$,
\begin{equation} \label{e:5.8a}
\bar V_t= \exp{  \left(    -r (T-t)  \right)}  \int_0^\infty 
  \EE \left[  \psi  \big (x \exp\big(   r(T-t) -\si^2 y/2  + \si \sqrt{ y}Z    \big)   \big)  \right] 
   \PP\left(  L_{(T-t-\tilde a)^+}     \in d y \right)    ,
\end{equation} 
where $Z$ is an independent random variable that is of standard Normal distribution with zero mean and unit variance.
 In particular, when $t=0$, 
 \begin{equation}\label{e:5.9}
 V_0= \bar V_0=\exp{  \left(    -r T \right)}  \int_0^\infty 
  \EE \Big [ \psi  \big(  \sS_0 \exp\big(   rT -\si^2 y/2  + \si \sqrt{ y}Z   \big)  \big)  \Big] 
  \PP\left(  L_{T - a }     \in d y \right)   
 \end{equation} 
 gives the price of the European option with payoff $\psi(\sS_T)$ at initial time $t=0$.
\end{theorem}

\begin{proof}      By \eqref{e:4.14n},   
$$ A_t^{-1} \sS_t =  \sS_0  \exp\left(  \sigma  \widetilde B_{L_{(t-a)^+}} - \frac 12   \si^2  L_{(t-a)^+}   \right)
\quad \hbox{for } t\in [0, T], 
$$
and so 
\begin{equation}\label{sol-Q-S}
\sS_t  = \sS_0 \exp\left(   rt   +  \si  \widetilde B_{L_{(t-a)^+}} - \frac 12    \si^2   L_{(t-a)^+}    \right).
\end{equation}
 We have from \eqref{sol-Q-S} that for any $t\in [0, T)$, 
\begin{eqnarray}\label{sol-Q-S-1}
\sS_T= \sS_t\exp\left(   r(T-t)+  \si   \left(\widetilde B_{L_{(T-a)^+}}-  \widetilde B_{L_{(t-a)^+}}     \right)- \frac 1 2 \si ^2    \left(L_{(T-a)^+}   -    L_{(t-a)^+} \right)   \right).
\end{eqnarray}
 It is shown in \cite{ZC}   (see two lines above (3.3) there) that 
$$ 
L_{(T-a)^+} - L_{(t-a)^+}  =L_{(T-t-R_t)^+} \circ \theta_{L_{(t-a)^+} }, 
$$
where $R_t= a +S_{L_{(t-a)^+}}-t$.
Thus by  \cite[Theorem 3.1]{ZC} (note that $\widetilde B_{L_{(\cdot-a)^+}}$ has the same law under $\QQ^*$ as $B_{L_{(\cdot-a)^+}}$ under $\PP$),  conditioning on   $(\sS_t, R_t)  = (x, \tilde a)$,  
\begin{eqnarray}\label{V-Q-1}
    \bar V_t   
  &=&  \exp{  \left(    -r (T-t)   \right) }  \EE^{{\QQ^*}}  \left[    \psi   (\sS_T)   \mid  \sF^{a}_t  \right]  \nonumber\\
  &=&  \exp{  \left(    -r (T-t)\right)} \EE \bigg[ \psi  \big(x \exp\big(   r (T-t)+  \si  B_{L_{(T-t -\tilde a)^+}}  - \tfrac 1 2 \si ^2   L_{(T-t -\tilde a)^+}  \big) \big)  \bigg]  \\
 &=& \exp{  \left(    -r (T-t)  \right)}  \int_0^\infty 
  \EE \left[  \psi  \big (x \exp\big(   r(T-t) -\si^2 y/2  + \si \sqrt{ y}Z    \big)   \big)  \right] 
   \PP\left(  L_{(T-t-\tilde a)^+}     \in d y \right)    \nonumber  \\  
 &=& \exp{  \left(    -r (T-t) \right)}  \int_0^\infty 
  \EE \bigg[ \psi  \big( x \exp\big(   r(T-t) -\si^2 y/2  + \si \sqrt{ y}Z  \big)  \big)  \bigg] 
  d_y  \PP\left(  S_y \geq (T-t-\tilde a)^+ \right)   ,   \nonumber     
       \end{eqnarray}   
 where the last equality is due to \cite[Lemma 2.1]{C1}. This establishes the desired formula \eqref{e:5.8a}.
\end{proof}
 
 \medskip

  \begin{example} \rm
Two-sided estimates on the density of a class of   subordinators  are derived  in \cite[Theorem 4.4]{CKKW2}.
When $S$ is a $\beta$-stable subordinator, its density and the density of its inverse can be found in 
\cite[Propositions 3.1 and 3.2]{GK}.  In particular,  for each $t>0$, 
$  \PP( L_t     \in d x ) = h_\beta(t,x)dx $ on $[0, \infty)$, where   
$$
 h_\beta(t, x) = \frac1{\pi} \sum_{k=0}^\infty \frac{(-x)^k}{k!} \Gamma (\beta (k+1)) t^{-\beta (k+1)}  \sin (\beta (k+1)\pi),
 \quad x>0.
 $$
 This can be used in \eqref{e:5.8a} and \eqref{e:5.9}  to ``explicitly" compute the value process $V_t$ and 
 the current  price $V_0$. 
  \end{example}
  
  \medskip
  
  \begin{remark} \rm 
 When the subordinator $S$ has unit drift and null L\'evy measure (that is, $S_t=t$ for all $t\geq 0$) and $a=0$, $L_{(t-a)^+}=t$ and the sub-diffusion $B_{L_{(t-a)^+}}=B_t$ 
 becomes the standard Brownian motion.  In this case, \eqref{e:5.9} gives 
 \begin{equation}\label{e:5.16}
 V_0=\exp{  \left(    -r T \right)} 
  \EE \Big [ \psi  \big(\sS_0 \exp\big(  (r-\si^2 /2)T   + \si \sqrt{ T}Z  \big)  \big)  \Big] ,
 \end{equation} 
 In particular, when $\psi  (x)=(x-\mathsf K)^+$, 
 \begin{equation}\label{e:5.17}
 V_0=\exp{  \left(    -r T \right)} 
  \EE \bigg[\Big(\sS_0 \exp\big(  (r-\si^2 /2)T   + \si \sqrt{ T}Z  \big)  -\mathsf K\Big) ^+ \bigg] ,
 \end{equation} 
 which is the price of the European call option in the classical Black-Scholes model where 
 $$
 d{\sS_t}= {\sS_t} \left( \mu  dt +\sigma  d B_t \right)
 $$
 and $dA_t=r A_t\,dt$. 
 Thus Theorem \ref{T:5.1} covers the classical Black-Scholes model as a special case.
 Hence we can call  \eqref{e:5.8a} and \eqref{e:5.9} the Black-Scholes formula for European call options
 in our  sub-diffusion  spot market model.  \qed
\end{remark}

 \subsection{Time-fractional Black-Scholes  PDE}\label{S:5.5}

In this section, we derive a time-fractional Black-Scholes  PDE for the extended value process $\bar V$ defined by
\eqref{e:5.17n}. 
Recall that $\kappa \geq 0$ and $\nu$ are the drift and L\'evy measure of the subordinator $S$. 
For $x>0$, define $w(x):= \nu ( (x, \infty))$.
Following \cite{C1, C2}, we define the time-fractional derivative of  a function $f$ on $(0, \infty)$ by
$$
\partial^w_t f(t):= \frac{d}{dt} \int_0^t w(s) (f(t-s)-f(0))ds , t > 0,
$$
whenever the derivative exists. 

Consider the following time-fractional partial differential equation (PDE)
 \begin{eqnarray} \label{e:5.18a}
\begin{cases}
 (\kappa \frac{ \partial}{\partial t} + \partial^w_t) u(t, x) =  \frac{ 1}2 \sigma^2 x^2 \  \frac{\partial^2}{\partial x^2}   u(t, x)
&\hbox{for } x>0,  t>0  \smallskip \\
u(0, x) =\psi  (x)   &\hbox{for } x>0   .
\end{cases}
\end{eqnarray}
Note that   the geometric Brownian motion 
 \begin{equation} \label{e:geoBM}
  Y_t:=Y_0 \exp\big(    \si  B_t  - \tfrac 1 2 \si ^2  t \big), \quad t\geq 0,  
 \end{equation}
 satisfies SDE $dY_t = \sigma Y_t dB_t$. Thus $Y_t= Y_0 \exp\big(    \si  B_t  - \tfrac 1 2 \si ^2  t \big) $ is a conservative continuous Feller process on $(0, \infty)$ with generator  $\frac{ 1}2 \sigma^2 x^2\frac{d^2}{dx^2}$.   
 Denote by $(\sL, {\rm Dom}(\sL))$ the Feller generator of $Y$ in the space of continuous functions $C_\infty ((0, \infty))$ that vanishes at infinity
 equipped with the uniform norm $\| \cdot \|_\infty$. 
  When $\psi  \in {\rm Dom}(\sL)$, it is first shown in \cite[Theorem 2.1]{C1} and then further strengthened in
 \cite[Theorem 3.1]{C2} that $u(t, x):= \EE_x [ \psi  (Y_{L_t}) ]$ is the unique strong solution to \eqref{e:5.18a}; see  the statement
 of \cite[Theorem 3.1]{C2} for its precise definition.  However, for European call option, we need to take $\psi  (x)=(x-\mathsf K)^+$
 which is not in ${\rm Dom}(\sL)$ so the aforementioned result can not be applied directly. 
 To circumvent this difficulty, we first write the time-fractional PDE in its integral form
 \begin{equation}\label{e:integral} 
\kappa (u(t, x)-\psi  (x)) + \int_0^t w(s) (u(t-s, x)-\psi  (x))ds  =  \frac{ 1}2 \sigma^2 x^2    \frac{\partial^2}{\partial x^2}  \int_0^t  u(s, x) ds .
 \end{equation}
 For $\lambda >0$, define the Laplace transform of $t\mapsto u(t, x)$ by $\bar u(\lambda, x)$; that is, 
 $\bar u(\lambda, x):= \int_0^\infty e^{-\lambda t} u(t, x) dt$. Observe that by the Tonelli theorem, 
 $$
 \int_0^\infty e^{-\lambda t} \left( \int_0^t  u(s, x) ds \right) dt 
 =\int_0^\infty u(s, x) \left( \int_s^\infty e^{-\lambda t} dt \right) ds  = \frac{\bar u(\lambda, x) }{\lambda}, 
 $$
 and   
 $$
 \int_0^\infty e^{-\lambda s} w(s) ds = \int_0^\infty 
 \left(\int_0^\infty e^{-\lambda s} \1_{\{z\geq s\}} \nu (dz) \right) ds
 = \frac1{\lambda} \int_0^\infty  (1-e^{-\lambda z}) \nu (dz) = \frac{\phi (\lambda) - \kappa \lambda}{\lambda},
 $$
 where $\phi (\lambda)$ is the L\'evy exponent \eqref{e:2.2} of the subordinator $S$ and $\kappa\geq 0$ is its drift.  
   Taking   Laplace transform in $t$ on both sides of \eqref{e:integral} yields that for $\lambda > 0$, 
 $$
\frac{\phi (\lambda)}{\lambda} \Big(\bar u(\lambda, x) - \frac{\psi  (x)}{\lambda}\Big)
 =  \frac{ 1}2 \sigma^2 x^2    \frac{\partial^2}{\partial x^2}  \frac{\bar u(\lambda, x)}{\lambda}   , 
 $$
that is,
 \begin{equation}\label{e:Laplace} 
    \frac{ 1}2 \sigma^2 x^2    \frac{\partial^2}{\partial x^2}   \bar u(\lambda, x) 
    -    \phi (\lambda)  \bar u(\lambda, x)  + \frac{\phi (\lambda)}{\lambda} \psi  (x) =0
    \qquad \hbox{for } x>0.
 \end{equation}
 Since   a function on $[0, \infty)$  is uniquely characterized by its Laplace transform,  this leads us to the following definition.

 \begin{definition} \rm A function $u(t, x)$ is said to be a solution to the time-fractional PDE \eqref{e:5.18a}  
 if there is some $\lambda_0\geq 0$
 so that for any $\lambda>\lambda_0 $, the Laplace transform $\bar u(\lambda, x)$ of $u(t, x)$ in $t$-variable exists for all $x>0$ and 
 $x\mapsto \bar u(\lambda, x)$ has at most linear growth and is a weak solution of the PDE  \eqref{e:Laplace} for every $\lambda >\lambda_0$. 
 \end{definition} 
 
 \begin{theorem} \label{T:5.5}
 Let $\mathsf K>0$. There is a unique solution to the following  time-fractional PDE
    \begin{eqnarray} \label{e:5.21}
\begin{cases}
 (\kappa \frac{ \partial}{\partial t} + \partial^w_t) u(t, x) =  \frac{ 1}2 \sigma^2 x^2 \  \frac{\partial^2}{\partial x^2}   u(t, x)
\quad &\hbox{for } x>0 \hbox{ and }  t>0,    \smallskip \\
u(0, x) = (x-\mathsf K)^+     & \hbox{for } x>0 . 
\end{cases}
\end{eqnarray}
This unique solution is given by $u(t, x)=\EE_x [  (Y_{L_t} -\mathsf K)^+]$. 
 \end{theorem}
 
 \begin{proof} Let $\{\psi_n; n\geq 1\}$ be a sequence of $C^2$-smooth functions on $(0, \infty)$ with compact support so that 
 $\psi_n =0$ on $(0, \mathsf K/2)$, $0\leq\psi_n(x)\leq 1+x$ and $\lim_{n\to\infty}\psi_n (x)=(x-\mathsf K)^+$ for $x\in(0 , 3\mathsf K/2)$, $\psi_n (x) = x-\mathsf K$ on $[3\mathsf K/2, 4\mathsf K+n]$  
 and $\lim_{n\to\infty}\psi_n(x)=x-\mathsf K$, increasingly in $n$, for $x>2\mathsf K$.
 Since $\psi_n \in {\rm Dom} (\sL)$, the time-fractional PDE \eqref{e:5.18a} with $\psi_n$ in place of $\psi$ 
 has a unique strong solution  
 \begin{equation}\label{e:un}
 u_n (t, x):= \EE_x [ \psi_n (Y_{L_t})] , \quad x>0
 \end{equation} 
  in the Banach space $(C_\infty ((0, \infty)), \| \cdot \|_\infty)$ in the sense of \cite[Theorem 3.1]{C2}.
 By the argument from \eqref{e:integral} to \eqref{e:Laplace}, for every $\lambda >0$ its  Laplace transform $\bar u (\lambda, x)$
 is a weak solution to 
  \begin{equation}\label{e:Laplace2} 
    \frac{ 1}2 \sigma^2 x^2    \frac{\partial^2}{\partial x^2}   \bar u_n (\lambda, x) 
    -    \phi (\lambda)  \bar u_n (\lambda, x)  = -  \frac{\phi (\lambda)}{\lambda} \psi_n  (x) , \quad x> 0 .
 \end{equation}
    that is, $u_n$ is a solution to \eqref{e:5.18a} with $\psi_n$ in place of $\psi$.
 We can solve \eqref{e:Laplace2}  by a change of variable. Define $\wt u_n (\lambda, x)= \bar u_n (\lambda, e^x)$ for $x\in \RR$ and $\lambda>0$. 
 Then 
\begin{equation}\label{e:chv}
 \frac{\partial }{\partial x} \wt u_n (\lambda, x) = e^x   \frac{\partial }{\partial x} \bar u_n (\lambda, \cdot) (e^x)
 \quad \hbox{and} \quad  
 \frac{\partial^2 }{\partial x^2} \wt u_n  (\lambda, x) = (e^x)^2  \frac{\partial^2 }{\partial x^2} \bar u_n (\lambda, \cdot) (e^x) + e^x 
 \frac{\partial }{\partial x} \bar u_n (\lambda, \cdot) (e^x) .
 \end{equation}
 Hence \eqref{e:Laplace2} is equivalent to $\wt u(\lambda, x)$ satisfying
  \begin{equation}\label{e:Laplace3} 
     \left(   \frac{ 1}2 \frac{\partial^2}{\partial x^2}    
    -   \frac{ 1}2 \frac{\partial }{\partial x }  -    \frac{\phi (\lambda)}{\sigma^2}   \right) \wt  u_n (\lambda, x) 
       =-  \frac{\phi (\lambda) \psi_n  (e^x) }{\lambda\sigma^2}     \quad \hbox{on } \RR.  
 \end{equation}
  Observe that for each $\lambda >0$, $\frac{ 1}2 \frac{\partial^2}{\partial x^2}    
    -   \frac{ 1}2 \frac{\partial }{\partial x }  -    \frac{\phi (\lambda)}{\sigma^2}   $
  is the infinitesimal generator of $W_t-(t/2)$ 
    killed at rate  $\frac{\phi (\lambda)}{\sigma^2} $,
  where $W_t$ is a one-dimensional Brownian motion. 
  Hence
  \begin{equation}\label{e:5.25}
  \wt u_n (\lambda, x)= \frac{\phi (\lambda)}{\lambda\sigma^2}\int_0^\infty  e^{- (\phi (\lambda)/\sigma^2)t} \int_{\RR} \frac{1}{\sqrt{2\pi t}}
   \exp \left( - \frac{ (y-x +(t/2) )^2}{2t} \right)  \psi_n  (e^y)   dy\, dt.
  \end{equation}
  
  As $\psi_n (x) \1_{(0, 2\mathsf K]} (x) $ converges to $(x-\mathsf K)^+  \1_{(0, 2\mathsf K]} (x) $ boundedly and 
  $\psi_n (x) \1_{( 2\mathsf K, \infty)} (x) $ increases to $(x-\mathsf K)^+  \1_{( 2\mathsf K, \infty) } (x) $ as $n\to \infty$,
  we have from  \eqref{e:un} and the  bounded and monotone convergence theorem  that 
  \begin{equation} \label{e:5.26}
  \lim_{n\to \infty} u_n(t, x)= \EE_x [  (Y_{L_t} -\mathsf K)^+] =: u(t, x), \quad x>0.
  \end{equation}
  Consequently, $\lim_{n\to \infty} \bar u_n(\lambda , x) = \bar u (\lambda, x)$, where 
  $\bar u (\lambda , x):= \int_0^\infty e^{-\lambda t} u(t, x) dt$.
  Set $\wt u(\lambda, x):= \bar u (\lambda , e^x)$.  Clearly we have
   $\lim_{n\to \infty} \wt  u_n(\lambda , x) = \wt  u (\lambda, x)$,
  On the other hand, taking $n\to \infty$ in \eqref{e:5.25}, we have
  \begin{align}
  \wt u  (\lambda, x) &=  \frac{\phi (\lambda)   }{\lambda\sigma^2}   \int_0^\infty  e^{- (\phi (\lambda)/\sigma^2)t} \int_{\RR} \frac{1}{\sqrt{2\pi t}}
   \exp \left( - \frac{ z^2}{2t} \right)    (e^{x+z-(t/2)}-\mathsf K)^+  dz\,  dt \label{e:5.27}\\
   &\leq  \frac{\phi (\lambda)   }{\lambda\sigma^2}   \int_0^\infty  e^{- (\phi (\lambda)/\sigma^2)t}\,  e^{x}\,   dt  = \frac{e^x}{\lambda}. \label{e:5.28}
  \end{align}
  We used $\EE[e^{W_t-t/2}]=1$ for a standard Brownian motion $W$.
By \eqref{e:5.28}, for every $\lambda >0$, $x\mapsto \wt u(\lambda, x)$ is finite and continuous on $\RR$.
Moreover, $\bar u(\lambda,x)\leq x/\lambda$, so that $x\mapsto \bar u (\lambda, x)$ has linear growth.
Since $\EE_x[Y_{L_t}]=x$, we also have $0\leq u(t,x)\leq x$ for $x>0$ and $t>0$.
By \eqref{e:5.27}, for each $\lambda >0$, 
  $x \mapsto \wt  u(\lambda , x)$ is a weak solution to 
   \begin{equation}\label{e:Laplace4} 
     \left(   \frac{ 1}2 \frac{\partial^2}{\partial x^2}    
    -   \frac{ 1}2 \frac{\partial }{\partial x }  -    \frac{\phi (\lambda)}{\sigma^2}   \right) \wt  u  (\lambda, x) 
       =-  \frac{\phi (\lambda)   (e^x-\mathsf K)^+ }{\lambda\sigma^2}     \quad \hbox{on } \RR.  
 \end{equation}
 A change of variable  as in \eqref{e:chv} shows that $\bar u(\lambda , x)$ is a weak solution to 
    \begin{equation}\label{e:Laplace5} 
    \frac{ 1}2 \sigma^2 x^2    \frac{\partial^2}{\partial x^2}   \bar u  (\lambda, x) 
    -    \phi (\lambda)  \bar u  (\lambda, x)  = -  \frac{\phi (\lambda)}{\lambda} (x-\mathsf K)^+ , \quad x> 0 .
 \end{equation}
 This proves that $u(t, x)=\EE_x [  (Y_{L_t} -\mathsf K)^+]$ is a solution to the time-fractional PDE \eqref{e:5.21}.
 
 \smallskip
 
 The uniqueness part is easier. Suppose $v(t, x)$ is another solution of \eqref{e:5.21}. Then $\Psi(t, x):=u(t, x)-v(t, x)$
 is a solution to 
    \begin{eqnarray*}  
\begin{cases}
 (\kappa \frac{ \partial}{\partial t} + \partial^w_t) \Psi(t, x) =  \frac{ 1}2 \sigma^2 x^2 \  \frac{\partial^2}{\partial x^2}   \Psi(t, x)
\quad &\hbox{for } x>0 \hbox{ and } t>0,    \smallskip \\
\Psi(0, x) =0     & \hbox{for } x>0 . 
\end{cases}
\end{eqnarray*}
By definition, there is some $\lambda_0 >0 $ so that for any $\lambda >\lambda_0$,
  $\bar \Psi(\lambda, x):=\int_0^\infty e^{-\lambda t} \Psi(t, x) dt$ is a weak solution for 
  $$
     \frac{ 1}2 \sigma^2 x^2    \frac{\partial^2}{\partial x^2}   \bar \Psi(\lambda, x) 
    -    \phi (\lambda)  \bar \Psi(\lambda, x)   =0
    \qquad \hbox{for } x>0
$$
and $| \bar \Psi (\lambda, x)| \leq c_\lambda (|x|+1)$ on $(0, \infty)$ for some $c_\lambda >0$. 
Then  $\wt \Psi(\lambda, x):= \bar \Psi (\lambda, e^x)$ is a weak solution for 
 $$
     \left(   \frac{ 1}2 \frac{\partial^2}{\partial x^2}    
    -   \frac{ 1}2 \frac{\partial }{\partial x }  -    \frac{\phi (\lambda)}{\sigma^2}   \right) \wt  \Psi  (\lambda, x) 
       =0    \quad \hbox{on } \RR
 $$
 with $|\wt \Psi(\lambda, x)| \leq c_\lambda (1+e^x)$ on $\RR$. 
 In the class of functions with $|\wt \Psi(\lambda, x)|\leq c_\lambda (1+e^x)$, this equation has only the zero solution, and so $\wt \Psi (\lambda, x)\equiv \bar \Psi(\lambda , x) \equiv 0$.  Consequently, $\Psi (t, x)\equiv 0$. 
 This establishes the uniqueness.
\end{proof}

  \medskip

 \begin{theorem} \label{T:5.6}
 The  process $\bar V_t$  of \eqref{e:5.17n} 
 is given by 
 $$
 \bar V_t= e^{-r(T-t)} u(T-t, e^{r(T-t)}\sS_t, R^a_t)  \quad \hbox{for } t\in [0, T], 
 $$
 where  for each $a\in [0, T)$,
  $(t, x)\mapsto u(t, x, a)$ is the unique solution
 to the following time-fractional PDE:
  \begin{eqnarray} \label{e:fBS}
\begin{cases}
 (\kappa \frac{ \partial}{\partial t} + \partial^w_t) u(t, x, a) =  \frac{ 1}2 \sigma^2 x^2 \  \frac{\partial^2}{\partial x^2}   u(t, x, a)
&\hbox{for } t>a,    \smallskip \\
u(t, x, a) = (x-\mathsf K)^+    &\hbox{for } t\in [0, a].
\end{cases}
\end{eqnarray}
 \end{theorem}
 
 \begin{proof} 
 Note that by \eqref{V-Q-1},  conditioning on   $(\sS_t, R^a_t)  = (\tilde x, \tilde a)$,  
  \begin{eqnarray*}
 \bar V_t &=&   \exp{  \left(    -r (T-t)\right)} \EE \bigg[\bigg(\tilde x \exp\big(   r (T-t)+  \si  B_{L_{(T-t -\tilde a)^+}}  - \tfrac 1 2 \si ^2   L_{(T-t -\tilde a)^+}  \big) -\mathsf K\bigg) ^+ \bigg] \\
  &=& \exp{  \left(    -r (T-t)\right)} \EE \bigg[\bigg(   e^{r(T-t)} \tilde x \exp\big(    \si  B_{L_{(T-t -\tilde a)^+}}  - \tfrac 1 2 \si ^2   L_{(T-t -\tilde a)^+}  \big) -\mathsf K\bigg) ^+ \bigg] 
 \end{eqnarray*}
 Define 
 \begin{equation}\label{e:5.33} 
 u(t, x, a):= \EE \bigg[\bigg(x \exp\big(    \si  B_{L_{(t-  a)^+}}  - \tfrac 1 2 \si ^2   L_{(t -  a)^+}  \big) -\mathsf K\bigg) ^+ \bigg]  
 \end{equation}
  Then by \eqref{V-Q-1}, 
 \begin{equation} \label{e:5.18}
 \bar V_t = e^{-r(T-t)} u(T-t, e^{r(T-t)}\sS_t, R^a_t).
 \end{equation}

 We can rewrite $u(t, x, a)$ as 
 $$
 u(t, x, a)= \EE \left[ (Y_{L_{(t-a)^+}}  -\mathsf K)^+|Y_0=x \right],
 $$
 where $Y$ is the geometric Brownian motion defined by \eqref{e:geoBM}. 
 By Theorem \ref{T:5.5}, for each $a \geq 0$, $(t, x)\mapsto u(t, x, a)$ is the unique solution for
  \begin{eqnarray} \label{e:5.19}
\begin{cases}
 (\kappa \frac{ \partial}{\partial t} + \partial^w_t) u(t, x, a) =  \frac{ 1}2 \sigma^2 x^2 \  \frac{\partial^2}{\partial x^2}   u(t, x, a)
&\hbox{for } t>a,    \smallskip \\
u(t, x, a) = (x-\mathsf K)^+    &\hbox{for } t\in [0, a].
\end{cases}
\end{eqnarray}
This completes the proof of the theorem.
\end{proof} 
 
 \medskip

 \begin{remark} \label{R:5.7}  \rm 
  When the subordinator $S_t\equiv t$ (that is, when $\kappa =1$ and $\nu=0$) and $a=0$, 
 $(t, x)\mapsto u(t, x, 0)$ satisfies 
 \begin{eqnarray} \label{e:5.20}
\begin{cases}
  \frac{ \partial}{\partial t}   u(t, x, 0) =  \frac{ 1}2 \sigma^2 x^2 \  \frac{\partial^2}{\partial x^2}   u(t, x, 0)
&\hbox{for } t>0,   \smallskip \\
u(0, x, 0) = (x-\mathsf K)^+     . 
\end{cases}
\end{eqnarray}
In this case,  $\{\sF^{0, \prime}_t\}_{t\geq 0}=  \{\sF^{0}_t\}_{t\geq 0}$ and
   $V_t=\bar V_t= v(t, \sS_t)$, where     
   $$
   v(t, x):= e^{-r(T-t)}u(T-t, e^{r(T-t)}x, 0).
   $$
    Note that  
\begin{eqnarray*}  
\begin{cases}
  \frac{ \partial}{\partial t}   v(t, x) =  e^{-r(T-t)} (ru -    \frac{ \partial}{\partial t} u  - r e^{r(T-t)}  x \partial_x u )  (T-t, e^{r(T-t)}x, 0)     \smallskip \\
 \partial_x   v(t, x) =   ( \frac{ \partial}{\partial x} u ) (T-t, e^{r(T-t)}x, 0)  \smallskip \\
    \partial^2_x   v(t, x)  =  e^{r(T-t)}    ( \frac{ \partial^2 }{\partial x^2} u ) (T-t, e^{r(T-t)}x, 0)  . 
\end{cases}
\end{eqnarray*}
Thus   it follows from \eqref{e:5.20} that $v(t, x)$ satisfies  
\begin{eqnarray*}  
\begin{cases}
\frac{ \partial}{\partial t}   v(t, x)  +\frac12  \sigma^2 x^2    \frac{ \partial^2 }{\partial x^2} v (t, x) + r  x   \frac{ \partial}{\partial x }  v (t, x) - r v (t, x)   =0
&\hbox{for } (t, x)\in (0, T) \times (0,\infty), \smallskip  \\
v(T, x)= (x-\mathsf K)^+. 
\end{cases}
\end{eqnarray*}
 This is the classical Black-Scholes  PDE for the European call option where the stock price is modeled by SDE  \eqref{e:1.2} 
 driven by Brownian motion.
 For this reason,  we  call \eqref{e:5.19} the time-fractional Black-Scholes  PDE for the price process $\bar V$ of the European stock call option.
\qed
\end{remark}

 \subsection{Computing the value process}\label{S:5.6}
 
 In this subsection, we indicate how to compute $V_t$ from $\bar V_t$ by \eqref{e:5.18n}:
 $$
 V_t = \EE^{\QQ^*} \big[ \bar V_t \mid \sF^{a, \prime}_t ] =  e^{-r(T-t)}  \EE^{\QQ^*} \big[ u(T-t, e^{r(T-t)}\sS_t, R^a_t)  \mid \sF^{a, \prime}_t  ] .
 $$
 As mentioned earlier, given $\sF^{a, \prime}_t $,  we know $\sS_t$ and $L_{(t-a)^+}$. 
 Note that the distribution of the subordinator $S$ and hence its inverse is invariant under $\QQ^*$ and $\PP$.  
 So 
 \begin{equation}\label{e:5.45} 
 V_t = 
 \begin{cases}
    e^{-r(T-t)}   u(T-t, e^{rT}\sS_0, a-t ) \quad &\hbox{for } t\in [0, a], \\
 e^{-r(T-t)} \int_{[0, \infty)} u(T-t, e^{r(T-t)}\sS_t, \wt a) \,  \PP \big(R^a_t \in d \wt a \mid  L_{(t-a)^+} \big)
 \quad &\hbox{for } t\in (a, T].
 \end{cases}
 \end{equation} 
 
 To find the formula for the conditional law in \eqref{e:5.45}, we assume in the following that 
 $\phi_0 (\xi) :=\int_{(0, \infty)} (1-\cos (\xi z) ) \nu (dz)$ satisfies the 
 Hartman--Wintner  condition
\begin{equation} \label{e:HW} 
\lim_{|\xi|\to\infty}\frac{ \phi_0 (\xi)}{\log(1+|\xi|)}=\infty ,
\end{equation}
This mild condition ensures that the subordinator $S$ has a jointly continuous density function $\bar p (s, x)$ with respect to the Lebesgue measure, that is, $\PP(S_s\in dx)=\bar p(s,x)\,dx$.
 
 The conditional law of the overshoot process $R^0_t$ given $L_t$ under condition \eqref{e:HW} 
 is   known to experts in L\'evy process
 so the right hand side of \eqref{e:5.45} can be computed out `explicitly'. 
 However, there is no ready reference we can cite directly. 
 For the reader's convenience, we give a proof. Recall $\kappa$ and $\nu$ are the drift and the L\'evy measure 
 of the subordinator $S$, and $\phi$ is its Laplace exponent. 
 Recall that $w(x)= \nu ( (x, \infty))$ for $x>0$. Let 
 $$
 \mu_r (dx):= \PP ( S_r \in dx) \quad \hbox{and} \quad  U(dx) = \int_0^\infty \mu_r (dx) dr.
 $$
 Since our subordinator is not a compounded Poisson, the potential measure 
$U$ is diffusive in the sense that $U(\{x\})=0$ for every $x\in [0, \infty)$;
see \cite[Proposition I.15 on p.30]{Ber1}.  
It is known that $U$ has finite measure on compact subsets of $[0, \infty)$.
 In fact, by \cite[Proposition III.1 on p.74]{Ber1}, there are positive constants $c_1$ and $c_2$ so that 
 \begin{align} \label{e:BeIII1}
\frac{c_1}{ \phi (1/t) }  \leq U ([0, t] ) \leq  \frac{ c_2}{\phi (1/t)}  \quad  \text{ for every  } t>0.
\end{align}
When $\kappa =0$, it is known (see \cite[Theorem III.4 on p.77]{Ber1}) 
that 
\begin{equation} \label{e:5.48}
 \PP (S_{L_t} >t)= 1 \quad \hbox{for every } t>0.
\end{equation} 
 When the subordinator $S$ has positive drift $\kappa > 0$, by
a result due to by a result due to  Neveu (see \cite[Proposition 1.7]{Ber2}),
  $U(dx)$ is absolutely continuous with respect to the Lebesgue measure on $[0, \infty)$,
  has a strictly positive continuous density function ${\vartheta} (x)$ on $[0, \infty)$ with ${\vartheta} (0)=1/\kappa$,
  and
\begin{equation}\label{e:3.1a}
  {\vartheta}(x) =  \kappa^{-1}  \, \PP (\hbox{there is some $t\geq 0$ so that } S_t=x)  \quad \hbox{for }  x\geq 0.
 \end{equation}
 Furthermore,  one has  (see, e.g., \cite[Theorem III.5 on p.79]{Ber1}) ,
\begin{equation}\label{e:3.1}
  \PP (S_{L_x} =x) = \kappa \,  {\vartheta}(x)  \quad \hbox{for  every }  x > 0.
  \end{equation}
    
   For notational convenience, let $O_t:=R^0_t=S_{L_t}-t$. 
 
\begin{theorem} \label{thm:triple}
\begin{enumerate}
\item[\rm (i)] (Overshoot by a jump) 
For $s>0$, $0\le y\le t$ and $r>0$,
\begin{equation}\label{e:5.49} 
\PP\bigl(L_t\in ds,\ S_{L_t-}\in dy,\ O_t\in dr\bigr)=ds\,\mu_s(dy)\,\nu (t-y+dr).
\end{equation} 

\item[\rm (ii)] ($O_t=0$ case)  Suppose the Hartman--Wintner condition \eqref{e:HW} holds.
When $\kappa >0$,  
\begin{equation}\label{e:5.51} 
\PP\bigl(L_t\in ds,\ O_t=0\bigr)=\kappa \, \bar p (s, t)\,ds \quad \hbox{for Lebesgue-almost every} t>0. 
\end{equation} 

\item[\rm (iii)] Suppose the Hartman--Wintner condition \eqref{e:HW} holds.
For Lebesgue-almost every $t>0$, the inverse subordinator $L_t$ has a density function with respect to the Lebesgue 
measure on $[0, \infty)$:
\begin{equation}\label{e:5.52} 
\rho_{L_t} (s)  =\int_{[0,t]} w(t-y) \mu_s(dy)\ +\kappa\bar p (s, t),    \quad s>0,
\end{equation} 
Moreover, for a.e. $s>0$, 
\begin{equation}\label{e:5.53} 
\PP\bigl(O_t\in dr\mid L_t=s\bigr)
=\frac{\displaystyle \1_{\{r>0\}}\int_{[0,t]}\mu_s(dy)\, \nu (t-y+dr)\;+\;\kappa\, \bar p (s, t)\,\delta_{\{0\}} (dr)}
{\displaystyle \int_{[0,t]}  w (t-y) \mu_s(dy)   +\kappa \bar p( s, t)} .
\end{equation} 
\end{enumerate} 
\end{theorem}

\begin{proof}
(i) Since $S$ is non-decreasing, $S_u\le t$ for $u<L_t$ so $S_{L_t-}\le t$. Thus $O_t>0$ holds if and only if there is a (necessarily unique) time $s$ with $S_{s-}\le t<S_{s-}+\Delta S_s$, and in that case $L_t=s$. Therefore for any bounded continuous function $f$ on $\RR_+^3$,
\[
f(L_t,S_{L_t-},O_t)\1_{\{O_t>0\}}=\sum_{s>0} f\bigl(s,S_{s-},S_{s-}+\Delta S_s-t\bigr)\,\1_{\{S_{s-}\le t<S_{s-}+\Delta S_s\}} .
\]
The jumps $(s,\Delta S_s)$ form a Poisson point process with intensity $ds\otimes\nu(dx)$.
So \[
\EE[ f (L_t,S_{L_t-},O_t)\1_{\{O_t>0\}} ]
 =\EE\int_0^\infty ds\int_{(0,\infty)}\ \nu (dx)\,f(s,S_s,S_s+x-t)\,\1_{\{S_s\le t<S_s+x\}}.
\]
 This establishes \eqref{e:5.49}. 
 
 \smallskip
 
 (ii) When $\kappa >0$, we have $dL_t = \kappa^{-1} \1_{\{O_t =0\}} dt $ by \cite[Proposition 3.2]{ZC}.
 Hence by a change of variable, we have for any continuous functions $f, g$ with compact support on  $(0, \infty)$
 \[
\int_0^\infty g(t)f(L_t)\1_{\{O_t=0\}}\,dt =\kappa  \int_0^\infty g(S_{L_t})f(L_t) \,dL_t 
=\kappa\int_0^\infty g(S_u)f(u)\,du .
\]
Taking expectation on both sides yields 
$$
 \int_0^\infty g(t) \EE \big[ f(L_t)\1_{\{O_t=0\}} \big]\,dt = \kappa \int_0^\infty \int_0^\infty g(t) \bar p(u, t) f(u) du dt 
 $$
 Hence, for Lebesgue-almost every $t>0$,
 $$  \EE \big[ f(L_t)\1_{\{O_t=0\}} \big] = \kappa \int_0^\infty  \bar p(u, t) f(u) du \quad \hbox{for almost every} t>0.
 $$
 As the above holds for every continuous function $f$  with compact support in $(0, \infty)$, we get \eqref{e:5.51}.

 \smallskip
 
 (iii) \eqref{e:5.52} and \eqref{e:5.53} follow readily from \eqref{e:5.48}, \eqref{e:5.49} and \eqref{e:5.51}
\end{proof} 
 
 \medskip
 
 As $R^a_t = S_{L_{(t-a)^+}} +a -t$, one readily derives the conditional law
 $\PP\bigl(R^a_t \in dr\mid L_{(t-a)^+}=s\bigr)$ from that of $\PP\bigl(O_t \in dr\mid L_t=s\bigr)$ in \eqref{e:5.53}.
 Indeed, for $t >a$, as $R^a_t=O_{t-a}$, 
 $$
 \PP\bigl(R^a_t \in dr\mid L_{(t-a)^+}=s\bigr) 
 = \PP\bigl(O_{t-a} \in dr  \mid L_{t-a}=s\bigr).
 $$
  This together with \eqref{e:5.45} gives, for almost every $t\in(a,T]$, the `explicit' formula for the pricing (or cost) process $V_t$ of \eqref{e:5.45}
 for the European call option.

\end{document}